\documentclass[mnsc,nonblindrev,copyedit]{informs3}

\usepackage{endnotes}
\usepackage{enumerate}
\usepackage{subfigure}
\usepackage{amsmath}
\usepackage{bm}
\usepackage{mathtools}
\usepackage[square,sort,comma,numbers]{natbib}
\usepackage{yhmath}
\usepackage[normalem]{ulem}
\usepackage{verbatim}

\usepackage{algorithm}
\usepackage{algorithmic}

\usepackage{setspace,graphicx,booktabs,color,url}
\definecolor{DarkBlue}{rgb}{0,0.08,0.45}
 \bibpunct[, ]{(}{)}{,}{a}{}{,}%
 \def\bibfont{\small}%
 \def\newblock{\ }%
 \def\BI\!Band{and}%

 \newcommand{\dist}{\text{dist}}

 \newcommand{\bean}{\begin{eqnarray*}}
 \newcommand{\eean}{\end{eqnarray*}}
 \newcommand{\bea}{\begin{eqnarray}}
 \newcommand{\eea}{\end{eqnarray}}
 \newcommand\defeq{\;\stackrel{\mathclap{\text{def}}}{=}\;}
 \usepackage{algorithm} 
 \usepackage{algorithmic}

\theoremstyle{EX}

\newtheorem{heuristic}{Heuristic}
\TheoremsNumberedThrough     
\ECRepeatTheorems

\EquationsNumberedThrough    

\MANUSCRIPTNO{} 

\begin{document}


\RUNAUTHOR{}

\RUNTITLE{}


\TITLE{\large{Dynamic Type Matching}}


\ARTICLEAUTHORS{
\AUTHOR{Ming Hu\\
Rotman School of Management, University of Toronto, Toronto, Ontario, Canada M5S 3E6\\
\href{mailto:ming.hu@rotman.utoronto.ca}{ming.hu@rotman.utoronto.ca}\\
\bigskip
Yun Zhou\\
DeGroote School of Business, McMaster University, Hamilton, Ontario, Canada L8S 4L8\\
\href{mailto:zhouy185@mcmaster.ca}{zhouy185@mcmaster.ca}\\
}
\AFF{Oct 21, 2018}
}

\ABSTRACT{
We consider an intermediary's problem of dynamically matching demand and supply of heterogeneous types in a periodic-review fashion. More specifically, there are two disjoint sets of demand and supply types, and a reward associated with each possible matching of a demand type and a supply type. In each period, demand and supply of various types arrive in random quantities. The platform's problem is to decide on the optimal matching policy to maximize the total discounted rewards minus costs, given that unmatched demand and supply will incur waiting or holding costs, and will be carried over to the next period (with abandonment). For this dynamic matching problem, we provide sufficient conditions on matching rewards such that the optimal matching policy follows a priority hierarchy among possible matching pairs. We show those conditions are satisfied by vertically and unidirectionally horizontally differentiated types, for which quality and distance determine priority, respectively. As a result of the priority property, the optimal matching policy boils down to a match-down-to threshold structure when considering a specific pair of demand and supply types in the priority hierarchy.
}
\maketitle

\linespread{1.5}

\section{Introduction}

Operations management is about managing the process of matching supply with demand. We consider a firm that periodically manages the matching between demand and supply.
In each period, demand and supply of various types arrive in \emph{random} quantities.
Each ``type'' represents a distinct set of characteristics of demand or supply.
The matching between demand and supply generates type- and time-dependent reward. With unmatched demand and supply fully or partially rolled over to the next period, the firm aims to maximize the total expected rewards (minus costs of waiting compensation for demand and inventory holding for supply).

The problem we describe above is crucial to many intermediaries who \emph{centrally} manage matchings in a sharing economy. 
Sharing economy platforms often use crowdsourced supply and match it dynamically with customer demand.
For example, 
commuter carpooling platforms such as UberCommute match a driver heading to a destination with a rider to the same destination (or in the same direction).
Amazon crowdsources inventories of an identical item from third-party merchants to its warehouses, to fulfill online orders.
The nonprofit organization, United Network for Organ Sharing (UNOS), allocates donated organs to patients in need of transplantation.
In the center of those business and nonprofit sharing-economy models, a platform is developed and maintained by an intermediary to enable sharing-economy activities. Those models have the following features.

\emph{Heterogeneous demand and supply types.} From the intermediary firm's perspective, matching between demand and supply of different characteristics often generates distinct rewards (or equivalently, mismatch costs).
We refer to demand/supply of different characteristics as different \emph{types} of demand/supply,
and consider two possible ways in which demand/supply types differ from each other.
In particular, types can be horizontally or vertically differentiated.
{\label{marker:hor-ver-exp} Horizontal differentiation means that the characteristics of a type are not always superior or inferior to those of another type (regarding generating matching rewards).
Instead, the matching reward between a demand type and a supply type is determined by the two's idiosyncratic taste on each other.
For example, for a ride-hailing platform, riders and drivers are characterized by their locations, with the matching between a pair closer to each other generating a higher reward (i.e., a shorter waiting time for the rider and shorter idle time for the driver).
Vertical differentiation means quality differences in the demand/supply types.
Under vertical differentiation, a particular type is always superior or inferior to a different type regarding generating matching rewards.
For example, from the perspective of UNOS, patients and organs may differ in their health condition. A patient/donor in a better health condition, in general, leads to a better transplant outcome. }

\emph{Time-variant uncertainty on both sides of the market.}
In contrast to conventional business models where supply is often treated as a decision (e.g., inventory replenishment decision) or a fixed capacity (e.g., in revenue management problems), crowdsourced supply in sharing economy activities may arrive at the system randomly and dynamically.
For example, in ride-hailing activities, drivers decide, on their own, when and how much time they make themselves available to provide service.
In Amazon's inventory commingling program, third-party merchants use their own inventory-regulating policies and may be subject to various time-varying supply shocks.

We use a finite-horizon stochastic dynamic program to formulate the problem with the features mentioned above.
Next, we present an overview of the main results of the paper, as well as the applications and implications of the model and the results.

\subsection*{Main results, applications and implications}

{A key result of the paper is the establishment of the modified Monge conditions. 
    \label{marker:key-result}
Under those conditions, a particular pair of demand and supply types should have ``priority" over a neighboring pair (i.e., a pair sharing the same demand or supply type) in the optimal matching policy.
This allows us to simplify the matching decision within a period, and focus on the trade-off between matching in the current period and that in the future. }

{\label{marker:trade-off0}
The optimal matching policy is complicated even for the static problem.
For example, consider a specific period $t$ without accounting for future arrivals of demand and supply.
On the one hand, one may want to prioritize the matching between a type $i$ demand and a type $j$ supply if the unit matching reward $r_{ij}^t$ is high.
On the other hand, matching $i$ with $j$ may prevent both matching $i$ with another type $j'$ supply and matching another type $i'$ demand with $j$. 
If $r_{ij}^t< r_{ij'}^t + r_{i'j}^t$, it may be undesirable to prioritize matching $i$ with $j$.
Moreover, the optimal matching policy is further complicated by possibly saving a demand or supply type for the current period and matching it with future supply or demand.  In other words, there are trade-offs \emph{within} a period, as well as \emph{across} the current period and future periods.}

{
Under the conditions we establish in this paper, we are able to prioritize the matching within a period.
This allows us to focus on the trade-off between the current period and future periods. 
}



Then we study two special versions of the model, namely, the model with horizontally differentiated types (in short, the horizontal model) and the model with vertically differentiated types (in short, the vertical model). Both satisfy the established modified Monge condition.

{\bf
The horizontal model.} We consider demand and supply types distributed in a metric space $C$ (which can be considered as the space of characteristics of demand/supply).  The matching reward between a demand type and a supply type depends on the ``distance'' between the two. The shorter the distance, the higher the reward. We start by studying the case with two demand types and two supply types. {In that case, a perfect pair (i.e., type 1 demand with type 1 supply, or type 2 demand with type 2 supply, both associated with the highest unit matching rewards) should be prioritized and matched greedily, whereas an imperfect pair (i.e., type 1 demand with type 2 supply, or type 2 demand with type 1 supply, which has a lower unit matching reward compared with a perfect pair) should be considered only when the corresponding demand and supply types have sufficiently high levels (after the greedy matching of the perfect pairs) and matched down to some threshold level. \label{marker:perfect-pair}}(We will define perfect and imperfect pairs formally in Section \ref{sec:struc-properties}.)
{\label{marker:trade-off1} Therefore, the main trade-off is between a lower reward from matching an imperfect pair in the current period and a possible higher reward by reserving demand/supply to form perfect pairs in a future period.}

When there are multiple demand and supply types, we focus on the \emph{unidirectional} case in which $C$ is a directed line segment, and the supply travels along a given direction to reach the demand for the matching. 
We show that a shorter distance implies a higher priority, i.e., the optimal policy would assign a demand type to the closest available supply type. 

The horizontal model has the following applications.

\emph{Capacity management with upgrading.} Upgrading uses a high-class supply to fulfill a  low-class demand, which is widely adopted in travel industries (see, e.g.,  \citealt{yu2014dynamic}) and in production/inventory settings (see, e.g., \citealt{bassok1999single}). 
\citet{shumsky2009dynamic} study a revenue management problem with fixed initial capacities of various supply types, and demand types can only be upgraded \emph{one-level up}.
\citet{yu2014dynamic} study the general upgrading problem, allowing demand types upgradable to be matched with a generally higher-quality supply type. The upgrading reward structure in \cite{yu2014dynamic} is a special case of unidirectionally horizontal types located along a line.  Thus our results apply to a generalized capacity management problem with general upgrading and random replenishment. The feature of random supply is desirable for upgrading, even for those revenue management settings, not to mention for the production/inventory settings.  For example, in car rental, car returns can be random, and in airline ticket selling,  early cancellations or airplane swaps can result in random capacity changes.


\emph{Commuter carpooling along a fixed route.} 
Carpooling platforms specifically designed for commuters,  such as UberCommute and GrabHitch, match riders heading to the same destination (or in the same direction). In those cases, the matching reward has two additive components: The first one is a disutility associated
with the distance traveled along the fixed route from the driver's current location to pick up the demand. 
The second is a utility associated with traveling along the route from the demand's pick-up location to its drop-off location. The former is the unidirectionally horizontal case, whereas the latter is a vertically differentiated attribute because, given the same pick-up location, it is more desirable if the demand's travel
distance is longer. We show that if riders and drivers head to the same destination at the end of the route, 
a shorter distance to pick up a rider
on the way has a higher priority in matching. 


{\bf The vertical model.} Each demand and supply type is associated with a quality level, with higher quality types leading to higher matching rewards. 
In particular, we focus on the case where the reward of matching a pair is the sum of the contributions brought in by its components, which are increasing in quality. 
 Then the optimal matching policy follows a simple structure, which we call \emph{top-down matching} (in an economic term, assortative mating): line up demand and supply types in descending order of their ``quality" levels from high to low; match them from the top, down to some level. Thus, the optimal matching policy in any period can be entirely determined by a total matching quantity. This result is generalizable to the case where the matching reward is non-linear but not far from being additive.

{\label{marker:trade-off2}
In the vertical model, 
the main trade-off is again between the current period and future periods.
The optimal policy will reserve some (lower-quality) demand or supply type(s), to reduce the chance of losing or delaying the matching of potential high-quality types arriving in the future. 
}

For two special cases, namely, the case with patient demand and supply and the case with impatient demand and patient supply, we further derive monotonicity properties of the optimal total matching quantity with respect to the state of demand and supply. 
We also propose a one-step-ahead (OSA) heuristic policy, which is guaranteed to perform better than greedy matching, and significantly reduce the degree of state-dependency.

The vertical model may shed light to the following applications. 

\emph{Online dating.} 
In the settings of assortative mating such as online dating platforms, the participants of matching have vertically distributed attributes such as wealth and education. 
\citet[Chapter 4, p. 31, Eq. (4.2)]{becker2009social} assume that in a decentralized marriage market the output of a marriage is the sum of the marital incomes of male and female.
We consider the same reward structure, but with dynamic and random arrivals of males and females and from a centralized perspective. The top-down matching structure in our vertical model implies that a centralized dating agency (or even a decentralized dating platform) may want to limit the number of matching pairs at any time, in anticipation of future arrivals of higher-quality participants. 

\emph{Organ allocation.}
Organ allocation decisions involve many factors, such as the efficient use of organs and health conditions of the patients.
On the one hand, organs differ in their quality (which can be determined by risk factors such as age and cause of the donor's death).
Higher quality of the organ in general leads to better post-transplantation health outcomes.
On the other hand,  patients differ in their health condition. Those who are sicker suffer lower quality of life and greater risk of death, and thus receive a higher benefit from transplantation.
The top-down matching procedure in our vertical model suggests that organs of higher quality levels and patients in worse health conditions should receive higher matching priority, and it can be optimal to reject some low-quality organs for patients in anticipation of high-quality organs arriving in the near future.\footnote{In addition to quality differences, the matching between a patient and an organ is subject to compatibility constraints. The top-down structure sheds light on the matching among patients and organs that are mutually compatible.}



 \section{Literature Review}

\begin{figure}[htb!]
\FIGURE
{\includegraphics*[scale=0.5]{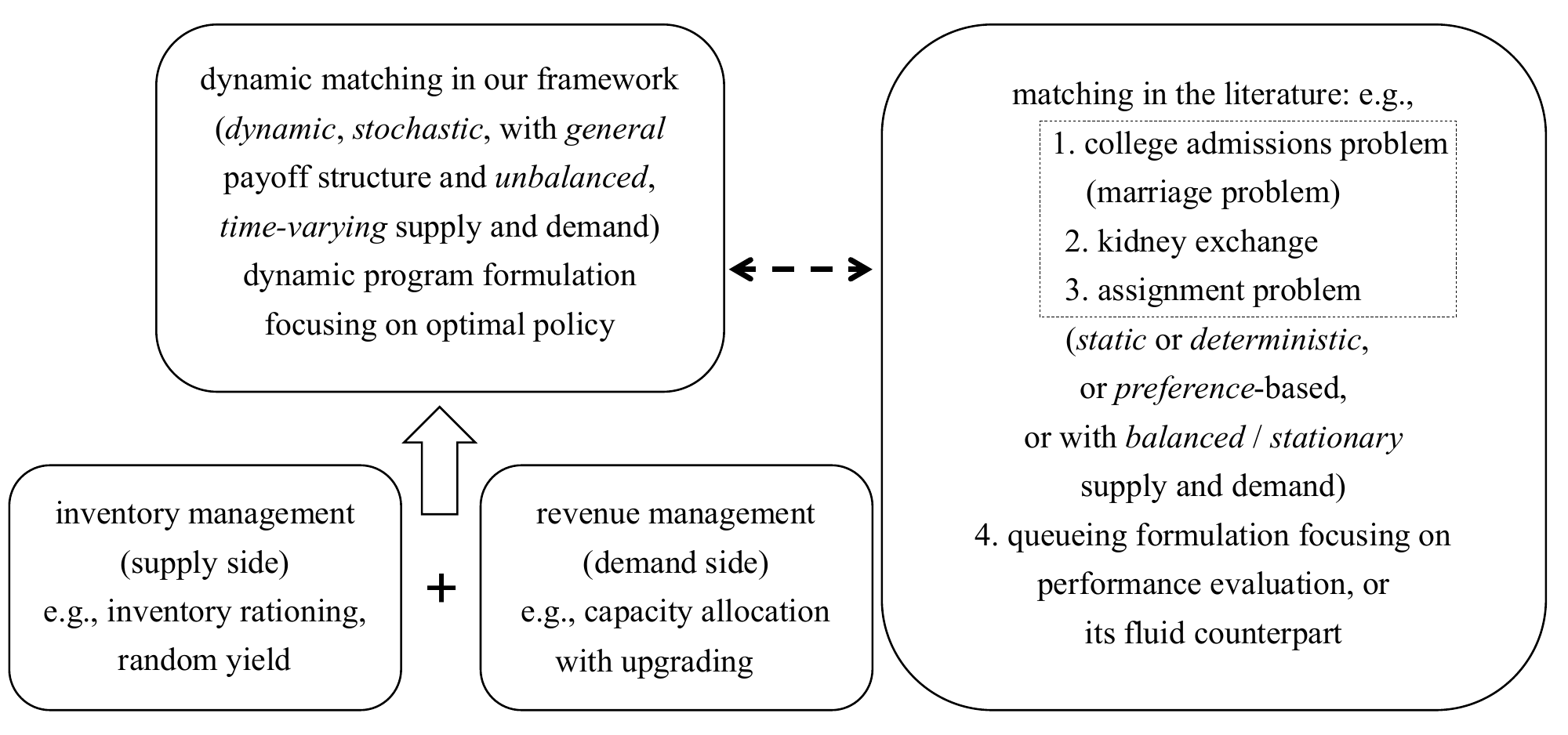}}
{Positioning in the literature.\label{fig:litrev}}
{}
\end{figure}

We illustrate the high-level positioning of our framework in Figure \ref{fig:litrev}. 
The proposed dynamic-matching framework can be viewed as a generalization of two foundations of operations management,  i.e., inventory management where the firm orders the supply centrally (\citealt{zipkin2000foundations}), and revenue management where the firm regulates the demand side with a fixed supply side (\citealt{talluri2006theory}),  and of a combination of the two, i.e., joint pricing and inventory control (\citealt{chen2012pricing}). Compared with existing work in inventory and revenue
management, the supply in the sharing economy is crowdsourced and hence has uncertainty. 

Driven by real-life applications, economists, computer scientists, and operations researchers have studied a variety of two-sided matching problems  (see, e.g., \citealt{RS90}, \citealt{abdulkadiroglu2013matching} for a survey), which include the college admissions problem (with the marriage problem as a special case), kidney exchange and the online bipartite matching problem. We compare our framework with those problems as follows. 

The college admissions problem and the marriage problem are \emph{preference-based}, and focus on finding stable matchings in a static and deterministic setting. In those problems, parties on both the demand and supply sides submit preferences over options (see, e.g., \citealt{ashlagi2014optimal}) to the matching agency.
As the matching outcomes (i.e., college admissions and marriages) can be life-changing events for the participants, serious efforts in soliciting preferences are necessary.
In contrast,  soliciting preferences may not be practical for day-to-day, or even real-time operations in sharing economy activities. For instance, when riders hail a car on Uber, they do not have the option, or may not even bother with which driver to serve them.
To handle such situations, we assign a ``monetary" contribution to the matching between a pair of demand and supply types, instead of adopting preferences by demand and supply.
For example, a lower reward will be generated if a farther-away car is dispatched.

In a typical situation of the kidney exchange,  patients and donors arrive in pairs, with an incompatible patient and donor in each pair.
Subject to compatibility constraints, researchers have designed efficient matching mechanisms based on cycles (e.g., two-way exchanges) or chains of patient-donor pairs (see, e.g., \citealt{RothSomnUnve:2004,RothSomnUnve:2007}) to maximize the number of matchings.
\citet{unver2010dynamic} studies dynamic kidney exchange with inter-temporal random arrivals of patient-donor pairs and attempts to maximize the number of matched compatible pairs.   Our model differs from his by allowing arbitrary \emph{unbalanced} arrivals of demand and supply, and considering the objective to maximize matching total reward minus cost (i.e., social welfare or profit).

Online bipartite matching problems have many applications such as allocation of display advertisements. Initiated by \citet{karp1990optimal}, the classic version considers a bipartite graph $G=(U,V,E)$, and assumes that the vertices in $U$ arrive in an ``online" fashion. That is, only when a vertex $u\in U$ (e.g., a web viewer) arrives, are its incident edges (e.g., his interests) revealed. Then $u$ can be matched to a previously unmatched adjacent vertex in $V$ (e.g., an advertiser). The objective is to maximize the number of matchings.  The problem has many variants, all with the focus on algorithms' competitive ratios (see \citealt{manshadi2012online} for a more recent literature review). The main difference from our model is the ``online" feature, other than that there is no explicit notation of inventory, with one side (e.g., advertisers) always there and the other (e.g., impressions) getting lost if not matched. Instead of worst-case analysis, we focus on the expected value optimization. 

Operations researchers have studied two-sided matching by the queueing approach or its fluid counterpart. 
\citet{arnostimanaging} study a decentralized two-sided matching market and show that limiting the visibility of applicants can significantly improve the social welfare. 
With a fluid approach of modeling stochastic systems, \citet{zenios2000dynamic} and \citet{su2006recipient} study kidney allocation by exploring the efficiency-equity trade-off, and \citet{AkanAlagAtaErenSaid:2012} study liver allocation by exploring the efficiency-urgency trade-off. 
Using double-sided queues, \citet{zenios1999modeling} studies the transplant waiting list and 
\citet{afeche2014double} study trading systems of crossing networks. 
\citet{su2004patient} analyze a queueing model with service discipline FCFS or LCFS to examine the role of patient choices in the kidney transplant waiting system.
\citet{adan2012exact} show that the stationary distribution of FCFS matching rates for two infinite multi-type sequences is of product form. 
\citet{gurvich2014dynamic} study the dynamic control of matching queues with the objective of minimizing holding costs. 
Focusing on the fluid approximation and its asymptotic optimality, the authors observe that in principle, the controller may choose to wait until some ``inventory" of items builds up to facilitate more rewardable matches in the future. We also make a similar observation.  
\cite{KanSab:2018} study a dynamic fluid matching model in which agents on one side receive proposals from those on the other side and determine whether they would pay screening cost to discover the value of the proposing agent. They show that, suitable restriction imposed by the matching platform on the searching of the agents can reduce wasted search effort.
In contrast to the above papers, we focus on the stochastic model (vs. the fluid counterpart) and optimal decision making (vs. performance evaluation).

\section{The Model}

\begin{figure}[htp!]
\begin{center}
    {\includegraphics*[scale=0.2]{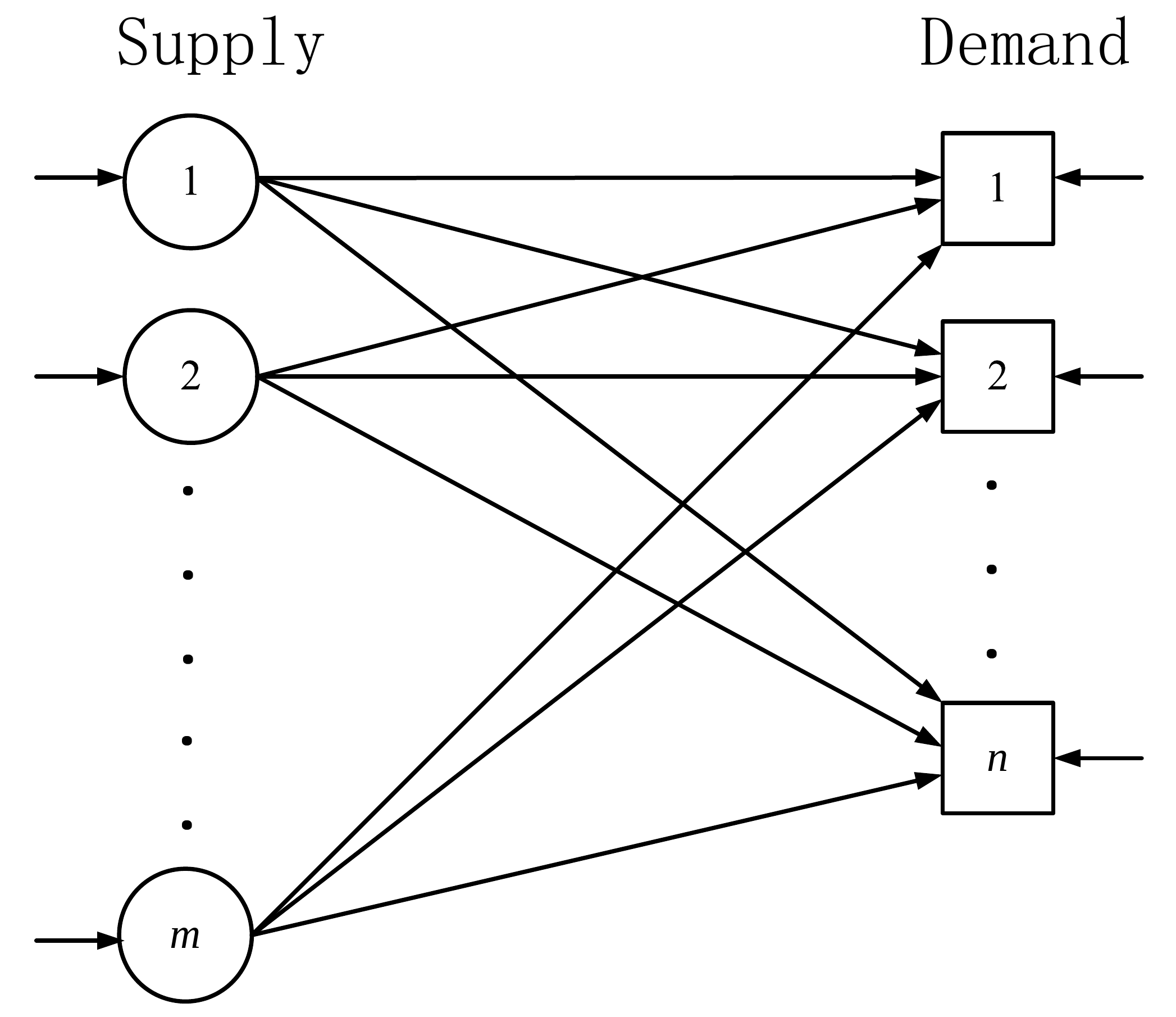}}
    \caption{Pairs of demand and supply.\label{fig:supply-demand}}
\end{center}
\end{figure}

Consider a finite horizon with a total number of $T$ periods. 
At the beginning of each period, $m$ types of demand and $n$ types of supply arrive in \emph{random} quantities.
Let $\mathcal{D}$ be the set of demand types and $\mathcal{S}$ be the set of supply types.
With a slight abuse of notation, we write $\mathcal{D}=\left\{ 1,2,\dots,m \right\}$ and $\mathcal{S}=\left\{ 1,2,\dots,n \right\}$, noting that $\mathcal{D}$ and $\mathcal{S}$ are disjoint sets. 
We use $i$ to index a demand type and $j$ to index a supply type. The pairs of demand and supply are shown in Figure \ref{fig:supply-demand} as a bipartite graph.
An arc $(i,j)$ represents the matching of type $i$ demand and type $j$ supply. 
Without loss of generality, we consider a complete bipartite graph in the base model. In other words, any demand type can potentially be matched with any supply type, apparently with different rewards (or equivalently, mismatch costs). If a demand type is not allowed to pair with supply type $j$, we can just set the matching reward between the two to zero.
We denote the complete set of arcs by $\mathcal{A}= \left\{ (i,j)\mid 1\le i\le m, 1\le j\le n \right\}$.

The state for a given period $t$ comprises the demand and supply levels of various types before matching but after the arrival of random demand $\mathbf{D}^t\in\mathbb{R}_+^{m}$ and supply $\mathbf{S}^t\in\mathbb{R}_+^{n}$ for that period. 
The distributions of supply and demand in one period can be \emph{exogenously} correlated with those in another period. But our model does not account for endogenized correlations among distributions of demand and supply, e.g., a driver's current pickup of a customer may affect future supply at the place where the driver drops off the customer. In other words, we assume away the possible dependence of future distributions of demand and supply on the current matching decisions.   

We denote, as the system state, the demand vector by $\mathbf{x}=(x_1,\dots,x_m)\in\mathbb{R}_+^{m}$ and the supply vector by $\mathbf{y}=(y_1,\dots,y_n)\in\mathbb{R}_+^{n}$,
where $x_i$ and $y_j$ are the quantity of type $i$ demand and type $j$ supply available to be matched. Although we assume that the states and the demand and supply arrivals are continuous quantities (and therefore so are the matching decisions), our results can be readily replicated if those quantities are discrete. 
On observing the state $(\mathbf{x},\mathbf{y})\in\mathbb{R}_+^{m+n}$,  the firm decides on the quantity $q_{ij}$ of type $i$ demand to be matched with type $j$ supply, for any $i\in\mathcal{D}$ and $j\in \mathcal{S}$. For conciseness, we write the decision variables of matching quantities in a matrix form as $\mathbf{Q}=(q_{ij})\in\mathbb{R}_+^{m\times n}$, with $\mathbf{Q}_i$ its $i$-th row (as a row vector) and $\mathbf{Q}^j$ its $j$-th column (as a column vector).  
There is a reward $r_{ij}^t$ for matching one unit of  type $i$ demand and one unit of type $j$ supply for all $i, j$.\footnote{We can account for the case with forbidden arcs. If $(i,j)\notin \mathcal{A}$, we can let $r_{ij}^t$ be zero or a negative number.}
We can write the rewards in a matrix form as $\mathbf{R}^t=(r_{ij}^t)\in\mathbb{R}^{m\times n}$. 
Thus the total matching reward is linear in the matching quantities. That is,  $\mathbf{R}^t\circ \mathbf{Q}\equiv \sum_{i=1}^{m}\sum_{j=1}^{n}r_{ij}^t q_{ij},$ where ``$\circ$" gives the sum of elements of the Hadamard product of two matrices. 
The \emph{post-matching levels} of type $i$ demand and type $j$ supply are given by $u_i=x_i -\bm{1}^m \mathbf{Q}_i^{\tt{T}}= x_i-\sum_{j'=1}^{m}q_{ij'}$ and $v_j=y_j- \bm{1}^n \mathbf{Q}^j=y_j-\sum_{i'=1}^{n}q_{i' j}$, respectively. That is, 
$\mathbf{u}=\mathbf{x}-\bm{1}^m\mathbf{Q}^{\tt{T}}$ and $\mathbf{v}=\mathbf{y}-\bm{1}^n\mathbf{Q}$. The post-matching levels cannot be negative; i.e., $\mathbf{u}\ge\bm{0}$, $\mathbf{v}\ge \bm{0}$.

The unmatched demand and supply at the end of a period carry over to the next period with a fraction of $\alpha$ and $\beta$, respectively. 
In other words, $(1-\alpha)$ fraction of demand and $(1-\beta)$ fraction of supply leave the system. Without loss of generality, we assume they leave the system with zero surpluses. 
{The carry-over fractions $\alpha$ and $\beta$ can be time-dependent (in which case they should be written as $\alpha_t$ and $\beta_t$). But because such time dependency would not affect our results, for ease of notation, we suppress the subscript $t$.}


The firm's goal is to determine a matching policy $\mathbf{Q}^*=(q^*_{ij})$ that maximizes the expected total discounted surplus (i.e., reward minus cost). (Our perspective is social-welfare maximization. Alternatively, the formulation can account for profit maximization if $r_{ij}^t$ is interpreted as the revenue collected from a matching.) 
Let $V_t(\mathbf{x},\mathbf{y})$ be the optimal expected total discounted surplus given that it is in period $t$ and the current state is $(\mathbf{x},\mathbf{y})$.
We formulate the finite-horizon problem by using the following stochastic dynamic program:
\begin{eqnarray}
    V_t(\mathbf{x},\mathbf{y}) & = & \max_{\mathbf{Q}\in \{\mathbf{Q}\ge\bm{0}\mid \mathbf{u}\ge\bm{0}, \mathbf{v}\ge \bm{0}\}} \quad H_t(\mathbf{Q},\mathbf{x},\mathbf{y}), \nonumber \\
    H_t(\mathbf{Q},\mathbf{x},\mathbf{y}) & = & \mathbf{R}^t\circ \mathbf{Q} + EV_{t+1}(\alpha \mathbf{u} + \mathbf{D}^t,\beta \mathbf{v} +\mathbf{S}^t).   \label{eqn:opt}
  \end{eqnarray}
 The boundary conditions are $V_{T+1}(\mathbf{x},\mathbf{y})=0$ for all $(\mathbf{x},\mathbf{y})$, without loss of generality. In other words, at the end of the horizon, all unmatched demand and supply leave the system with zero surpluses. 
 Note that we do not explicitly discount future rewards in (\ref{eqn:opt}) because discouting is implicitly accounted for by using time-dependent rewards.

 A matching policy $P=\left\{ \mathbf{Q}^t(\mathbf{x},\mathbf{y}) \right\}_{t=1,\ldots,T}$ consists of $T$ mappings, $\mathbf{Q}^t:\mathbb{R}^{m+n}\rightarrow \mathbb{R}^{m\times n}$, where $\mathbf{Q}^t(\mathbf{x},\mathbf{y})=(q_{ij}^t(\mathbf{x},\mathbf{y}))_{i=1,\ldots,m, \ j=1\ldots,n}$ is a feasible matching decision in period $t$ for state $(\mathbf{x},\mathbf{y})$.

As mentioned, the state of the system $(\mathbf{x},\mathbf{y})$ is assumed to be real-valued without loss of generality. Nevertheless, our formulation (\ref{eqn:opt}) applies to integer-valued states (with $\alpha$ and $\beta$ equal to either 0 or 1; i.e., each demand/supply type is either completely patient or completely inpatient).

{\label{marker:wait-cost} We can account for waiting costs of those demand and supply types that are not immediately matched by incorporating those costs into the matching rewards. 
 Suppose that a demand type $i$ (resp., supply type $j$) incurs a per-unit waiting cost $c_i^t$ (resp., $h_j^t$) in period $t$ if unmatched. In Online Appendix \ref{sec_appendix:waiting}, we prove that the problem with costs shares the same optimal matching policy as that without costs but with an updated per-unit benefit of matching type $i$ demand with type $j$ supply in period $t$ as $\bar{r}_{ij}^t:=r_{ij}^t+ \sum_{\tau=t}^{T} \alpha^{\tau-t}c_j^\tau+\sum_{\tau=t}^{T} \beta^{\tau-t}h_j^\tau$. To see this intuitively, 
if a unit of type $i$ demand (resp., type $j$ supply) is never matched, its total waiting cost from period $t$ to period $T$ is $ \sum_{\tau=t}^{T}  \alpha^{\tau-t}c_j^\tau$ (resp.,  $\sum_{\tau=t}^{T} \beta^{\tau-t}h_j^\tau$),
which could be saved if this unit for type $i$ demand is matched in period $t$. 
Then, the per-unit benefit of matching type $i$ demand with type $j$ supply in period $t$ becomes $\bar{r}_{ij}^t$, which can be used as the unit matching reward in place of $r_{ij}^t$.}

The existence of an optimal matching policy is resolved by the following proposition.

\begin{proposition}
    \label{lem:concave}
    The functions $H_t(\mathbf{Q},\mathbf{x},\mathbf{y})$ and $V_t(\mathbf{x},\mathbf{y})$ are continuous and concave. There exists an optimal matching policy $P^*=\{ \mathbf{Q}_t^*(\mathbf{x},\mathbf{y})\}_{t=1,\ldots,T}$. 
\end{proposition}

Note that the continuity and concavity in Proposition \ref{lem:concave} hold only when the states and decisions take continuous values.
For the problem with integer-valued states and decisions, concavity is undefined in high dimensional spaces.
Nevertheless, an optimal decision still exists, and all subsequent results still hold.
In general, we expect the optimal policy to be \emph{state-dependent} and extremely complex. In the next section, we characterize some of its structural properties.

\section{Structural Properties of the Optimal Policy}
\label{sec:struc-properties}

We are interested in matching policies with natural properties, e.g., matching an ``essential" pair of a demand type and a supply type before matching any less so pairs.

In particular, we compare two neighboring pairs of demand and supply types (i.e., two arcs in the bipartite network that are incident to a common vertex) and provide sufficient conditions for one pair to be more ``essential'' than another.

We first define a partial relation $\succ_{\mathcal{M}}$ to compare two pairs of demand and supply.

\begin{definition}\label{def:Neighbor-dominance} 
    (Weak Modified Monge Condition)
    We say that $(i,j)\succ_{\mathcal{M}} (i',j)$ if (i) $r_{ij}^{t}\ge r_{i'j}^{t}$; and (ii)  $r_{ij}^t-r_{i'j}^t\ge \alpha (r_{ij''}^{t+1}-r_{i'j''}^{t+1})$ for all $t=1,\ldots,T$ and all $j''\in\mathcal{S}$.
    Similarly, we say that $(i,j)\succ_{\mathcal{M}} (i,j')$ if (i) $r_{ij}^{t}\ge r_{ij'}^{t}$; and (ii)  $r_{ij}^t-r_{ij'}^t\ge \beta (r_{i''j}^{t+1}-r_{i''j'}^{t+1})$ for all $t=1,\ldots,T$ and all $i''\in\mathcal{D}$.
\end{definition}

Let us examine the condition $r_{ij}^t-r_{ij'}^t\ge \beta (r_{i''j}^{t+1}-r_{i''j'}^{t+1})$ in Definition \ref{def:Neighbor-dominance} (to which, the other condition, $r_{ij}^t-r_{i'j}^t\ge \alpha (r_{ij''}^{t+1}-r_{i'j''}^{t+1})$, is symmetric). The condition is easy to satisfy if suppliers are impatient (i.e., 
$\beta$ is small). If suppliers are relatively patient (i.e., $\beta$ is close to 1), for the inequality to hold, the differences in contributions to the matching brought by $j$ and $j'$ should not increase notably over time. 

We further define a stronger partial relation for comparing two pairs of demand and supply.

\begin{definition}\label{def:Neighbor-dominance-strong}
    (Modified Monge Condition)
    We denote the partial relation $\succ_{\mathcal{M}}$ by $\succ_{\mathcal{M}_s}$ if 
    for any $(i,j)\succ_{\mathcal{M}} (i',j)$ and $(i,j)\succ_{\mathcal{M}} (i,j')$, it holds for all $t=1,\ldots,T$ that 
    \begin{equation}
    r_{ij}^t+r_{i'j'}^t\ge r_{ij'}^t+r_{i'j}^t. \label{eqn:ddd}
    \end{equation}
\end{definition}

We will show that there exists an optimal matching policy that is consistent 
with the partial relation $\succ_{\mathcal{M}_s}$.
For a pair $(i,j)\in\mathcal{A}$, we consider the two sets of pairs $\mathcal{B}_{ij,L}:=\left\{ (i,j')\mid (i,j)\succ_{\mathcal{M}_s} (i,j'),j'\neq j \right\}$ and $\mathcal{B}_{ij,R}=\left\{ (i',j)\mid (i,j)\succ_{\mathcal{M}_s} (i',j), i'\neq i \right\}$, which contain neighboring pairs dominated by $(i,j)$ under $\succ_{\mathcal{M}_s}$. 
Note that any $(i,j'')\notin\mathcal{B}_{ij,L}$ (or $(i'',j) \notin \mathcal{B}_{ij,R}$) is not dominated by $(i,j)$, i.e.,  either dominates, or is incomparable with, $(i,j)$ under $\succ_{\mathcal{M}_s}$.
We also note that $(i,j)$ does not belong to $\mathcal{B}_{ij,L}$ or $\mathcal{B}_{ij,R}$ by definition. 

Given the state $(\mathbf{x},\mathbf{y})$ and feasible matching decision $\mathbf{Q}^t$ in period $t$, the quantity 
$
a_i^{t}:=x_i-\sum_{j'': (i,j'')\notin \mathcal{B}_{ij,L}}q_{ij''}^t
$
represents the remaining quantity of type $i$ demand, after matching $i$ through pairs non-dominated by or incomparable with $(i,j)$.
Likewise, 
$
b_j^{t}:=y_j-\sum_{i'': (i'',j)\notin \mathcal{B}_{ij,R}}q_{i''j}^t
$
is the remaining quantity of type $j$ supply after matching it through pairs non-dominated by or incomparable with $(i,j)$.
We then define a class of policies that respect a partial relation $\succ_{\mathcal{M}_s}$. 

\begin{definition} (Compatibility) 
	We say that an optimal matching policy $\{\mathbf{Q}^{t} \}_{t=1,\ldots,T}$ respects the partial relation $\succ_{\mathcal{M}_s}$ if  
	(i) for all $(i,j)\succ_{\mathcal{M}_s}(i',j)$ and all $t=1,\ldots,T$, either $q_{i'j}^{t}=0$ or $a_i^{t}=0$;
	(ii) for all $(i,j)\succ_{\mathcal{M}_s}(i,j')$ and all $t=1,\ldots,T$, either $q_{ij'}^{t}=0$ or $b_j^{t}=0$.
    \label{def:mms-policy}
\end{definition}

Part (i) of Definition \ref{def:mms-policy} says, unless $a_i^{t}=0$ (i.e., $x_i$ is completely consumed by $(i,j)$ and those non-dominated by or incomparable with $(i,j)$ so that the further matching between $i$ and $j$ is impossible), any $i'$ such that $(i',j)\in\mathcal{B}_{ij,R}$ (i.e., any pair $(i',j)$ dominated by $(i,j)$) will not be matched with $j$. 
In other words, matching of $(i,j)$ is prioritized over all $(i',j)$ such that $(i,j)\succ_{\mathcal{M}_s} (i',j)$.
Similarly, part (ii) implies that matching of $(i,j)$ is prioritized over $(i,j')$ if $(i,j)\succ_{\mathcal{M}_s} (i,j')$.

The following theorem demonstrates the existence of an optimal policy that respects $\succ_{\mathcal{M}_s}$. Then, in that optimal policy, matching of $(i,j)$ is prioritized over $(i',j)$ in any period if $(i,j)\succ_{\mathcal{M}_s}(i',j)$.

\begin{theorem}\label{cor:Neighboring-priority-strong}
    There exists an optimal matching policy that respects $\succ_{\mathcal{M}_s}$.
\end{theorem}

If only the weak modified Monge condition $\succ_{\mathcal{M}}$ is satisfied, we show in the appendix that there exists an optimal policy compatible with the partial relation $\mathcal{M}$ in a weaker sense than Definition \ref{def:mms-policy}.

\subsection*{Optimality of greedily matching a pair of demand and supply.} 

We further provide sufficient conditions for greedy matching between a pair $(i,j)$ to be optimal.

\begin{proposition}\label{prop:greedy-optimal}
    Suppose that the pair $(i,j)$ dominates all its neighboring pair by $\succ_{\mathcal{M}_s}$ (i.e., $(i,j)\succ_{\mathcal{M}_s} (i',j)$ and $(i,j)\succ_{\mathcal{M}_s} (i,j')$ for all $i'\in\mathcal{D}$ and all $j'\in\mathcal{S}$).
    Also suppose $r_{ij}^{t}\ge \max\left\{ \alpha,\beta \right\} r_{ij}^{t+1}$.
    Then, greedy matching between $i$ and $j$ is optimal in all periods. In other words, in any period $t$ with any state $(\mathbf{x},\mathbf{y})$, the optimal matching quantity between $i$ and $j$ is $q_{ij}^{t*}=\min\left\{ x_i,y_j \right\}$.
\end{proposition}

{In the rest of the paper, we refer to the pair $(i,j)$ as a perfect pair if it dominates all its neighboring pair by $\succ_{\mathcal{M}_s}$.
Any other pair is referred to as an imperfect pair. \label{marker:perfect-pair2} }

As an immediate application of Proposition \ref{prop:greedy-optimal}, consider demand and supply types that are specified by their locations in an Euclidean space (e.g., Uber drivers and riders in different locations; products and customers in different locations for Amazon's inventory commingling program). In each period, the reward of matching supply with demand is a fixed prize minus the disutility proportional to the Euclidean distance between the demand location and the supply location (i.e., $r_{ij}^t=R_t - \gamma_t \text{dist}_{ij}$, where $\text{dist}_{ij}$ represents the Euclidean distance between $i$ and $j$). 
If the parameter $\gamma_t$ is decreasing in time, we can verify that a demand type and a supply type from the \emph{same} location forms a perfect pair, and by Proposition \ref{prop:greedy-optimal}, they should be matched as much as possible.

\begin{corollary}\label{cor:euclidean-dist}
    Suppose that the demand and supply types are uniquely characterized by their spatial locations.
    The per-unit matching reward in period $t$ between $i\in\mathcal{D}$ and $j\in\mathcal{S}$ is $r_{ij}^t=R_t-\gamma_t \dist_{ij}$, where $\dist_{ij}$ is the Euclidean distance between $i$ and $j$'s locations. If both $R_t$ and $\gamma_t$ are decreasing in $t$, $i$ and $j$ should be matched greedily in any period.
\end{corollary}

The partial relations defined in Definitions \ref{def:Neighbor-dominance} and \ref{def:Neighbor-dominance-strong} are reminiscent of the classic \emph{Monge sequence} discovered by Gaspard Monge, a French mathematician, in 1781.
\citet{Hoff:1961} provides a necessary and sufficient condition for a static transportation problem to be solvable by a greedy algorithm, in which a permutation (referred to as the Monge sequence) is followed.
The Monge condition provides a priority sequence for \emph{all} the arcs (i.e., demand-supply pairs) in the bipartite network and requires only condition (\ref{eqn:ddd}) of Definition \ref{def:Neighbor-dominance-strong}.
Our Definitions \ref{def:Neighbor-dominance} and \ref{def:Neighbor-dominance-strong} compare two neighboring arcs to determine their priorities in the setting with the dynamic and stochastic arrival of demand and supply types over time. 
Naturally, our conditions may appear more restrictive than the requirements of the Monge sequence because our problem is more complex. 
In particular, to compare $(i,j)$ and $(i',j)$ we require the inequality  $r_{ij}^t-r_{i'j}^t\ge \alpha (r_{ij''}^{t+1}-r_{i'j''}^{t+1})$ to hold for all $j''\in\mathcal{S}$.
(Similarly, to compare $(i,j)$ and $(i,j')$ we require $r_{ij}^t-r_{ij'}^t\ge \beta (r_{i''j}^{t+1}-r_{i''j'}^{t+1})$ to hold for all $i''\in\mathcal{D}$.)
Nevertheless, in subsequent sections, we show that those conditions are satisfied by two classes of problems, namely the \emph{horizontal} model and the \emph{vertical} model.


\begin{remark}
	Our model and results can be generalized to the case with time-dependent carry-over rates.
	Suppose that in period $t$, a fraction $\alpha_t$ of the unmatched demand and a fraction $\beta_t$ of the unmatched supply will carry over to the next period $t+1$, for any type of demand and supply. 
	Then, in Definition \ref{def:Neighbor-dominance}, the conditions $r_{ij}^t-r_{i'j}^t\ge \alpha (r_{ij''}^t-r_{i'j''}^t)$ and $r_{ij}^t-r_{ij'}^t\ge \beta (r_{i''j}^t-r_{i'' j'}^t)$ should be replaced with $r_{ij}^t-r_{i'j}^t\ge \alpha_t (r_{ij''}^t-r_{i'j''}^t)$ and $r_{ij}^t-r_{ij'}^t\ge \beta_t (r_{i''j}^t-r_{i'' j'}^t)$, respectively.
	All subsequent results remain true.
    For example, Proposition \ref{prop:greedy-optimal} still holds if we replace the condition $r_{ij}^{t}\ge \max\left\{ \alpha,\beta \right\} r_{ij}^{t+1}$ with $r_{ij}^{t}\ge \max\left\{ \alpha_t,\beta_t \right\} r_{ij}^{t+1}$.
\end{remark}


\section{Horizontally Differentiated Types}

Consider demand and supply types located in a space $C$. 
Each point in $C$ represents the characteristics of the corresponding (demand/supply) type.
A shorter distance between $i\in\mathcal{D}$ and $j\in\mathcal{S}$ implies a higher unit matching reward in each period.
Thus, the types are ``horizontally'' distributed.

\subsection{Two demand types and two supply types}
\label{subsec:2x2}

We begin with the space $C$ consisting of two distinct locations, namely, locations 1 and 2.
There are two demand types and two supply types, $\mathcal{D}=\left\{ 1,2 \right\}$ and $\mathcal{S}=\left\{ 1,2 \right\}$.
Type 1 demand and type 1 supply share location 1, while type 2 demand and type 2 supply co-locate at location 2. For $k=1,2$, we denote the other index in the set $\left\{ 1,2 \right\}$ by $-k$, i.e., $-k=3-k$.  
Since a shorter distance implies a higher reward, we make the following two assumptions for the rest of this subsection.

\begin{assumption}\label{as:reward-22}
    $r_{kk}^t\ge \max\{r_{k,-k}^t,r_{-k,k}^t\}$, for $k=1,2$. 
\end{assumption} 

The next assumption further compares the unit matching rewards across different periods.

\begin{assumption}
    For any $k\in\left\{ 1,2 \right\}$, all $i\in\mathcal{D}$, $j\in\mathcal{S}$ and $t=1,\ldots,T$, $r_{kk}^t-r_{k,-k}^t\ge r_{ik}^{t+1}-r_{i,-k}^{t+1}$ and $r_{kk}^t-r_{-k,k}^t\ge r_{kj}^{t+1}-r_{-kj}^{t+1}$.\label{as:horizontal22b}
\end{assumption}

From Assumptions \ref{as:reward-22} and \ref{as:horizontal22b} it is straightforward to verify that $(k,k)\succ_{\mathcal{M}_s}(k,-k)$ and $(k,k)\succ_{\mathcal{M}_s}(-k,k)$ for $k=1,2$.
In other words, demand type 1 and supply type 1 form a perfect pair, and so do demand type 2 and supply type 2, while $(k,-k)$ is an imperfect pair, for $k=1,2$. As an application,
consider a premier service and a regular service (e.g., luxury vs. economy car services) provided by crowdsourced suppliers. 
The fares for the two services are $f_p$ and $f_r$, respectively.
The intermediary firm pays the two types of suppliers $c_p$ and $c_r$, respectively.
If the firm offers the premier service to a customer requesting the regular service, the customer will only pay the regular fare (i.e., free upgrading). However, the intermediary firm still needs to pay the premier wage $c_p$ to the premier service provider.
If a customer originally requesting the premier service is offered the regular service, s/he also pays the regular fare, with a possible penalty cost $\pi$ incurred to the firm (monetary compensation, loss of goodwill, etc.)
It is natural to assume that $f_p>f_r$, $c_p>c_r$ and that the margin of the premier service $f_p-c_p$ is higher than that of the regular service $f_r-c_r$. 
Then, the reward for matching a premier customer with a premier supplier (i.e., $f_p-c_p$) is higher than that for matching a premier customer with a regular supplier (i.e., $f_r-c_r-\pi$), and also higher than that for matching a regular customer with a premier supplier (i.e., $f_r-c_p$).
Likewise, matching a regular customer with a regular supply (i.e., $f_r-c_r$) generates more reward than matching a regular customer with a premier supplier (i.e., $f_r-c_p$), and than matching a premier customer with a regular supplier (i.e., $f_r-c_r-\pi$). This verifies Assumption \ref{as:reward-22}, and as a result, Assumption \ref{as:horizontal22b} trivially holds when the parameters are assumed to be time-independent. 

 It follows directly from Proposition \ref{prop:greedy-optimal} that type 1 demand should be matched with type 1 supply as much as possible, before we match type 1 demand with type 2 supply, or type 2 demand with type 1 supply. 
 Likewise, type 2 demand should be matched with type 2 supply greedily. 
 Clearly, after greedy matching between the pair $(k,k)$, there cannot be any positive remaining quantity for both demand type $k$ and supply type $k$ ($k=1,2$).
This observation allows us to collapse the state space:
In period $t$ with the (original) state $(\mathbf{x},\mathbf{y})=(x_1,x_2,y_1,y_2)$, we define the new state as $\mathbf{z}:=(z_1,z_2)$, where $
    z_1=x_1-y_1,
    z_2=y_2-x_2.
$
The quantity $z_1$ describes the imbalance between type 1 demand and type 1 supply. 
A nonnegative $z_1$ represents the remaining quantity of type 1 demand after greedy matching with type 1 supply in period $t$ (the remaining quantity of type 1 supply will be zero).
For a negative value of $z_1$, $z_1^-=-z_1$ is the remaining quantity of type 1 supply after greedy matching with type 1 demand.
Similarly, $z_2^+$ is the remaining quantity of type 2 supply after greedy matching with type 2 demand, whereas $z_2^-$ is the remaining quantity of type 2 demand after greedy matching with type 2 supply.

After the first round of greedy matching in period $t$, if there are remaining type $k$ demand and type $-k$ supply ($k=1,2$) simultaneously (i.e., either $z_1>0$ and $z_2>0$, or $z_1<0$ and $z_2<0$), we will match the two with each other, but not necessarily in a greedy way. 
The intermediary may withhold some type $k$ demand in order to match it with type $k$ supply in a future period (or withhold type $-k$ supply to match with type $-k$ demand in the future).
The amount of type $k$ demand to withhold generally depends on the available amount of type $-k$ supply.
For example, if there is a high level of type $-k$ supply in the current period, it is  unlikely for all of those supply to meet type $-k$ demand (i.e., its best match) in the future, and we may, therefore, use more type $k$ demand to match with type $-k$ supply.  
Symmetrically, the amount of type $-k$ supply to withhold depends on the available type $k$ demand.
Thus, the matching between an imperfect pair is governed by state-dependent match-down-to target levels, where the state-dependency is one-dimensional (e.g., the target level for type $k$ demand depends only on the available type $-k$ supply). \label{discussion:2x2}

To formalize the above discussion, we define $I\!B:=z_1-z_2=x_1+x_2-y_1-y_2$ as the aggregate imbalance between demand and supply.
We describe the structure of the optimal policy as follows.

 \begin{proposition}
    The optimal policy performs two rounds of matching in each period $t$.
\begin{itemize}
    \item {Round 1: Matching of perfect pairs.} 
    
    \quad For $k=1,2$, match type $k$ demand with type $k$ supply greedily.

    \item {Round 2: Matching of an imperfect pair.} 
    \begin{enumerate}[(i)]
    \item No matching in round 2 if $z_1z_2\le 0$. 

    \item If $z_1>0$ and $z_2>0$, match type 1 demand and type 2 supply. There exist protection levels $p_{d,+}^{t,I\!B}$ and $p_{s,+}^{t,I\!B}$ dependent on the imbalance $I\!B$, such that $p_{d,+}^{t,I\!B}-p_{s,+}^{t,I\!B}=I\!B$, and the matching between the pair $(1,2)$ reduces type 1 demand to $\min\left\{ z_1,p_{d,+}^{t,I\!B} \right\}$ and type 2 supply to $\min\left\{ z_2,p_{s,+}^{t,I\!B} \right\}$.

    \item If $z_1<0$ and $z_2<0$, match type 2 demand and type 1 supply. There exist  protection levels $p_{d,-}^{t,I\!B}$ and $p_{s,-}^{t,I\!B}$ dependent on $I\!B$, such that $p_{d,-}^{t,I\!B}-p_{s,-}^{t,I\!B}=I\!B$, and that the matching between the pair $(2,1)$ reduces type 2 demand to $\min\left\{ -z_2, p_{d,-}^{t,I\!B}\right\}$ and type 1 supply to $\min\left\{ -z_1, p_{s,-}^{t,I\!B}\right\}$.
    \end{enumerate}
    \end{itemize}
    \label{prop:opt-horizontal22}
 \end{proposition}

 According to Proposition \ref{prop:opt-horizontal22}, the matching of round 2 is dependent on the state $\mathbf{z}$.
 When $z_1 z_2\le 0$, after round 1, either both type 1 and type 2 supply are depleted, or both type 1 and type 2 demand are depleted. With neither supply nor demand is available, there is no matching in round 2.

 When $z_1>0$ and $z_2>0$, we have remaining quantities of type 1 demand and type 2 supply.
 Part (ii) of Proposition \ref{prop:opt-horizontal22} shows that the matching between the pair $(1,2)$ is characterized by the protection levels $p_{d,+}^{t,I\!B}$ and $p_{s,+}^{t,I\!B}$, which are the target levels to reduce type 1 demand and type 2 supply to, respectively. 
 In the beginning of round 2, if the quantity $z_1$ of available type 1 demand is above $p_{d,+}^{t,I\!B}$, the optimal policy will reduce it to $p_{d,+}^{t,I\!B}$ (by the quantity $z_1-p_{d,+}^{t,I\!B}$) by matching it with type 2 supply. 
 In the mean time, the relation $p_{d,+}^{t,I\!B}-p_{s,+}^{t,I\!B}=I\!B$ guarantees that type 2 supply will be reduced to $z_2-(z_1-p_{d,+}^{t,I\!B})=p_{d,+}^{t,I\!B}-I\!B=p_{s,+}^{t,I\!B}$. 
 If $z_1$ is below $p_{d,+}^{t,I\!B}$,  there is no matching and type 1 demand remains at the level of $z_1$.
 (Note that $z_1\le p_{d,+}^{t,I\!B}$ implies $z_2\le p_{s,+}^{t,I\!B}$, hence type 2 supply remains at the level of  $z_2$.)

 The state-dependent protection levels $p_{d,+}^{t,I\!B}$ and $p_{s,+}^{t,I\!B}$ only depend on the one-dimensional quantity $I\!B=z_1-z_2$ rather than on the full, two-dimensional state $\mathbf{z}$.
The case of $z_1<0$ and $z_2<0$ is symmetric to the case of $z_1>0$ and $z_2>0$. 
We further consider two special cases, for which we will characterize the properties of the protection levels with respect to the state.

\subsection*{Patient demand and supply types}
\label{subsec:patient-ds}

Consider $\alpha=\beta=1$, i.e., demand and supply are infinitely patient and stay until they are matched. 


\begin{proposition}\label{prop:protection-monotone}
    The protection levels $p_{d,+}^{t,I\!B}$ and $p_{d,-}^{t,I\!B}$ for round 2 matching are increasing in the aggregate imbalance $I\!B$. 
    The protection levels $p_{s,+}^{t,I\!B}$ and $p_{s,-}^{t,I\!B}$ are decreasing in $I\!B$. 
    Moreover, the decreasing and increasing rates are no higher than 1.
\end{proposition}

Proposition \ref{prop:protection-monotone} examines the monotonicity of the protection levels with respect to the aggregate imbalance.
We interpret the proposition as follows.

When $I\!B\ge 0$, demand is in excess. A higher value of $I\!B$ suggests more demand over supply.
The chance of a demand type meeting a better match in a future period becomes smaller.
Therefore it becomes more imperative to consume more demand by lowering the protection level for supply.
As a result, the protection levels $p_{s,+}^{t,I\!B}$ and $p_{s,-}^{t,I\!B}$ decrease as $I\!B$ increases.
The rate of decrease, however, is no higher than 1, which implies that the increment in $I\!B$ (i.e., extra demand more than supply) will not be entirely matched in the current period, through reducing the protection level for supply.
The relations $p_{d,+}^{t,I\!B}-p_{s,+}^{t,I\!B}=I\!B$ and $p_{d,-}^{t,I\!B}-p_{s,-}^{t,I\!B}=I\!B$ then immediately imply that $p_{d,+}^{t,I\!B}$ and $p_{d,-}^{t,I\!B}$ are increasing in $I\!B$ with the increasing rates capped by 1.

When $I\!B<0$, supply is in excess. A larger $I\!B$ suggests less supply in excess of demand.
Thus, it is less imperative to consume the excess in supply, implying a higher protection level for demand.

Proposition \ref{prop:protection-monotone} is particularly helpful when demand and supply quantities take integer values. 
In that case, once we obtained the value of the protection level $p_{d,+}^{t,I\!B}$, the protection level $p_{d,+}^{t,I\!B+1}$ is either $p_{d,+}^{t,I\!B}$ or $p_{d,+}^{t, I\!B}+1$, whichever yields higher matching rewards. 



Although we have assumed $\alpha=\beta=1$, Proposition \ref{prop:protection-monotone} is generalizable to the case with arbitrary values of $\alpha$ and $\beta$ as long as the two carry-over rates are equal to each other (i.e., $\alpha=\beta$).

\subsection*{Impatient demand types and patient supply types}
\label{subsec:patient-s-impatient-d}

Consider  $\alpha=0$ and $\beta=1$.
In this case, demand is impatient and is lost if not matched in the current period. 
Thus we only need to record supply levels as the system state.
Then, round 2 matching is fully characterized by protection levels on the supply side only, as shown in the following proposition.
For ease of notation, let $a\wedge b:=\min\left\{ a,b \right\}$ be the smaller of two numbers $a$ and $b$. 

\begin{proposition}\label{prop:protection-lost-demand}
    There exist \emph{state-independent} protection levels $p_{s,+}^t$ and $p_{s,-}^t$ such that in round 2 matching of period $t$,  

    (i)  if $z_1>0$ and $z_2>0$, the optimal matching policy reduces type 2 supply as close to the protection level $p_{s,+}^t$ as possible; the post-matching level of type 2 supply is $\max\left\{ z_2-z_1,z_2\wedge p_{s,+}^t \right\}$;

    (ii)  if $z_1<0$ and $z_2<0$, the optimal matching policy reduces type 1 supply as close to the  to protection level $p_{s,-}^t$ as possible; the post-matching level of $\max\left\{ z_2-z_1,(-z_1)\wedge p_{s,-}^t \right\}$.
\end{proposition}

Proposition \ref{prop:protection-lost-demand} shows that the optimal policy always aims to reduce type 1 supply to the protection level $p_{s,-}^t$, and type 2 supply to $p_{s,+}^t$. This result is generalizable to the case with $\beta\in(0,1)$. 

More specifically, consider the case with $z_1>0$ and $z_2>0$.
In this case we match type 1 demand with type 2 supply in round 2.
According to the proposition, if type 1 demand is ample, the optimal policy will reduce type 2 supply to $z_2\wedge p_{s,+}^t$ (i.e., to $p_{s,+}^t$ if the quantity $z_2$ of available type 2 supply is above $p_{s,+}^t$, or there is no matching if $z_2$ is already no more than $p_{s,+}^t$).
If there is a low level of type 1 demand, however, type 2 supply can be reduced at most by $z_1$ (when all available type 1 demand is matched with type 2 supply) to the level $z_2-z_1$.
The case of $z_1<0$ and $z_2<0$ is symmetric to the case of $z_1>0$ and $z_2>0$.

\subsection{Multiple demand types and supply types}
\label{subsec:multi-type-hor}

We now study the more general case with $m$ demand types and $n$ supply types, all located in the space $C$. 
Here we consider the case where $C$ is a line segment, with its two endpoints denoted by $o$ and $d$, respectively. 
The fitness of matching a demand type $i$ and a supply type $j$ is determined by the \emph{distance} between $i$ and $j$ on $C$.
We consider two distance metrics, but focus on the \emph{directed} distance in this subsection.

\underline{\emph{Undirected} distance}.
This is the shortest distance between the location of $i\in\mathcal{D}$ and $j\in\mathcal{S}$ on  $C$.

\underline{\emph{Directed} distance}.
Suppose that $C$ is endowed with a direction, say, from endpoint $o$ to endpoint $d$ (in short, $o\to d$).
If the location of $i\in\mathcal{D}$ can be reached from the location of $j\in\mathcal{S}$ by traveling along the given direction $o\to d$ (i.e., $i$ is located between $j$ and endpoint $d$),
the distance between $i$ and $j$, denoted by $\dist_{i\leftarrow j}$, is defined as the distance to be travelled by $j$ along the given direction $o\to d$ to reach the location of $i$.  

We now focus on the directed distance, and assume that the unit matching reward between $i\in\mathcal{D}$ and $j\in\mathcal{S}$ is a linearly decreasing function of the distance $\dist_{i\leftarrow j}$ if $j$ can reach $i$ by traveling along the direction $o\to d$, i.e., 
$r_{ij}^t=R^t-\dist_{i\leftarrow j}$.
If $j$ cannot reach $i$ by traveling along the direction $o\to d$, the unit reward is $r_{ij}^t=0$.

It is clear that the optimal matching quantity $q_{ij}^{t*}=0$ if $j$ cannot reach $i$ along the direction $o\to d$.
Next, we compare two pairs of demand and supply, for both of which the supply type can reach the demand type along the direction $o\to d$.

\begin{lemma}
    (i)  Suppose that supply type $j$ can reach both $i$ and $i'$ along the direction $o\to d$.
    Then, $(i,j)\succ_{\mathcal{M}_s} (i',j)$ if and only if along the direction $o\to d$, the distance from $j$ to $i$ is shorter than the distance from $j$ to $i'$.

    (ii)  Suppose that both supply types $j$ and $j'$ can reach type $i$ demand along the direction $o\to d$. 
    Then, $(i,j)\succ_{\mathcal{M}_s} (i,j')$ if and only if along the direction $o\to d$, $j$ is closer to $i$ than $j'$.
    \label{lem:horizontal-dominance}
\end{lemma}

Lemma \ref{lem:horizontal-dominance} suggests that for two neighboring pairs of demand and supply, the pair with a shorter, unidirectional distance should have a higher priority.
It follows from this lemma that any two neighboring pairs are comparable by $\succ_{\mathcal{M}_s}$.

\begin{proposition}\label{prop:horizontal-structure}
    (i)  If $\dist_{i\leftarrow j}< \dist_{i\leftarrow j'}$, the optimal policy matches $(i,j)$ before $(i,j')$.
    If $\dist_{i\leftarrow j}<\dist_{i'\leftarrow j}$, the optimal policy matches $(i,j)$ before $(i',j)$.

    (ii)  Suppose that $R^t$ decreases in $t$. If there are no other demand or supply types located between $i\in\mathcal{D}$ and $j\in\mathcal{S}$ on $C$, $i$ and $j$ should be matched with each other greedily, i.e., $q_{ij}^{t*}=\min\left\{ x_i,y_j \right\}$.
\end{proposition}

Proposition \ref{prop:horizontal-structure} \label{after_prop:horizontal-structure} prescribes a priority hierarchy for the optimal matching policy, by classifying the pairs of demand and supply into priority tiers. 
Let the set of tier 0 pairs, denoted by $\mathcal{A}_0$, be those not dominated by any neighboring pair under $\succ_{\mathcal{M}_s}$.
Recursively, we can define tier $k$ pairs, denoted by $\mathcal{A}_k$, as those pairs that belong to $\mathcal{A}\backslash \bigcup_{\ell=1}^{k-1} \mathcal{A}_{\ell}$ and are not dominated by any other neighboring pairs in $\mathcal{A}\backslash \bigcup_{\ell=1}^{k-1} \mathcal{A}_{\ell}$.
Suppose that there are a total number of $K$ tiers.
The optimal policy always matches the pairs in $\mathcal{A}_{k-1}$ before it moves on to match the pairs in $\mathcal{A}_k$, for $k=1,\ldots,K$.
Moreover, if a pair $(i,j)\in\mathcal{A}_{k-1}$ is not matched to the full extent (i.e., there are remaining quantities of both type $i$ demand and type $j$ supply), any pair of the form $(i',j)$ or $(i,j')$ in $\mathcal{A}\backslash \bigcup_{\ell=1}^{k-1} \mathcal{A}_{\ell}$ will not be matched (i.e., with a zero matching quantity) in the optimal policy.

\subsection*{A heuristic idea}

\label{marker:horizontal-heu}

Proposition \ref{prop:horizontal-structure} \label{after_prop:horizontal-structure} provides a partial characterization of the optimal policy with respect to the priority structure, but does not prescribe how much to match for each pair of demand and supply types.
Motivated by this proposition, we briefly describe a heuristic idea to compute the optimal matching decisions. 
For a given period $t$, we consider $i\in\mathcal{D}$ and $j\in\mathcal{S}$ both located on the line segment $C$ such that $i$ is accessible from $j$ along the given direction.  
When matching $i$ with $j$, we may want to reserve some type $i$ demand (resp., type $j$ supply) for future supply types (resp., demand types) located between $i$ and $j$ on the line segment $C$. 
But we may not want to reserve type $i$ demand (resp., type $j$ supply) for any supply type $j'$ (resp. demand type $i'$) located outside the segment between $i$ and $j$, due to the lower priority of the pair $(i,j')$ (resp., $(i',j)$) than the pair $(i,j)$ (see Proposition \ref{prop:horizontal-structure}).
As a heuristic, we determine the matching between $i$ and $j$ by considering a subproblem P$(i,j)$ that comprises only demand type $i$, supply type $j$ and the types located between $i$ and $j$ on $C$.
According to Proposition \ref{prop:horizontal-structure}, we should not match $i$ with $j$ until there is no remaining quantity for any demand type $i'$ or supply type $j'$ located between $i$ and $j$.
Thus, we assume that in the subproblem P$(i,j)$, all types except demand type $i$ and supply type $j$ have zero remaining quantity.
Analogous to the $2\times 2$ model in Section \ref{subsec:2x2}, we can show that the optimal matching between $i$ and $j$ is characterized by a protection level $p_{ij,d}^{t,I\!B}$ on type $i$ demand and $p_{ij,s}^{t,I\!B}$ on type $j$ supply, with both protection levels dependent on $I\!B$ (which is the imbalance between type $i$ demand and type $j$ supply) and $p_{ij,d}^{t,I\!B}-p_{ij,s}^{t,I\!B}=I\!B$. 
More specifically, we will match $i$ with $j$ until type $i$ demand is reduced to $p_{ij,d}^{t,I\!B}$ and type $j$ supply is reduced to $p_{ij,s}^{t,I\!B}$, or as close as possible.

Next, we outline the heuristic matching procedure for a period $t$, assuming that the protection levels $p_{ij,d}^{t,I\!B}$ and $p_{ij,s}^{t,I\!B}$ are already obtained for all $i\in\mathcal{D}$ and $j\in\mathcal{S}$.

\begin{heuristic}
   (Prioritized matching for the horizontal model) 
\begin{algorithmic}[]
    \FOR{$k=1$ \TO $K$}
        \FOR{each pair $(i,j)$ in priority tier $k$}
            \STATE Match $i$ with $j$ until type $i$ demand is reduced to $p_{ij,d}^{t,I\!B}$ and type $j$ supply is reduced to $p_{ij,s}^{t,I\!B}$, or as close as possible
        \ENDFOR
    \ENDFOR
\end{algorithmic}
\label{heu:horizontal}
\end{heuristic}

Within each priority tier, it does not matter which pair we match first, because the matching of one pair does not affect the subproblem for another pair within the same tier. 
The computation of the protection levels for each subproblem P$(i,j)$, however, remains challenging. In the appendix, we discuss a heuristic method that converts the subproblem P$(i,j)$ to a $2\times 2$ model by consolidating demand and supply types.
Next, we discuss a couple of applications of the horizontal model.

\subsection*{Commuter car pooling platforms}\label{subsec:carpool}

\begin{figure}[htpb]
\FIGURE
{   \includegraphics*[scale=0.5]{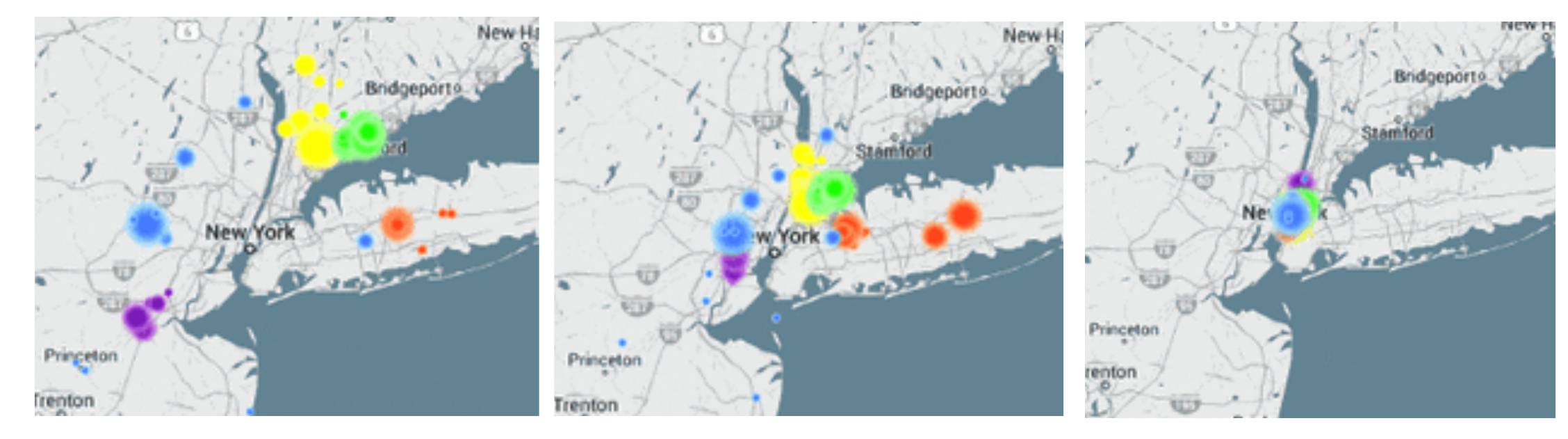}}
  {Commuting pattern in the mornings of weekdays, from suburban areas to NYC. \label{fig:Commute_pat}}
  {
  Each colored dot represents a group of commuters on the same commuting route.  
  }
\end{figure}

Carpooling platforms such as iCarpool and UberPool match a driver heading to a destination with several riders to the same destination (or in the same direction).
Commuting patterns of many cities indicate that drivers and riders often share the same destination. For example, Figure \ref{fig:Commute_pat} displays the New York City commuting pattern in the mornings of weekdays, from which we see that commuters travel from different suburban areas in the same direction to the city.
In this case, the directed line segment $C$ is corresponding to the route that starts from a suburban area (i.e., endpoint $o$) and ends in the city (i.e., endpoint $d$).
Drivers, who may be commuters themselves, pick up riders along the route.

Moreover, if all riders share the same destination (e.g., the city) and a driver picks up riders along the way to the destination, the closer a rider to the driver, the shorter the waiting time for the rider and the higher the payment for the ride (due to the longer distance travelled by the rider).
To formalize this intuition, we generalize the reward function $r_{ij}^t=R^t-\dist_{i\leftarrow j}$ mentioned above as follows.
If $j$ can reach $i$ along the direction $o\to d$, the unit reward of matching $i$ with $j$ is $r_{ij}^t=R_i^{t}-\dist_{i\leftarrow j}$. Otherwise, $r_{ij}^t=0$.
Here, $R_i^t$ represents the reward resulting from the match and is dependent on the attribute of type $i$ demand, e.g., the travel distance of the rider (from rider $i$'s initial location to the destination). 
The second term $-\dist_{i\leftarrow j}$ represents the disutility proportional to the traveling distance by the driver for the pickup, e.g., a longer distance implies a longer roaming time for the driver and longer waiting time for the rider.
Following similar analysis, we can show that $(i,j)\succ_{\mathcal{M}_s} (i',j)$ if and only if $i$ is closer to $j$ than $i'$ along the direction of the route, and that $(i,j)\succ_{\mathcal{M}_s} (i,j')$ if and only if $j$ is closer to $i$ than $j'$ along the direction of the route.


\subsection*{Product/Service general upgrading}

\begin{figure}
\centering
\begin{minipage}{.5\textwidth}
  \centering
 \includegraphics*[scale=0.2]{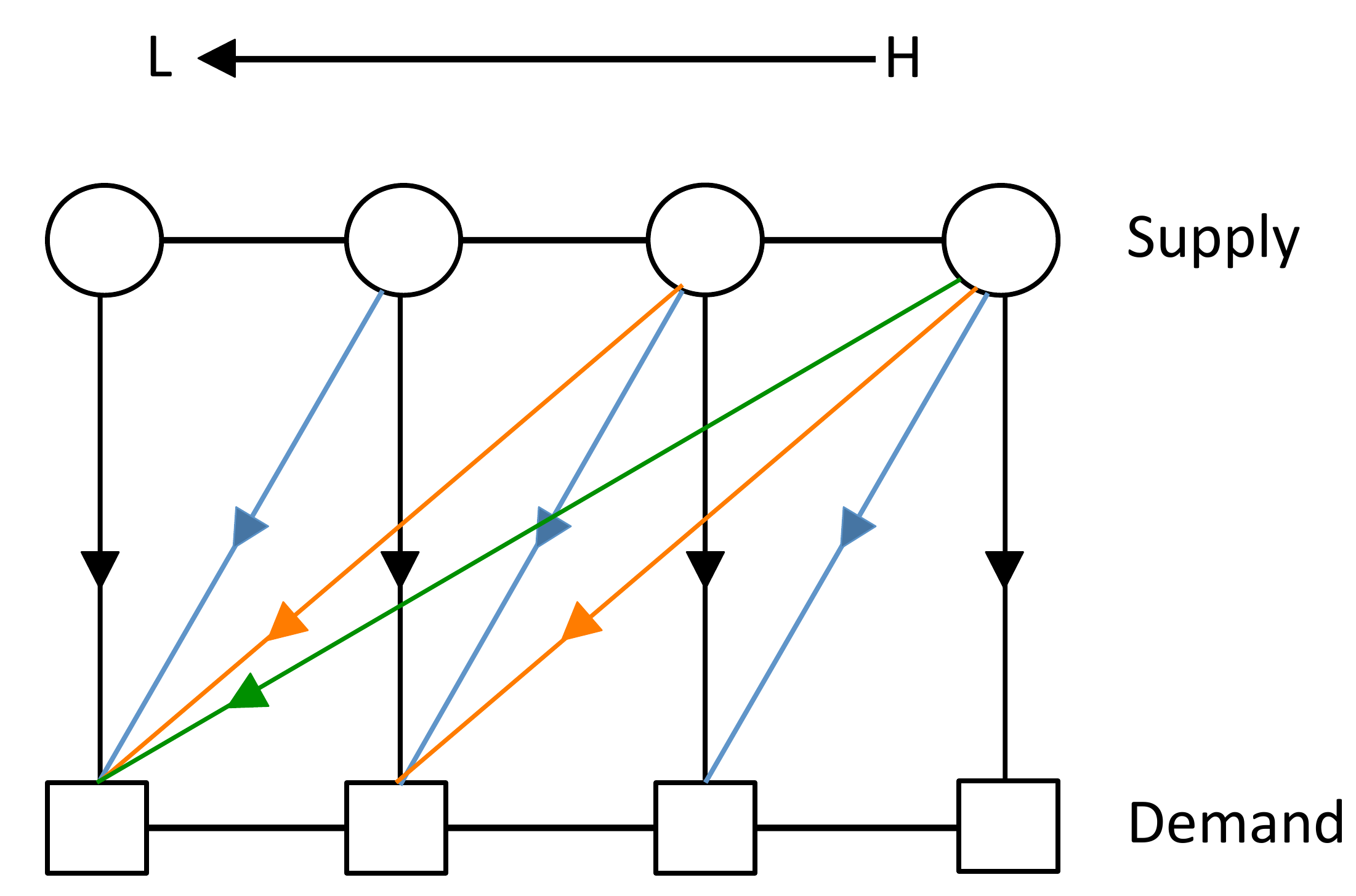}
    \caption{Product upgrade.}
    \label{fig:upgrade2}
    \end{minipage}%
\begin{minipage}{.5\textwidth}
  \centering
    \includegraphics[scale=0.2]{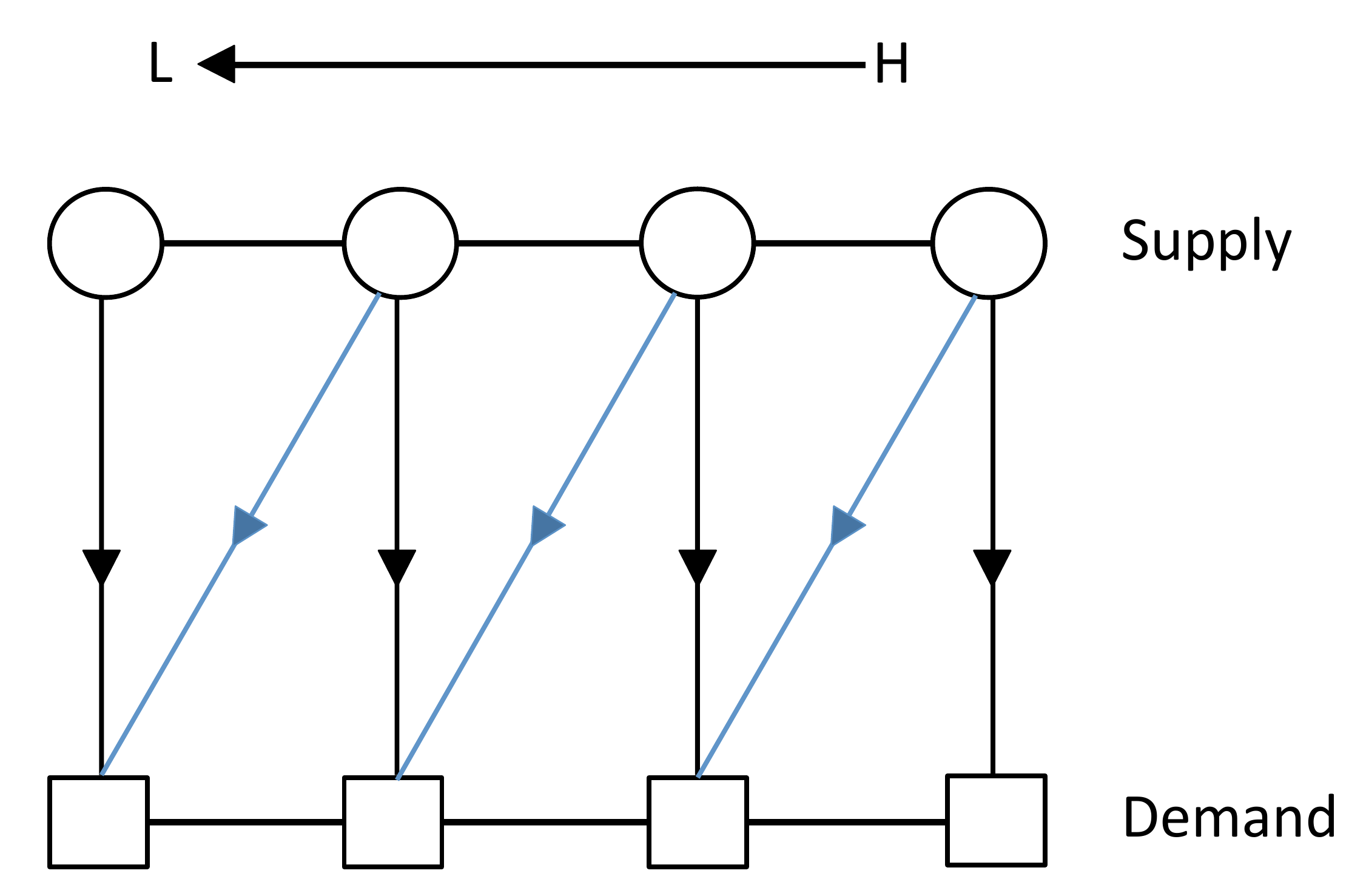}
    \caption{One-level upgrading.}
    \label{fig:upgrade1}
\end{minipage}
\end{figure}

Upgrading uses a high-class supply to fulfill a  low-class demand, which is widely adopted in the business practice, e.g., in travel industries (see, e.g.,  \citealt{yu2014dynamic}) and in production/inventory settings (see, e.g., \citealt{bassok1999single}). 
Figure \ref{fig:upgrade2} illustrates such a model that allows \emph{general} upgrading (see \citealt{yu2014dynamic}).
In this model, product classes $1,\ldots,n$ are indexed according to the descending order of quality. 
Class $i$ products are intended for the customer segment $i$. Thus it is mostly desirable to satisfy a class $i$ customer demand using a class $i$ product,  more desirable to satisfy a class $i$ demand using a class $k$ product than using a class $k-1$ product ($1<k\le i$), and infeasible to satisfy a class $i$ demand using a class $\ell$ product ($\ell>i$).

In contrast to the existing works in the literature where the supply side is either fixed or controlled through replenishing decisions, there are many settings in which new supply arrives randomly. 
For example, ride-hailing platforms such as Uber  randomly have new drivers coming online or existing drivers completing a service and becoming available, who provide differentiated types of service (UberX, UberSELECT, UberBLACK, etc.; a more premium vehicle can be used to serve a less premium customer class through upgrading).
Car rental companies may have random supply levels due to early/late return of cars by customers.
Airlines and hotels can also have random ``arrival'' of supply due to customer cancellations.

The problem of \emph{general} upgrading has the structure of a directed line segment in the product line space. \label{marker:general-ug-explain}
Class $i$ demand and class $i$ product share the same location on the line segment, and the lower the class index, the closer the class is located to the endpoint $o$.
For $j\le i$, let $c_j$ be the unit purchase cost for class $j$ product and $f_i^t$ be the fare paid by class $i$ customers in period $t$.  
The unit profit for assigning $i$ to $j$ is then $r_{ij}^t=f_i^t-c_j=f_i^t-c_i-(c_j-c_i)$.
If we define $\dist_{i\leftarrow j}:=c_j-c_i$ as the distance between $i$ and $j$ and $R_t^i:=f_i^t-c_i$ as the unit profit from a class $i$ customer being satisfied by a class $i$ product, then the reward structure reduces to the one we have already considered in the previous application, i.e., $r_{ij}^t=R_i^t-\dist_{i\leftarrow j}$.
Then the optimal policy will satisfy a customer with a product class that is the same as or closer to the originally requested product, and assign a product to a customer the same as or closer to the customer class that the product is intended for. 

\subsection*{One-level product upgrading}


\citet{shumsky2009dynamic} study a capacity management problem in which each customer class can only be upgraded one level higher. 
Figure \ref{fig:upgrade1} demonstrates the structure of such a problem. 
Again, we can think of the customer classes and product classes located on a line segment $C$, where the class $j$ product and its intended customer class $j$ share the same location. 

The infeasibility of upgrading with more than one level makes the problem structurally different from the general upgrading problem. 
In the one-level-up upgrading problem, the reward structure is the same as before for any $i, j=i-1$, i.e., $r_{i, i-1}^t=R_i^t-\dist_{i\leftarrow i-1}$, with $R_i^t$ decreasing in time.  But for any $i, j\neq i-1$, $r_{ij}^t=0$, different from the general upgrading problem. As a result, in the one-level-up upgrading problem, two neighboring pairs of demand and supply, $(i,i)$ and $(i,i-1)$, are not necessarily comparable under $\succ_{\mathcal{M}_s}$. Specifically, $r_{ii}^t-r_{i,i-1}^t=\dist_{i\leftarrow i-1} - \dist_{i\leftarrow i} = \dist_{i\leftarrow i-1}$,
and $r_{i+1,i}^t-r_{i+1,i-1}^t=r_{i+1,i}^t=R_i^t-\dist_{i+1\leftarrow i}$.
Therefore, a necessary condition for $(i,i)\succ_{\mathcal{M}_s} (i,i-1)$ is that $R_i^t\le \dist_{i+1\leftarrow i}+\dist_{i\leftarrow i-1}=\dist_{i+1\leftarrow i-1}$, which may not be guaranteed in general.

The above argument implies that it may \emph{not} be optimal to prioritize the matching between a pair of demand and supply intended for each other (i.e., demand type $i$ with supply type $i$), over upgrading (i.e., demand type $i$ with supply type $i-1$, or demand type $i+1$ with supply type $i$), when the supply is random. 
In the followings, we investigate the loss of optimality caused by enforcing the aforementioned priority structure.

\begin{remark}
    Our priority structure for the general upgrading problem is consistent with \cite{bassok1999single}, who consider a single-period version of the problem with general upgrading. They prove that greedy matching along the specified priority structure (i.e., a product-customer pair has a higher priority if they are closer to each other) is optimal by showing that such a priority structure leads to a classical Monge sequence. However, similar to our arguments above, there no longer exits a Monge sequence when only one-level upgrading is allowed, even
    in the single-period problem considered by \cite{bassok1999single}. 
\end{remark}

Let $\mathcal{P}^{\text{IOU}}$ be the set of matching policies that prioritizes \underline{i}ntended pairs \underline{o}ver \underline{u}pgrading.
More specifically, a policy belongs to $\mathcal{P}^{\text{IOU}}$ if and only if it matches $(i,i)$ greedily before $(i,i-1)$ and $(i+1,i)$.
That is, unless $q_{ii}^t=\min\left\{ x_i,y_j \right\}$, the policy has $q_{i,i-1}^t=q_{i+1,i}^t=0$.
In the following proposition, we show that by enforcing the best policy in $\mathcal{P}^{\text{IOU}}$, the optimality loss is no more than 50\%.


\begin{proposition}\label{prop:one-level}
    Suppose that $r_{ii}^t\ge \max\left\{ \alpha,\beta \right\}r_{ii}^{t+1}$ for any $i=1,\ldots,n$ and $t=1,\ldots,T-1$.\footnote{For time-dependent carry-over rates, we replace the condition $r_{ii}^t\ge \max\left\{ \alpha,\beta \right\}r_{ii}^{t+1}$ with $r_{ii}^t\ge \max\left\{ \alpha_t,\beta_t \right\}r_{ii}^{t+1}$.} There exists a matching policy belonging to $\mathcal{P}^{\text{IOU}}$, such that it retains at least 50\% of the total expected reward under the optimal matching policy.
\end{proposition}

\section{Vertically Differentiated Demand and Supply Types}

In this section, we consider \emph{vertically} differentiated demand and supply types.
Each demand/supply type is associated with a ``quality'' level, and generates a higher reward if it is matched with a supply or demand type of a higher quality.
In other words, we have the reward function $r_{ij}^t=f^t(a_i,b_j)$ increasing in $a_i$ and $b_j$, where $a_i$ represents the quality of demand type $i$ and $b_j$ represents the quality of supply type $j$.
For simplicity, we consider a linearly additive reward function $r_{ij}^t=f_d^t(a_i)+f_s^t(b_j)$, where $f_d^t$ and $f_s^t$ are increasing in $a_i$ and $b_j$, respectively.
We write $r_{id}^t:=f_d^t(a_i)$ and $r_{js}^t:=f_s^t(b_j)$. Later we consider the generalization in which the reward structure can be nonlinear (see \S\ref{subsec:nonadd}).


Without loss of generality, we assume that a demand/supply type with a smaller index has a higher quality. 
That is, $r_{1d}^t>r_{2d}^t>\cdots>r_{md}^t$ and $r_{1s}^t>r_{2s}^t>\cdots>r_{ns}^t$.
In addition, we make the following assumption, which requires the quality difference between types to weakly decrease in time.
For convenience of notation, we define $r_{m+1,d}^t=r_{n+1,s}^t:=0$ for all period $t$.

\begin{assumption}\label{as:vert-quality}
    For any $t=1,\ldots,T-1$, $i=1,\ldots,m$ and $j=1,\ldots,n$, we assume that $r_{id}^t-r_{i+1,d}^t\ge \alpha(r_{id}^{t+1}-r_{i+1,d}^{t+1})$ and $r_{js}^t-r_{j+1,s}^t\ge \beta(r_{js}^{t+1}-r_{j+1,s}^{t+1})$.
\end{assumption}

Assumption \ref{as:vert-quality} enables us to compare neighboring pairs of demand and supply under $\succ_{\mathcal{M}_s}$.

\begin{lemma}\label{lem:vert-comp}
    Under Assumption \ref{as:vert-quality}, $(i,j)\succ_{\mathcal{M}_s} (i,j')$ for all $j<j'$, and $(i,j)\succ_{\mathcal{M}_s} (i',j) $ for all $i<i'$.
\end{lemma}

\begin{figure}[htb!]
\FIGURE
{\includegraphics*[scale=0.2]{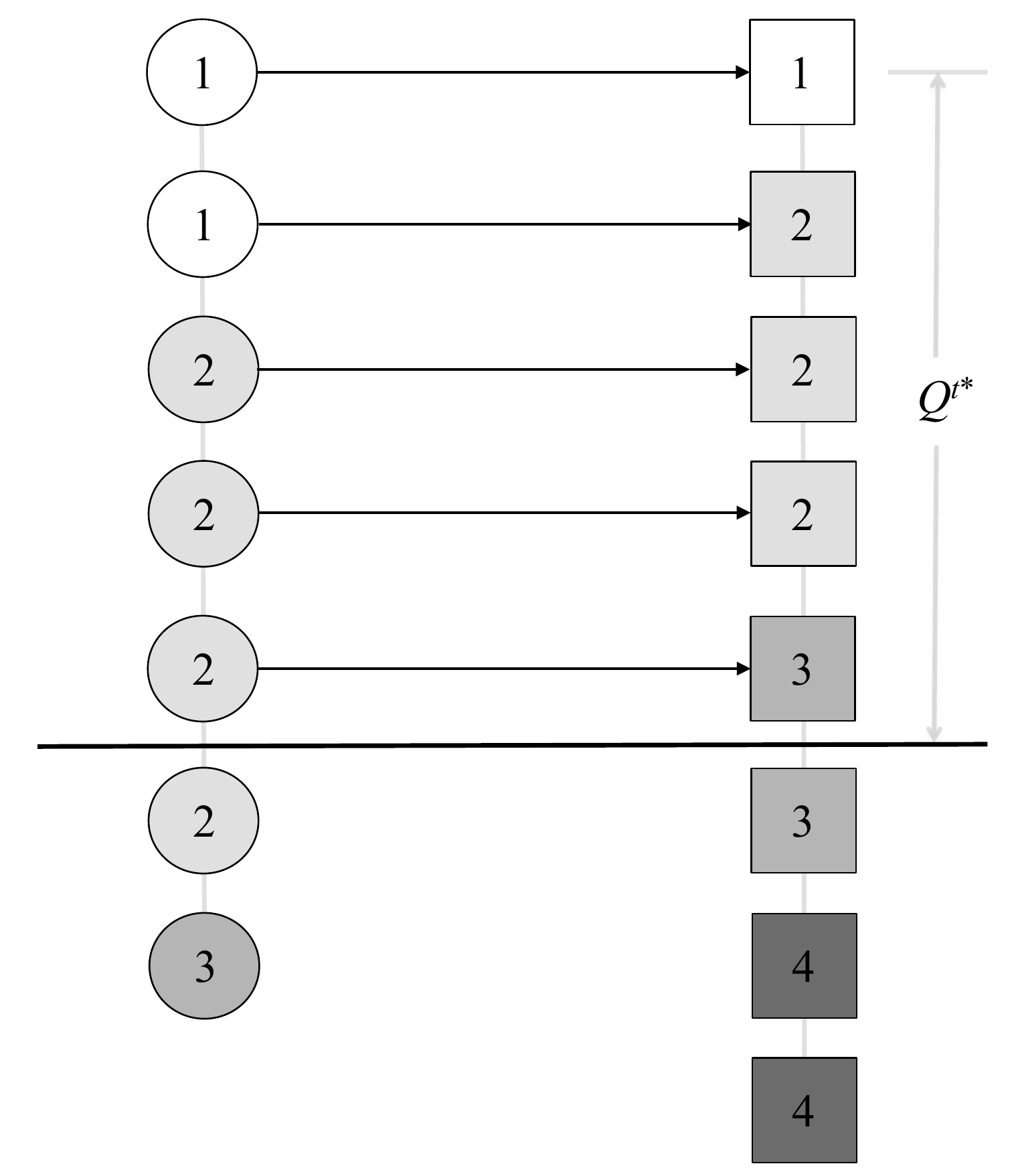}}
{Line up, match up (to a ``match-down-to" level).\label{fig:lineup}}
{}
\end{figure}

It follows from Lemma \ref{lem:vert-comp} and Theorem \ref{cor:Neighboring-priority-strong} that a higher-quality supply type $j$ will prioritize over a lower-quality supply type $j'$ for matching with any demand type $i$. (Symmetrically, a higher-quality demand type $i$ will prioritize over a lower-quality demand type $i'$ for matching with any supply type $j$.) 
Then if we line up demand types and supply types separately in ascending order of their indices (i.e., descending order of their quality levels), the optimal policy will match demand with supply from the top and down to some level (see Figure \ref{fig:lineup}), referred to as ``top-down" matching.

\vspace{-2em}

\begin{proposition}\label{prop:vert-top-down}
  The top-down matching is optimal, and the optimal matching quantities in period $t$ are fully determined by the optimal total matching quantity $Q^{t*}:=\sum_{i\in\mathcal{D},j\in\mathcal{S}} q_{ij}^{t*}$.
\end{proposition}

\begin{remark}
    With the additive reward $r_{ij}^t=r_{id}^t + r_{js}^t$, 
   matching a unit of type $i$ demand with any supply type always contributes the reward $r_{id}^t$, regardless of which supply type it matches with. 
    Likewise, a unit of type $j$ supply always contributes $r_{js}^t$ upon matching.
    Thus, the optimal decision $\mathbf{Q}^{t*}$ is not unique. 
    As long as we determine the total quantity $Q^{t*}$ of demand and supply to match, the total matching reward is determined and independent of how we pair up the demand types and supply types.
    Nevertheless, the top-down matching with higher quality types prioritized over lower quality types offers a \emph{stable matching} (see \citealt{RS90}) in terms of incentive compatibility for each individual within the centrally determined matching pool of size $Q^{t*}$. 
\end{remark}

Like in the horizontal model, we can also view the optimal matching policy in the form of a ``match-down-to threshold'' structure.
For ease of notation, we define $\tilde{x}_i:=\sum_{k=1}^{i} x_k$ for $i=1,\ldots, m$ and $\tilde{y}_j:=\sum_{k=1}^{j} y_k$ for $j=1,\ldots,n$ (with $\tilde{x}_0:=0$ and $\tilde{y}_0:=0$) as the \emph{transformed} system state. (See Online Appendix \ref{app:vert-alt} for an alternative formulation of the vertical model based on the transformed state $(\tilde{\mathbf{x}},\tilde{\mathbf{y}}):=(\tilde{x}_1,\ldots,\tilde{x}_m,\tilde{y}_1,\ldots,\tilde{y}_n)$.)
Under the top-down matching, type $i$ demand is matched with type $j$ supply only when 
$\tilde{y}_j>\tilde{x}_{i-1}$ and $\tilde{y}_{j-1}<\tilde{x}_i$.\footnote{If $\tilde{y}_j\le \tilde{x}_{i-1}$, type $j$ supply has been fully consumed when the optimal policy starts to match type $i$ demand. Likewise, if $\tilde{y}_{j-1}\ge \tilde{x}_i$, type $i$ demand is fully consumed before the optimal policy starts to match type $j$ supply.}
When type $i$ demand matches with type $j$ supply, there exist state-dependent protection levels $p_{i,d}^t$ and $p_{j,s}^t$ such that by matching $i$ and $j$, the optimal policy reduces type $i$ demand to $p_{i,d}^t$ and type $j$ supply to $p_{j,s}^t$, or as close as possible.\footnote{Both $p_{i,d}^t$ and $p_{j,s}^t$ depend on $\tilde{x}_i-\tilde{y}_j$, $\tilde{x}_{i},\ldots,\tilde{x}_{m}$, and $\tilde{y}_{j},\ldots,\tilde{y}_n$. We suppress the dependency for ease of notation.}

The optimal total matching quantity $Q^{t*}$ (or the protection levels) is state-dependent and challenging to compute due to the high dimensionality of the problem.
Next, we propose a heuristic method to determine the matching decisions in each period, and also explore the structural properties of the optimal total quantity $Q^{t*}$ with respect to the state $(\mathbf{x},\mathbf{y})$ for two special cases.  


\subsection{The one-step-ahead policy}
\label{subsec:1sa}

Let the greedy policy be defined as one that retains the top-down structure and matches demand with supply as much as possible in every period.
We consider the following one-step-ahead policy, which improves on the greedy policy.
Under this policy, the intermediary assumes in every period $t$ that from the next period $t+1$ until the end of the time horizon the greedy policy will be enforced. 
(In the next period $t+1$, however, instead of using the greedy decisions the policy will use the ``optimal'' policy that maximizes the total expected reward-to-go, provided that greedy matching is enforced from period $t+2$ to the end of the horizon.)

It is well-known that high dimensional dynamic programs are extremely hard to solve due to the difficulty in evaluating the value functions.
The one-step-ahead policy provides an easy-to-compute approximation to the optimal value function (i.e., the optimal reward that can be received from the next period to the end of the horizon).
Consider period $t+1$ with state $(\mathbf{x},\mathbf{y})$.  
For any realization $\omega^{[t+1,T]}$ of demand and supply of each type from period $t+1$ to period $T$, we can calculate the total reward received under the greedy policy  (denoted by $R^{[t+1,T]}(\mathbf{x},\mathbf{y},\omega^{[t+1,T]})$).  
We randomly draw $N$ samples $\omega^{[t+1,T]}_1,\ldots,\omega_N^{[t+1,T]}$, and approximate the optimal value function by $V_{t+1}(\mathbf{x},\mathbf{y})\approx \frac{1}{N} \sum_{k=1}^{N} R^{[t+1,T]}(\mathbf{x},\mathbf{y},\omega_k^{[t+1,T]})$.
Based on this approximation, we can solve the total expected reward maximization problem in period $t$, and use its solution to approximate the optimal total matching quantity in period $t$. 
More technical details on the one-step-ahead policy are in Appendix \ref{sec_app:1sa}. Since the one-step-ahead policy improves upon the greedy policy, naturally we have the following: 

\begin{proposition}\label{prop:1sa}
   The one-step-ahead policy has a higher expected reward than the greedy policy. 
\end{proposition}

\subsection{Patient demand and supply types}

Consider the case with $\alpha=\beta=1$, in which all demand types and supply types are infinitely patient. (The results hold for the case with equal carry-over rates, $\alpha=\beta\in(0,1)$, with more complicated notation.)
We study the monotonicity properties of the optimal total matching quantity $Q^{t*}$ with respect to the state $(\mathbf{x},\mathbf{y})$.
We write $Q^{t*}$ as $Q^{t*}(\mathbf{x},\mathbf{y})$ to reflect its dependency on the state $(\mathbf{x},\mathbf{y})$.
We also define $\frac{\partial Q^{t*}(\mathbf{x},\mathbf{y})}{\partial x_i}:=\lim\sup_{\varepsilon\to 0+} [Q^{t*}(\mathbf{x}+\varepsilon\mathbf{e}_i,\mathbf{y}) - Q^{t*}(\mathbf{x},\mathbf{y}) ]/\varepsilon$ and $\frac{\partial Q^{t*}(\mathbf{x},\mathbf{y})}{\partial y_j}:=\lim\sup_{\varepsilon\to 0+} [Q^{t*}(\mathbf{x},\mathbf{y}+\varepsilon\mathbf{e}_j) - Q^{t*}(\mathbf{x},\mathbf{y}) ]/\varepsilon$ if $Q^{t*}$ is not differentiable.

\begin{proposition}
    The optimal total matching quantity $Q^{t*}(\mathbf{x},\mathbf{y})$ is increasing in the demand level $x_i$ and supply level $y_j$ for all $i\in\mathcal{D}$ and $j\in\mathcal{S}$, with the increasing rate no greater than one, i.e., $0\le \frac{\partial Q^{t*}(\mathbf{x},\mathbf{y})}{\partial x_i}\le 1$ and $0\le \frac{\partial Q^{t*}(\mathbf{x},\mathbf{y})}{\partial y_j}\le 1$ for all $i\in\mathcal{D}$ and $j\in\mathcal{S}$. 
    Moreover, $Q^{t*}(\mathbf{x},\mathbf{y})$ is more sensitive to demand and supply types of higher quality, i.e., $\frac{\partial Q^{t*}(\mathbf{x},\mathbf{y})}{\partial x_{i}}\ge \frac{\partial Q^{t*}(\mathbf{x},\mathbf{y})}{\partial x_{i+1}}$ for $i=1,\ldots,m-1$ and $\frac{\partial Q^{t*}(\mathbf{x},\mathbf{y})}{\partial y_j} \ge \frac{\partial Q^{t*}(\mathbf{x},\mathbf{y})}{\partial y_{j+1}}$ for $j=1,\ldots,n-1$.
    \label{prop:vertical-mono}
\end{proposition}


\subsubsection*{The one-step-ahead policy.} 
With additive rewards, a unit of type $i$ demand always generates the reward $r_{id}^{t}$ in period $t$ regardless of which supply type it pairs with.
Similarly, a unit of type $i$ demand always generates the reward $r_{js}^t$ in period $t$.
Therefore, it may be less crucial which demand/supply types we withhold for future supply/demand, but more important how much demand/supply we withhold.
When we match type $i$ demand with type $j$ supply, the trade-off is between receiving the unit matching reward $r_{ij}^t$ in the current period $t$, and withholding type $i$ demand (or type $j$ supply) for some future supply $j'<j$ (or some future demand $i'<i$) so that the latter will not be delayed in getting paired or be lost in the case of $\alpha=\beta<1$.  
Intuitively, the amount of demand/supply we want to withhold is determined by the imbalance between demand and supply, e.g., the more demand (supply) in excess of supply (demand), the less supply (demand) we want to withhold. This intuition is shown to hold for the one-step-ahead policy. 

Under the one-step-ahead policy, we show that the matching decision in each period can be described in terms of protection levels which only depend on the aggregate imbalance between demand and supply, when both demand and supply are patient. 

Let $I\!B:=\tilde{x}_m-\tilde{y}_n$ be the aggregate imbalance between demand and supply.

\begin{proposition}
    For any $i\in\mathcal{D}$, $j\in\mathcal{S}$ and $t=1,\ldots,T$, there exist protection levels $p_{ij,d}^{t,I\!B}$ and $p_{ij,s}^{t,I\!B}$ dependent on $I\!B$ such that $ p_{ij,d}^{t,I\!B} - p_{ij,s}^{t,I\!B} = \tilde{x}_m-\tilde{y}_n$, and that in period $t$: 

    (i)  The one-step-ahead policy matches $i$ and $j$ only if $\tilde{x}_i>\tilde{y}_{j-1}$ and $\tilde{y}_j>\tilde{x}_{i-1}$; 

    (ii)  When the one-step-ahead policy matches $i$ and $j$, it aims to reduce the total available demand to the protection level $p_{ij,d}^{t,I\!B}$ and the total available supply to  $p_{ij,s}^{t,I\!B}$ or as close as possible.

    (iii)  $ p_{ij,d}^{t,I\!B}$ is increasing in $I\!B$ with the increasing rate no greater than 1, and $p_{ij,s}^{t,I\!B}$ is decreasing in $I\!B$ with the decreasing rate no greater than 1.
    \label{prop:vert-1sa-structure}
\end{proposition}

Proposition \ref{prop:vert-1sa-structure} shows that the one-step-ahead policy has limited state-dependency.
Instead of depending on the full state $(\mathbf{x},\mathbf{y})$, the protection levels $p_{ij,d}^{t,I\!B}$ and $p_{ij,s}^{t,I\!B}$ only depend on the aggregate imbalance $I\!B$ between demand and supply.

As indicated by part (i) of Proposition \ref{prop:vert-1sa-structure}, type $i$ demand is matched with type $j$ supply only when  $\tilde{x}_i>\tilde{y}_{j-1}$ (in which case type $i$ demand is not fully consumed by types $1,\ldots,j-1$ demand) and $\tilde{y}_j>\tilde{x}_{i-1}$ (in which case type $j$ supply is not fully consumed by types $1,\ldots,i-1$ demand).

Under the top-down matching, type $i$ demand (or type $j$ supply) would never be used unless all higher-quality demand types (or supply types) have run out.
Thus, immediately prior to the matching between $i$ and $j$, a total quantity $\tilde{x}_{i-1}\vee \tilde{y}_{j-1}:=\max\{\tilde{x}_{i-1},\tilde{y}_{j-1}\}$ of demand and the same quantity of supply have been consumed.
There is a total quantity $\tilde{x}_m-\tilde{x}_{i-1}\vee \tilde{y}_{j-1}$ of remaining demand and a total quantity $\tilde{y}_n-\tilde{x}_{i-1}\vee \tilde{y}_{j-1}$ of remaining supply.

As we match along the top-down structure in period $t$, type $i$ demand and type $j$ supply would be matched to the maximum extent when the total quantity reaches $\tilde{x}_i\wedge \tilde{y}_i:=\min\{\tilde{x}_i,\tilde{y}_i\}$ (in that case, either type $i$ demand or type $j$ supply runs out).
This happens when the aggregate demand level reduces to $\tilde{x}_m-\tilde{x}_i\wedge \tilde{y}_j$, or equivalently, the aggregate supply level reduces to $\tilde{y}_n-\tilde{x}_i\wedge \tilde{y}_j$.

Overall, the one-step-ahead policy matches $i$ with $j$ in the following intuitive way:
\begin{itemize}
\item If $\tilde{x}_m-\tilde{x}_{i-1}\vee \tilde{y}_{j-1}$ is below $p_{ij,d}^{t,I\!B}$ (or equivalently, $\tilde{y}_n-\tilde{x}_{i-1}\vee \tilde{y}_{j-1}$ is below $p_{ij,s}^{t,I\!B}$), the aggregate demand level (or aggregate supply level) is already below the target level before $i$ matches with $j$. The one-step-ahead policy will not match $i$ and $j$, and neither any demand/supply of lower quality.

\item If $p_{ij,d}^{t,I\!B}\le \tilde{x}_m-\tilde{x}_i\wedge \tilde{y}_j$ (or equivalently $p_{ij,s}^{t,I\!B}\le \tilde{y}_n-\tilde{x}_i\wedge\tilde{y}_j$), either type $i$ demand is depleted before total demand reduces to $p_{ij,d}^{t,I\!B}$ or type $j$ supply is depleted before total supply reduces to $p_{ij,s}^{t,I\!B}$. The one-step-ahead policy will match $i$ with $j$ to the full extent.

\item If $\tilde{x}_m-\tilde{x}_i\wedge  \tilde{y}_j<p_{ij,d}^{t,I\!B}\le \tilde{x}_m-\tilde{x}_{i-1}\vee \tilde{y}_{j-1}$ and  $\tilde{y}_n-\tilde{x}_i\wedge\tilde{y}_j<p_{ij,s}^{t,I\!B}\le \tilde{y}_n-\tilde{x}_{i-1}\vee \tilde{y}_{j-1}$, the policy matches $i$ and $j$ until the total demand reduces to $p_{ij,d}^{t,I\!B}$ and total supply to $p_{ij,s}^{t,I\!B}$.
\end{itemize}


\subsection{Impatient demand and patient supply}

Consider the case with $\alpha=0$ and $\beta=1$, in which unmatched demand is lost at the end of each period and unmatched supply is fully carried to the next period. (The results hold for the case with an arbitrary supply carry-over rate $\beta\in(0,1)$.) 
Since demand does not carry over to the next period, we do not have the demand state $\mathbf{x}$ in the dynamic program and the optimal policy does not depend on $\mathbf{x}$.
Following similar analysis, we can show that Proposition \ref{prop:vertical-mono} remains true.

\subsection{Non-additive reward structure}
\label{subsec:nonadd}

Our results can be generalized to account for non-additive reward structures.
Instead of adopting the reward function $r_{ij}^t=r_{id}^t+r_{js}^t$, we consider the following assumption.

\begin{assumption}\label{as:vert-non-add}
    (i)  The unit matching reward $r_{ij}^t$ is decreasing in $i$ and $j$;

    (ii)  For $i=1,\ldots,m-1$ and $j=1,\ldots,n-1$, 
    $r_{ij}^t-r_{i+1,j}^t\ge \alpha (r_{ij''}^{t+1}-r_{i+1,j''}^{t+1})$ holds for all $j''\in\mathcal{S}$
    and $r_{ij}^t-r_{i,j+1}^t\ge \beta(r_{i''j}^{t+1}-r_{i'',j+1}^{t+1})$ holds for all $i''\in\mathcal{D}$;

    (iii)  $r_{ij}^t$ is supermodular with respect to $i$ and $j$,
    i.e., $r_{ij}^t-r_{i,j+1}^t \ge r_{i+1,j}^t-r_{i+1,j+1}^t$ for $i=1,\ldots,m-1$ and $j=1,\ldots,n-1$.
\end{assumption}

Under Assumption \ref{as:vert-non-add}, Propositions \ref{prop:vert-top-down} and  \ref{prop:vertical-mono} remain true.

In particular, part (i) of Assumption \ref{as:vert-non-add} suggests that a demand/supply type with a smaller index has a higher quality level.
Part (ii) generalizes Assumption \ref{as:vert-quality} and says that the quality difference between a high type demand (supply) and a low type demand (supply) is decreasing over time. 
Part (iii) of Assumption \ref{as:vert-non-add}  further ensures the condition (\ref{eqn:ddd}) in Definition \ref{def:Neighbor-dominance-strong}.
By Theorem \ref{cor:Neighboring-priority-strong}, there is an optimal policy that respects $\mathcal{M}_s$.
Thus, a high-quality demand/supply type has higher priority over a lower-quality type in such a policy. Consequently, Proposition \ref{prop:vert-top-down} holds under Assumption \ref{as:vert-non-add}.

Following similar analysis (see the proof of Proposition \ref{prop:vertical-mono} in Appendix \ref{App:proofs}), we can show that under Assumption \ref{as:vert-non-add} (i)--(iii), Proposition \ref{prop:vertical-mono} remains true for both the case with patient demand and supply and the case with patient supply but impatient demand.


To illustrate the conditions in Assumption \ref{as:vert-non-add}, we consider the reward structure $r_{ij}^t=a_i^t+b_j^t+\gamma a_i^t b_j^t$.
Parts (i) and (iii) of the assumption is satisfied when $a_i^t$ decreases in $i$ and $b_j^t$ decreases in $j$.
One can verify that part (ii) is satisfied if and only if:
\begin{align*}
    & \min_{i,i'\in\mathcal{D}} \frac{1+\gamma a_i^t}{ 1+ \gamma a_{i'}^t}\ge \beta\cdot \frac{b_j^{t+1}-b_{j+1}^{t+1}}{b_j^t-b_{j+1}^t}
    \text{ and }  \min_{j,j'\in\mathcal{S}} \frac{1+\gamma b_j^t}{ 1+ \gamma b_{j'}^t}\ge \alpha\cdot \frac{a_i^{t+1}-a_{i+1}^{t+1}}{a_i^t-a_{i+1}^t}.
\end{align*}

The above conditions are met if both $\beta\cdot \frac{b_j^{t+1}-b_{j+1}^{t+1}}{b_j^t-b_{j+1}^t}$ and $\alpha\cdot \frac{a_i^{t+1}-a_{i+1}^{t+1}}{a_i^t-a_{i+1}^t}$ are smaller than 1, and the parameter $\gamma$ is sufficiently small (i.e., the additive component of the reward is sufficiently more significant than the multiplicative component).
 

\section{Conclusion}

We consider a stochastic and dynamic matching framework with heterogeneous demand and supply types in the discrete-time setting.
We generalize the Monge sequence to establish conditions (which we call the (weak) modified Monge conditions) to prioritize demand-supply pairs optimally. 
Two reward structures satisfy the modified Monge condition for all neighboring pairs. In the unidirectionally horizontal reward structure, ``distance" determines priority,
and in the vertical reward structure, ``quality" determines priority.
Under both reward structures, the optimal matching proceeds along the priority structure, and when it comes to the matching between a specific pair, the optimal policy has a match-down-to threshold structure. 
This structural property of ``priority and thresholds" is a generalization of priority structures seen in the balanced and deterministic transportation problems, and the threshold-type policies seen in the inventory management (such as base-stock levels) and quantity-based revenue management (such as protection levels).

The proposed framework generalizes many classic problems. For example, we generalize inventory rationing problems and dynamic capacity allocation models with upgrading, by allowing for multiple exogenous supply streams and arbitrary substitution.
It also lays out a foundation for further research in the area of dynamic matching at the operational level. 
For example, one can consider joint pricing and matching decisions and competition among platforms.

\begin{APPENDIX}{}
	\subsection*{Weak compatibility}


\begin{definition} (Weak Compatibility)
	We say that an optimal matching policy $\{\mathbf{Q}^{t} \}_{t=1,\ldots,T}$ weakly respects $\succ_{\mathcal{M}}$, if  
(i) for all $(i,j)\succ_{\mathcal{M}} (i',j)$ and all $t=1,\ldots,T$, either $q_{i'j}^{t}=0$ or $u_i^{t}=0$;
    (ii) for all $(i,j)\succ_{\mathcal{M}}(i,j')$ and all $t=1,\ldots,T$, either $q_{ij'}^{t}=0$ or $v_j^{t}=0$.
    \label{def:mm-policy}
\end{definition}

If a policy weakly respects $\succ_{\mathcal{M}}$, then under this policy, a dominant pair of demand and supply types has higher ``priority" than a dominated pair (with the dominance relation determined by $\succ_{\mathcal{M}}$) in the following sense: 
If $(i,j)\succ_{\mathcal{M}}(i',j)$, unless there is no remaining type $i$ demand (i.e., $u_i^{t}=0$; in other words, it is impossible to further match $i$ with $j$), the optimal policy would not match type $i'$ demand with type $j$ supply.
We can verify that the weak compatibility can be inferred by the compatibility defined in Definition \ref{def:mms-policy}.

The following result studies the structure of the optimal policy when only the weak modified Monge condition is satisfied.


\begin{theorem}
    There exists an optimal matching policy that weakly respects $\succ_{\mathcal{M}}$. 
    \label{thm:Priority-neighboring}
\end{theorem}





Even though Theorem \ref{thm:Priority-neighboring} suggests that the optimal policy would not match $(i',j)$ unless type $i$ demand runs out (provided that $(i,j)\succ_{\mathcal{M}} (i',j)$), it does not necessarily mean that a dominant pair should always be matched \emph{before} a dominated pair; see the following example. 

\begin{example}\label{exa:priority-M}
Suppose that $(i,j)\succ_{\mathcal{M}} (i',j)$ and $(i,j)\succ_{\mathcal{M}} (i,j')$.
In period $t$, with the state $(\mathbf{x},\mathbf{y})$ such that $x_i=y_j=1$, the matching quantities $q_{ij}^{t*}=0$, $q_{i'j}^{t*}=q_{ij'}^{t*}=1$ can be optimal and weakly respect $\succ_{\mathcal{M}}$,  consistent with Theorem \ref{thm:Priority-neighboring}. (With $q_{i'j}^{t*}=q_{ij'}^{t*}=1$, there is no remaining type $i$ demand or type $j$ supply.) However, the matching over $(i,j)$ is not prioritized over $(i',j)$ and $(i,j')$.
%
\Halmos 
\end{example}

In contrast, the partial relation  $\succ_{\mathcal{M}_s}$ strengthens $\succ_{\mathcal{M}}$ to provide a sufficient condition that indeed ensures that a dominant pair is prioritized over a dominated pair in the optimal matching policy (Theorem \ref{cor:Neighboring-priority-strong}).


\begin{remark} \label{rem:robust}
  The partial relation $\succ_{\mathcal{M}}$ is not only sufficient but also \emph{robustly necessary} for Theorem \ref{thm:Priority-neighboring}. That is, if the conditions associated with $\succ_{\mathcal{M}}$ are not satisfied, one can construct an instance of demand and supply distributions such that the optimal policy does not satisfy the property in Definition \ref{def:mm-policy}. 
	Likewise, the partial relation $\succ_{\mathcal{M}_s}$ is robustly necessary for Theorem \ref{cor:Neighboring-priority-strong}. 
\end{remark}

\proof{Proof of Remark \ref{rem:robust}.}
To see the first claim, we show that there exists an instance in which the statement ``either $q_{ij'}^t=0$ or $u_j^t=0$'' does not hold under the optimal policy, in either of the following situations: 
(i) $r_{ij}^t < r_{ij'}^t$; (ii) $r_{ij}^t -r_{ij'}^t < \beta (r_{i''j}^{t+1}-r_{i''j'}^{t+1})$ for some $i''\in\mathcal{D}$.


For i), we consider the state with $x_i=1$, $y_j=y_{j'}=1$ and $x_{i'''}=y_{j'''}=0$ for all other $i'''\in\mathcal{D}$ and $j'''\in\mathcal{S}$ in period $t$.
Moreover, there is no new arrival of demand or supply from period $t$ to the end of the time horizon.
In this case, the problem reduces to a single-period problem. Given that $r_{ij}^t < r_{ij'}^t$, it is optimal to set $q_{ij}^{t}=0$ and $q_{ij'}^{t}=1$ in period $t$. 
Since $v_j^t=1>0$ under those matching quantities, the statement does not hold. 

For ii), we again consider the state with $x_i=1$, $y_j=y_{j'}=1$ and $x_{i'''}=y_{j'''}=0$ for all other $i'''\in\mathcal{D}$ and $j'''\in\mathcal{S}$ in period $t$.
In period $t+1$, 1 unit of type $i''$ demand arrives. Other than that, there is no new arrival from period $t$ to the end of time horizon.
The problem reduces to a two-period (i.e., period $t$ and period $t+1$) deterministic problem.
We either match $i$ with $j$ in period $t$ and match $i''$ with $j'$ in period $t+1$ (note that  quantity $\beta$ of type $j'$ supply carries over to period $t+1$),
or match $i$ with $j'$ in period $t$ and match $i''$ with $j$ in period $t+1$.
Since 
$r_{ij}^t -r_{ij'}^t < \beta (r_{i''j}^{t+1}-r_{i''j'}^{t+1})$, the latter option leads to a higher total reward.
Thus under the optimal policy, $q_{ij'}^t=1>0$ and $v_j^t=1>0$, implying that the statement does not hold.

To see the second claim, suppose that $r_{ij}^t+r_{i'j'}^t\ge r_{i'j}^t+r_{ij'}^t$ does not hold for some $(i,j)\succ_{\mathcal{M}}(i,j')$ and $(i,j)\succ_{\mathcal{M}}$.
Consider the state with $x_i=y_j=x_{i'}=y_{j'}=1$ in period $t$, and there is no new arrival of demand or supply from period $t$ to the end of the time horizon.
Then, the optimal decision in period $t$ is to match $i$ with $j'$ for 1 unit and match $i'$ with $j$ for 1 unit, which does not prioritize $(i,j)$ over $(i',j)$ or $(i,j')$.
\Halmos
\endproof

The weaker partial order $\succ_{\mathcal{M}}$ requires less on the reward structure than $\succ_{\mathcal{M}_s}$. But the stronger partial order $\succ_{\mathcal{M}_s}$ guarantees more, i.e., the optimal policy satisfies a priority hierarchy which can turn matching into a sequential procedure based on the partial order. In the paper, we have focused on problems in which the reward structure indeed leads to the strong partial relation $\mathcal{M}_s$ and obtain the optimal priority matching structure. But even if we only have the weak partial order $\succ_{\mathcal{M}}$ (i.e., the additional condition (\ref{eqn:ddd}) is not satisfied), we may adopt, the best policy within the sequential procedure based on the partial order, as a heuristic. In Online Appendix \ref{sec_appendix:weak--strong}, we show that any policy that weakly respects a partial order is only different from another policy that respects the same partial order by a single-period transportation problem.

\subsection*{Computation of protection levels in Heuristic \ref{heu:horizontal}}


We consolidate all demand types $i'\neq i$ into a single artificial demand type $i^c$, and all supply types $j'\neq j$ into a single artificial supply type $j^c$.
See Figure \ref{fig:heu-consol} for a demonstration.
\begin{figure}[htb]
    \centering
    \includegraphics[width=0.2\textwidth]{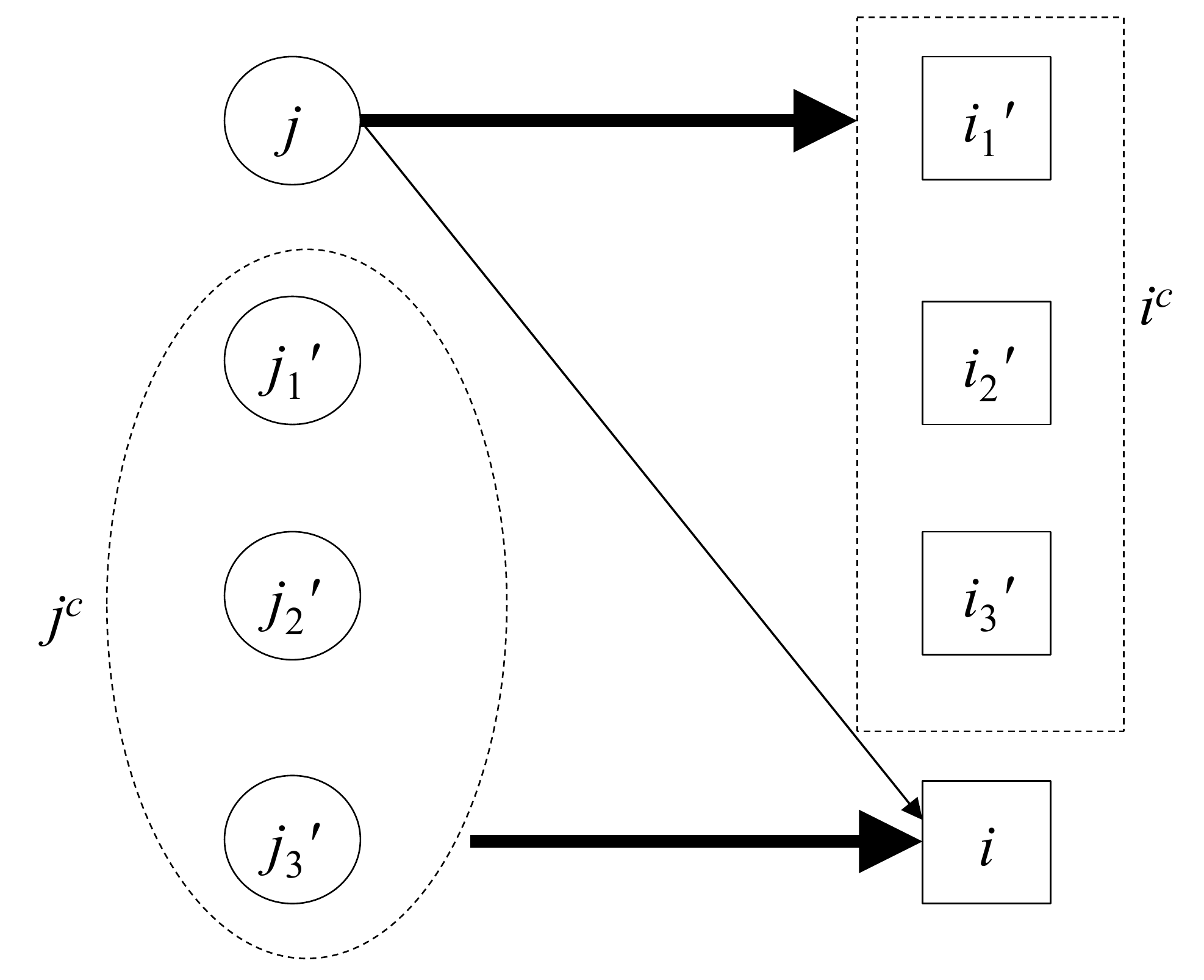}
    \caption{Consolidation of types.}
    \label{fig:heu-consol}
\end{figure}

Since each pair $(i,j')$ has a higher priority than $(i,j)$ and so does each pair $(i',j)$, 
both $(i,j^c)$ and $(i^c,j)$ are prioritized over $(i,j)$.
While the total demand of all type $i' (\neq i)$ demand combined is $\sum_{i'\neq i} D_{i'}^\tau$ in period $\tau$, only part of the amount are available for supply type $j$.
This is because a supply type $j''$ located between $i$ and $j$ on the line segment $C$ has a higher priority than $j$ in matching with some $i'$, and thus may consume some of type $i'$ demand.
Likewise, the combined demand of all type $j'$ supply may not be fully available for type $i$ demand.
Next, we define the ``effective'' demand of type $i^c$ and ``effective'' supply of type $j^c$.

For a period $\tau$ ($t\le \tau\le T$), let type $i'$ demand be realized as $d_{i'}^\tau$ and type $j'$ supply be realized as $s_{j'}$, for all types $i'$ and $j'$ involved in the subproblem P$(i,j)$. 
We also assume that those are all the available demand and supply in period $\tau$, for $i'\neq i$ and $j'\neq j$.
The pair $(i,j)$ has the lowest priory among all other pairs in the subproblem P$(i,j)$, because the distance between $i$ and $j$ is longer than the distance between any other pair. 
We perform greedy matching on the subproblem in period $\tau$ according to the distance-based priority structure, but terminate right before matching $i$ with $j$.
We consider the remaining demand and supply, denoted by $\tilde{d}_{i'}^\tau$ and $\tilde{s}_{j'}^\tau$ ($i'\neq i, j'\neq j$), as the type $i'$ demand and type $j'$ supply available to match with type $j$ supply and type $i$ demand, respectively.
We let $\tilde{D}_{i'}^\tau$ and $\tilde{S}_{j'}^\tau$ be the corresponding random variables to the realizations  $\tilde{d}_{i'}^\tau$ and $\tilde{s}_{j'}^\tau$, respectively.
Although it is difficult to derive the exact distributions of those random variables, we may sample a number of realizations and fit the distributions of $\tilde{D}_{i'}^\tau$ and $\tilde{S}_{j'}^\tau$ accordingly.

We then consolidate all the demand type $i'\neq i$ as a single demand type $i^c$ and all the supply type $j'\neq j$ as a single supply type $j^c$.
In each period $\tau$ ($t\le \tau\le T$), let $D_{i^c}^\tau:=\sum_{i'\neq i} \tilde{D}_{i'}^\tau$ be the quantity of type $i^c$ demand to arrive in period $\tau$, and $S_{j^c}^\tau:=\sum_{j'\neq j} \tilde{S}_{j'}^\tau$ the quantity of type $j^c$ supply to arrive in period $\tau$.
We assume that the unit matching reward between $i$ and the artificial state $j^c$ is a weighted average among all  of the state $j'\neq j$ being consolidated, and the same for the unit matching reward between $i^c$ and $j$. 
More specifically, we define the unit rewards $\tilde{r}_{ij}$ between the four types, $i,i^c,j,j^c$, as follows.
Let $\tilde{r}_{ij}^\tau:=r_{ij}^\tau$, $\tilde{r}_{ij^c}^\tau:= [\sum_{j'\neq j} r_{ij'}^\tau E\tilde{S}_{j'}^\tau]/ ES_{j^c}$, $\tilde{r}_{i^cj}^\tau := [\sum_{i'\neq i} r_{i'j}^\tau E\tilde{D}_{i'}^\tau]/ ED_{i^c}$ and $\tilde{r}_{i^cj^c}^\tau:=0$.
We can readily verify that $i$ and $j^c$ are perfect matches for each other, and so are $i^c$ and $j$.
The subproblem P$(i,j)$ then reduces to a $2\times 2$ model, from which we can obtain the protection levels.

\end{APPENDIX}

\bibliographystyle{ormsv080}
\bibliography{ref_proposal}

\newpage
 \setcounter{page}{1}
\renewcommand\thesection{\Alph{section}}
\renewcommand{\theequation}{\thesection.\arabic{equation}}
\setcounter{equation}{0}
\setcounter{section}{0}
\setcounter{lemma}{0}
\setcounter{example}{0}
\setcounter{theorem}{0}
\setcounter{corollary}{0}
\setcounter{proposition}{0}
\renewcommand\thelemma{\thesection.\arabic{lemma}}
\renewcommand\theproposition{\thesection.\arabic{proposition}}
\renewcommand\thetheorem{\Alph{section}}
\renewcommand\thecorollary{\Alph{section}}
\setcounter{definition}{0}

   
\begin{center}
       {\bf\large Online Appendices to ``Dynamic Type Matching"}
\end{center}

\section{Proofs}
\label{App:proofs}

\proof{Proof of Theorem \ref{thm:Priority-neighboring}.}
To prove Theorem \ref{thm:Priority-neighboring}, we need two lemmas.

\begin{lemma}
     The following statements hold for all periods.
\begin{enumerate}[(i)]
\item 
	For any $x_i>0$ and any $\varepsilon \in [0, x_i]$, there exists $( \lambda_1^\tau,\dots, \lambda_m^\tau )\ge\bm{0}$ for $\tau=t,\ldots,T$, such that $\sum_{\tau=t}^{T}\alpha^{-(\tau-t)}\sum_{j'=1}^{m}\lambda_{j'}^\tau\le \varepsilon$ and $V_t(\mathbf{x}-\varepsilon\mathbf{e}_i^n+\varepsilon\mathbf{e}_{i'}^n,\mathbf{y})-V_t(\mathbf{x},\mathbf{y})\ge -\sum_{\tau=t}^{T}\sum_{j'=1}^{m}\lambda_{j'}^\tau(r_{ij'}^\tau-r_{i'j'}^\tau)$.
\item 
	For any $y_j>0$ and any $\varepsilon\in[0, y_j]$, there exists $( \xi_1^\tau, \dots, \xi_n^\tau)\ge\bm{0}$ for $\tau=t,\ldots,T$, such that $\sum_{\tau=t}^{T}\beta^{-(\tau-t)}\sum_{i'=1}^{n}\xi_{i'}^\tau\le \varepsilon$ and $V_t(\mathbf{x},\mathbf{y}-\varepsilon\mathbf{e}_j^m+\varepsilon\mathbf{e}_{j'}^m)-V_t(\mathbf{x},\mathbf{y})\ge -\sum_{\tau=t}^{T}\sum_{i'=1}^{n}\xi_{i'}^\tau(r_{i' j}^\tau-r_{i' j'}^\tau)$.
    \end{enumerate}
    \label{lem:ineq-0}
\end{lemma}
\proof{Proof of Lemma \ref{lem:ineq-0}.}
We only need to prove part (i)  The proof of part (ii) is symmetric to part (i) 

The proof is based on induction.
The result holds for $t=T+1$. Because $V_{T+1}(\mathbf{x},\mathbf{y})\equiv 0$, we can simply set $\lambda_j^T$ to zero.
Suppose that it holds for period $t+1$.

Now consider period $t$.
Let $\hat{\mathbf{Q}}\in\arg\max_{\mathbf{Q}}H_t(\mathbf{Q},\mathbf{x},\mathbf{y})$ be an optimal decision in period $t$ under the state $(\mathbf{x},\mathbf{y})$ in period $t$.
We will construct a decision $\bar{\mathbf{Q}}$ that is feasible under the state $(\mathbf{x}-\varepsilon\mathbf{e}_i^n+\varepsilon\mathbf{e}_{i'}^n,\mathbf{y})$.

Under the new state $(\mathbf{x}-\varepsilon\mathbf{e}_i^n+\varepsilon\mathbf{e}_{i'}^n,\mathbf{y})$, the capacity of $i$ is reduced by $\varepsilon$ compared with the original state $(\mathbf{x},\mathbf{y})$.
We need to adjust the matching decision $\hat{\mathbf{Q}}$ accordingly to make it feasible for the new state. 
In particular, we reduce the matching quantity $\hat{q}_{ij}$ by $\mu_j$ for $j=1,\ldots,n$,
where the nonnegative numbers $\mu_1,\ldots,\mu_n$ are defined as follows. 
\begin{align*}
    \mu_j=& \min\{\hat{q}_{ij},(\varepsilon-\sum_{j'=1}^{j-1}\hat{q}_{ij'})^+\}, \text{ for } j=1,\ldots,n. 
\end{align*}
If $\sum_{j'=1}^{k-1}\hat{q}_{ij'}<\varepsilon\le \sum_{j'=1}^{k}\hat{q}_{ij'}$ for some $1\le k\le n$, then one can verify that $\mu_j=\hat{q}_{ij}$ for $j=1,\ldots,k-1$, $\mu_k= \varepsilon-\sum_{j'=1}^{k}\hat{q}_{ij'}$ and $\mu_j=0$ for $j=k+1,\ldots,n$.
In this case, $\sum_{j=1}^{n}\mu_j=\varepsilon$, and thus $\sum_{j=1}^{n} (\hat{q}_{ij}-\mu_j)=\sum_{j=1}^{n} \hat{q}_{ij} - \sum_{j=1}^{n} \mu_j=\sum_{j=1}^{n} \hat{q}_{ij}-\varepsilon \le x_i-\varepsilon$.

If $\varepsilon>\sum_{j'=1}^{n}\hat{q}_{ij'}$, then $\mu_j=\hat{q}_{ij}$ for all $j=1,\ldots,n$.
Therefore, we reduce the matching quantity $\hat{q}_{ij}$ starting from $j=1$, until either a total reduction $\varepsilon$ is reached or all quantities $\hat{q}_{ij}$ ($j=1,\ldots,n$) are reduced to 0.
In this case, $\sum_{j=1}^{n} (\hat{q}_{ij}-\mu_j)=0\le x_i-\varepsilon$.

On the other hand, under the new state $(\mathbf{x}-\varepsilon\mathbf{e}_i^n+\varepsilon\mathbf{e}_{i'}^n,\mathbf{y})$, the capacity of $i'$ is increased by $\varepsilon$.
This allows us to increase the matching quantity $\hat{q}_{i'j}$ by $\mu_j$ for all $j=1,\ldots,n$.

We define
\begin{align*}
    \bar{\mathbf{Q}} = \hat{\mathbf{Q}} - \sum_{j=1}^{n}  \mu_{j}\mathbf{e}_{ij}^{m\times n} + \sum_{j=1}^{n}  \mu_{j}\mathbf{e}_{i'j}^{m\times n},
\end{align*} 
which is feasible for the state $(\mathbf{x}-\varepsilon\mathbf{e}_i^n+\varepsilon\mathbf{e}_{i'}^n,\mathbf{y})$. To see this, we have
\begin{align*}
    \mathbf{1}_m  \bar{\mathbf{Q}} = &\mathbf{1}_m \hat{\mathbf{Q}} - \sum_{j=1}^{n}  \mu_{j} \mathbf{1}_m \mathbf{e}_{ij}^{m\times n} + \sum_{j=1}^{n}  \mu_{j} \mathbf{1}_m\mathbf{e}_{i'j}^{m\times n} = \mathbf{1}_m \hat{\mathbf{Q}} - \sum_{j=1}^{n}  \mu_{j}  \mathbf{e}_{j}^{n} + \sum_{j=1}^{n}  \mu_{j} \mathbf{e}_{j}^{n} = \mathbf{1}_m \hat{\mathbf{Q}} \le \mathbf{y}.
\end{align*}
Also,
\begin{align*}
    \bar{\mathbf{Q}} \mathbf{1}_n^{\tt{T}}=& \hat{\mathbf{Q}} \mathbf{1}_n^{\tt{T}}- \sum_{j=1}^{n}  \mu_{j}\mathbf{e}_{ij}^{m\times n} \mathbf{1}_n^{\tt{T}}+ \sum_{j=1}^{n}  \mu_{j}\mathbf{e}_{i'j}^{m\times n}\mathbf{1}_n^{\tt{T}}
    =\hat{\mathbf{Q}} \mathbf{1}_n^{\tt{T}}- \sum_{j=1}^{n}  \mu_{j}(\mathbf{e}_{i}^{m})^{\tt{T}} + \sum_{j=1}^{n}  \mu_{j}(\mathbf{e}_{i'}^{m})^{\tt{T}}.
\end{align*}
It follows that $(\bar{\mathbf{Q}} \mathbf{1}_n^{\tt{T}})_i=\sum_{j=1}^{n} \hat{q}_{ij} - \sum_{j=1}^{n} \mu_j\le x_i-\varepsilon$, $(\bar{\mathbf{Q}} \mathbf{1}_n^{\tt{T}})_{i'}=\sum_{j=1}^{n} \hat{q}_{i'j} + \sum_{j=1}^{n} \mu_j\le \sum_{j=1}^{n} \hat{q}_{i'j}+\varepsilon\le x_{i'}+\varepsilon$ and $(\bar{\mathbf{Q}} \mathbf{1}_n^{\tt{T}})_{i''}=\sum_{j=1}^{n}\hat{q}_{i''j} \le x_{i''} $ for all $i''\neq
i,i'$.
Thus, $\bar{\mathbf{Q}} \mathbf{1}_n^{\tt{T}}\le (\mathbf{x}-\varepsilon\mathbf{e}_i^n+\varepsilon\mathbf{e}_{i'}^n)^{\tt{T}}$.

Therefore, $\bar{\mathbf{Q}}$ is a feasible decision for the state $(\mathbf{x}-\varepsilon\mathbf{e}_i^n+\varepsilon\mathbf{e}_{i'}^n,\mathbf{y})$.
Under the decision $\bar{\mathbf{x}}$, the total reward received in period $t$ is 
\begin{align*}
    \mathbf{R}^t\circ \bar{\mathbf{Q}}= \mathbf{R}^t\circ (\hat{\mathbf{Q}}- \sum_{j=1}^{n}  \mu_{j}\mathbf{e}_{ij}^{m\times n} + \sum_{j=1}^{n}  \mu_{j}\mathbf{e}_{i'j}^{m\times n})
    =\mathbf{R}^t\circ \hat{\mathbf{Q}} - \sum_{j=1}^{n}  \mu_{j}r_{ij}^{t} + \sum_{j=1}^{n}  \mu_{j}r_{i'j}^{t}.
\end{align*}
The post-matching levels in period $t$ are
\begin{align*}
    \bar{\mathbf{u}} =& \mathbf{x} - \varepsilon\mathbf{e}_i^m +\varepsilon\mathbf{e}_{i'}^m-\mathbf{1}_n \bar{\mathbf{Q}}^{\tt{T}}\\
    =&  \mathbf{x} - \varepsilon\mathbf{e}_i^m +\varepsilon\mathbf{e}_{i'}^m - \mathbf{1}_n \hat{\mathbf{Q}}^{\tt{T}}  + \sum_{j=1}^{n}  \mu_{j}\mathbf{e}_{i}^{m} - \sum_{j=1}^{n}  \mu_{j}\mathbf{e}_{i'}^{m}\\
    =& \hat{\mathbf{u}} -  (\varepsilon-\sum_{j=1}^{n}  \mu_{j})\mathbf{e}_{i}^{m} +(\varepsilon- \sum_{j=1}^{n}  \mu_{j})\mathbf{e}_{i'}^{m},\\
    \bar{\mathbf{v}} =& \mathbf{y} - \mathbf{1}_m \bar{\mathbf{Q}} = \mathbf{y} - \mathbf{1}_m \hat{\mathbf{Q}} = \hat{\mathbf{v}}.
\end{align*}

Consequently,
\begin{align*}
    &V_t(\mathbf{x}-\varepsilon\mathbf{e}_i^n+\varepsilon\mathbf{e}_{i'}^n,\mathbf{y})-V_t(\mathbf{x},\mathbf{y})\\
    \ge& H_t(\bar{\mathbf{Q}},\mathbf{x}-\varepsilon\mathbf{e}_i^n+\varepsilon\mathbf{e}_{i'}^n,\mathbf{y}) - H_t(\hat{\mathbf{Q}},\mathbf{x},\mathbf{y})\\
    =& - \sum_{j=1}^{n}  \mu_{j} (r_{ij}^{t}-r_{i'j}^{t}) \\
    &+ EV_{t+1}(\alpha \hat{\mathbf{u}} -  \alpha(\varepsilon-\sum_{j=1}^{n}  \mu_{j})\mathbf{e}_{i}^{m} + \alpha(\varepsilon-\sum_{j=1}^{n}  \mu_{j})\mathbf{e}_{i'}^{m}+\mathbf{D}^{t+1},\beta \hat{\mathbf{v}}+\mathbf{S}^{t+1}) - EV_{t+1}(\alpha \hat{\mathbf{u}} +\mathbf{D}^{t+1},\beta \hat{\mathbf{v}}+\mathbf{S}^{t+1}).
\end{align*}

By the induction hypothesis, for each realization of $\mathbf{D}^{t+1}$ and $\mathbf{S}^{t+1}$, there exists $(\Lambda_{1}^\tau,\ldots,\Lambda_{n}^\tau)$ for $\tau=t+1,\ldots,T$ such that $\sum_{\tau=t+1}^{T} \alpha^{-(\tau-t-1)}\sum_{j=1}^{n} \Lambda_{j}^\tau\le \alpha(\varepsilon-\sum_{j=1}^{n}  \mu_{j})$ and  
\begin{align*}
   & V_{t+1}(\alpha \hat{\mathbf{u}} -  \alpha(\varepsilon-\sum_{j=1}^{n}  \mu_{j})\mathbf{e}_{i}^{m} + \alpha(\varepsilon-\sum_{j=1}^{n}  \mu_{j})\mathbf{e}_{i'}^{m}+\mathbf{D}^{t+1},\beta \hat{\mathbf{v}}+\mathbf{S}^{t+1}) - V_{t+1}(\alpha \hat{\mathbf{u}} +\mathbf{D}^{t+1},\beta \hat{\mathbf{v}}+\mathbf{S}^{t+1})\\
   \ge& -\sum_{\tau=t+1}^{T} \sum_{j=1}^{n} \Lambda_{j}^\tau (r_{ij}^\tau - r_{i'j}^\tau).
\end{align*}

Note that $\Lambda_{j}^\tau$ is a random variable due to its possible dependency on the random vectors $\mathbf{D}^{t+1}$ and $\mathbf{S}^{t+1}$.

It then follows that
\begin{align*}
    &V_t(\mathbf{x}-\varepsilon\mathbf{e}_i^n+\varepsilon\mathbf{e}_{i'}^n,\mathbf{y})-V_t(\mathbf{x},\mathbf{y})\ge - \sum_{j=1}^{n}  \mu_{j} (r_{ij}^{t}-r_{i'j}^{t}) - \sum_{\tau=t+1}^{T} \sum_{j=1}^{n} E\Lambda_{j}^\tau \cdot (r_{ij}^\tau - r_{i'j}^\tau).
\end{align*}

Since $\sum_{\tau=t+1}^{T} \alpha^{-(\tau-t-1)}\sum_{j=1}^{n} \Lambda_{j}^\tau\le \alpha(\varepsilon-\sum_{j=1}^{n}  \mu_{j})$, 
we have $\sum_{j=1}^{n}  \mu_{j} + \sum_{\tau=t+1}^{T} \alpha^{-(\tau-t)} \sum_{j=1}^{n} \Lambda_{j}^\tau \le \varepsilon $.
Let $\lambda_j^t=\mu_j$ for all $j=1,\ldots,n$, and $\lambda_j^\tau=E\Lambda_{j}^\tau$ for all $j=1,\ldots,n$ and $\tau=t+1,\ldots,T$. The proof of the lemma is then completed. \Halmos
\endproof

\begin{lemma}
    (i)  Suppose that $(i,j)\succ_{\mathcal{M}} (i',j)$.
    Then, transferring matching quantity from $(i',j)$ to $(i,j)$ weakly improves the total expected reward, i.e., $H_t(\mathbf{Q}+\varepsilon\mathbf{e}_{ij}^{m\times n}-\varepsilon\mathbf{e}_{i'j}^{m\times n},\mathbf{x},\mathbf{y})\ge H_t(\mathbf{Q},\mathbf{x},\mathbf{y})$, if $\mathbf{Q}+\varepsilon\mathbf{e}_{ij}^{m\times n}-\varepsilon\mathbf{e}_{i'j}^{m\times n}$ is a feasible decision under the state
    $(\mathbf{x},\mathbf{y})$.

    (ii)  Similarly, if $(i,j)\succ_{\mathcal{M}} (i,j')$, then $H_t(\mathbf{Q}+\varepsilon\mathbf{e}_{ij}^{m\times n}-\varepsilon\mathbf{e}_{ij'}^{m\times n},\mathbf{x},\mathbf{y})\ge H_t(\mathbf{Q},\mathbf{x},\mathbf{y})$.
    \label{lem:transfer}
\end{lemma}
\proof{Proof of Lemma \ref{lem:transfer}.}
    We prove part (i) only since part (ii) can be proved analogously. 
    The post-matching levels for using $\mathbf{Q}+\varepsilon\mathbf{e}_{ij}^{m\times n}-\varepsilon\mathbf{e}_{ij'}^{m\times n}$ are
    \begin{align*}
        \bar{\mathbf{u}} = & \mathbf{x} - \mathbf{1}_n  (\mathbf{Q}+\varepsilon\mathbf{e}_{ij}^{m\times n}-\varepsilon\mathbf{e}_{i'j}^{m\times n})^{\tt{T}} 
        = \mathbf{x} - \mathbf{1}_n  \mathbf{Q} - \varepsilon \mathbf{e}_i^m + \varepsilon \mathbf{e}_{i'}^m = \mathbf{u} - \varepsilon \mathbf{e}_i^m + \varepsilon \mathbf{e}_{i'}^m,\\
        \bar{\mathbf{v}} = & \mathbf{y} - \mathbf{1}_m  (\mathbf{Q}+\varepsilon\mathbf{e}_{ij}^{m\times n}-\varepsilon\mathbf{e}_{i'j}^{m\times n}) = \mathbf{y} - \mathbf{1}_m  \mathbf{Q} = \mathbf{v},
    \end{align*}
    where $(\mathbf{u},\mathbf{v})$ are the post-matching levels by using the decision $\mathbf{Q}$ in period $t$.
    Then,
\begin{align}
  & H_t(\mathbf{Q}+\varepsilon \mathbf{e}_{ij}^{m\times n} - \varepsilon \mathbf{e}_{ij'}^{m\times n},\mathbf{x},\mathbf{y}) - H_t(\mathbf{Q},\mathbf{x},\mathbf{y}) \nonumber \\
  =& \varepsilon( r_{ij}^t - r_{ij'}^t ) +   EV_{t+1}(\alpha \bar{\mathbf{u}}+\mathbf{D}^{t+1}, \beta \bar{\mathbf{v}} + \mathbf{S}^{t+1})
  -  EV_{t+1}(\alpha \mathbf{u}+\mathbf{D}^{t+1},\beta \mathbf{v} +\mathbf{S}^{t+1}).\label{eqn:H-diff}
\end{align}

By Lemma \ref{lem:ineq-0}, there exists $(\Lambda_{1}^\tau,\ldots,\Lambda_{n}^\tau)$ for $\tau=t+1,\ldots,T$ such that $\sum_{\tau=t+1}^{T} \alpha^{-(\tau-t-1)}\sum_{j'=1}^{n} \Lambda_{j'}^\tau\le \alpha\varepsilon$ and  
\begin{align*}
   V_{t+1}(\alpha \bar{\mathbf{u}}+\mathbf{D}^{t+1}, \beta \bar{\mathbf{v}} + \mathbf{S}^{t+1})
  -  V_{t+1}(\alpha \mathbf{u}+\mathbf{D}^{t+1},\beta \mathbf{v} +\mathbf{S}^{t+1})\ge -\sum_{\tau=t+1}^{T} \sum_{j'=1}^{n} \Lambda_{j'}^\tau (r_{ij'}^\tau - r_{i'j'}^\tau).
\end{align*}
Note that $\Lambda_j^\tau$ is a random variable since it may depend on $\mathbf{D}^{t+1}$ and $\mathbf{S}^{t+1}$.

$(i,j)\succ_{\mathcal{M}} (i',j)$ implies that $r_{ij'}^{\tau}-r_{i'j'}^{\tau} \le \alpha^{-(\tau-t)} (r_{ij}^{t}-r_{i'j}^{t})$ for all $j'\in\mathcal{S}$ and $\tau=t+1,\ldots,T$.
Thus,
\begin{align}
  & V_{t+1}(\alpha \bar{\mathbf{u}}+\mathbf{D}^{t+1}, \beta \bar{\mathbf{v}} + \mathbf{S}^{t+1})
  -  V_{t+1}(\alpha \mathbf{u}+\mathbf{D}^{t+1},\beta \mathbf{v} +\mathbf{S}^{t+1})\nonumber\\
  \ge&  -\sum_{\tau=t+1}^{T} \sum_{j'=1}^{n} \Lambda_{j'}^\tau (r_{ij'}^\tau - r_{i'j'}^\tau)\nonumber\\
	\ge& -(r_{ij'}^t - r_{i'j'}^t)\sum_{\tau=t+1}^{T} \alpha^{-(\tau-t-1)}\sum_{j'=1}^{n} \Lambda_{j'}^\tau \nonumber\\
    \ge & -\frac{1}{\alpha} (r_{ij}^{t}-r_{i'j}^{t}) \times \alpha\varepsilon \nonumber\\
    =& - (r_{ij}^{t}-r_{i'j}^{t})\varepsilon. \label{eqn:V-bound}
\end{align}
Combining (\ref{eqn:H-diff}) and (\ref{eqn:V-bound}), we have $H_t(\mathbf{Q}+\varepsilon \mathbf{e}_{ij}^{m\times n} - \varepsilon \mathbf{e}_{ij'}^{m\times n},\mathbf{x},\mathbf{y}) \ge H_t(\mathbf{Q},\mathbf{x},\mathbf{y})$.  \Halmos
\endproof

We now proceed to prove Theorem \ref{thm:Priority-neighboring}.

Let $\mathbf{Q}^*$ be an optimal decision in period $t$ under the state $(\mathbf{x},\mathbf{y})$.
Suppose that $\mathbf{Q}^*$ does not satisfy the desired. 
For $(i,j)\succ_{\mathcal{M}} (i',j)$, we will transfer a quantity $\varepsilon:=\min\left\{ q_{i'j}^{t*},u_i^{t*} \right\}$ from $(i',j)$ to $(i,j)$ if both $q_{i'j}^{t*}$ and $u_i^{t*}$ are positive.
Similarly, for $(i,j)\succ_{\mathcal{M}} (i,j')$, we will transfer a quantity $\varepsilon:=\min\left\{ q_{ij'}^{t*},v_j^{t*} \right\}$ from $(i,j')$ to $(i,j)$ if both $q_{ij'}^{t*}$ and $v_j^{t*}$ are positive.
If multiple transfers are possible, we choose the one that yields the greatest transferring quantity first.

After each transfer, we obtain a different matching decision that is feasible and weak dominates $\mathbf{Q}^*$ according to Lemma \ref{lem:transfer}. 
Let $\mathbf{Q}^k$ be the matching decision after $k$ transfers.
Then, $\mathbf{Q}^k$ is also optimal since it weakly dominates the optimal decision $\mathbf{Q}^*$
Since the transfers are unidirectional (i.e., from a dominated pair of demand and supply to a dominant one), either the transferring procedure ends in finite steps (say, in $K$ steps) or the quantity transferred converges to zero.
In the former case, we obtain an optimal decision $\mathbf{Q}^K$ that satisfies the desired properties.
In the latter case, $\mathbf{Q}^k$ converges to an optimal decision $\mathbf{Q}^K$ that satisfies the desired properties.
\Halmos
\endproof

\bigskip

\proof{Proof of Theorem \ref{cor:Neighboring-priority-strong}.}
The proof is analogous to Theorem \ref{thm:Priority-neighboring}.
In addition to the transfers we considered in the proof of Theorem \ref{thm:Priority-neighboring}, we also consider the following type of transfers.

Let $\mathbf{Q}^k$ be a feasible decision in period $t$ under the state $(\mathcal{x},\mathcal{y})$, and $(\mathcal{u},\mathcal{v})$ be the corresponding post-matching levels. Consider $(i,j)$, $(i',j)$ and $(i,j')$ such that $(i,j)\succ_{\mathcal{M}_s} (i',j)$ and $(i,j)\succ_{\mathcal{M}_s} (i,j')$.
If both $q_{i'j}^{t*}$ and $q_{ij'}^{t*}$ are positive,
we construct another feasible decision $\mathbf{Q}^{k+1}:=\mathbf{Q}^k+\varepsilon \mathbf{e}_{ij}^{m\times n} + \varepsilon\mathbf{e}_{ij}^{m\times n} - \varepsilon\mathbf{e}_{i'j}^{m\times n} - \varepsilon\mathbf{e}_{ij'}^{m\times n}$, where $\varepsilon:=\min\left\{ q_{i'j}^{t*},  q_{ij'}^{t*}\right\}$.

By repeatedly applying the transfers (as we did in the proof of Theorem \ref{thm:Priority-neighboring}), we eventually reaches a feasible decision $\mathbf{Q}$ either in finite steps or in the limit, with the following properties:

(i)  For $(i,j)\succ_{\mathcal{M}_s}(i',j)$, either $q_{i'j}^t=0$ or $u_i^t=0$.

(ii)  For $(i,j)\succ_{\mathcal{M}_s}(i,j')$, either $q_{ij'}^t=0$ or $v_j^t=0$.

(iii)  For $(i,j)\succ_{\mathcal{M}_s} (i,j')$ and $(i,j)\succ_{\mathcal{M}_s} (i',j)$, either $q_{i'j}^t=0$ or $q_{ij'}^t=0$.

We now show that the above properties lead to the corollary.
Suppose to the contrary that $q_{i'j}^t>0$ and $a_i^t>0$ for some $(i,j)\succ_{\mathcal{M}_s}(i',j)$.
According to property (i) above, $u_i^t=0$, which implies that $x_i-\sum_{j'': (i,j'')\in\mathcal{B}_{ij,L}} q_{ij''}^t - \sum_{j'': (i,j'')\notin\mathcal{B}_{ij,L}} q_{ij''}^t=0$.
By the definition of $a_i^t$, we know that $x_i- \sum_{j'': (i,j'')\notin\mathcal{B}_{ij,L}} q_{ij''}^t>0$.
It follows that $\sum_{j'': (i,j'')\in\mathcal{B}_{ij,L}} q_{ij''}^t$.
This further implies that there exists $(i,j')$ such that $(i,j)\succ_{\mathcal{M}_s}(i,j')$ and $q_{ij'}^t>0$. However, this contradicts property (iii) 
Thus, either $q_{i'j}^t=0$ and $a_i^t=0$ if $(i,j)\succ_{\mathcal{M}_s}(i',j)$.

Similarly, we can show that either $q_{ij'}^t=0$ and $b_j^t=0$ if $(i,j)\succ_{\mathcal{M}_s}(i,j')$.

Finally, if we started with an optimal decision, $\mathbf{Q}$ is also optimal since it weakly dominates the initial optimal decision. \Halmos
\endproof

\bigskip

\proof{Proof of Proposition \ref{prop:greedy-optimal}.}
We show that greedy matching between $i$ and $j$ is optimal by induction.
It is easy to verify that greedy matching between $i$ and $j$ is optimal in the final period $T$.
Suppose that it is also optimal in period $t+1$.

Let $\mathbf{Q}$ be an optimal decision in period $t$ under the state $(\mathbf{x},\mathbf{y})$.
Suppose that $q_{ij}<\min\left\{ x_i,y_j \right\}$.

We first show that by using $\mathbf{Q}$ both the post-matching levels $u_i^t$ and $v_j^t$ are positive if $(i,j)$ dominates all its neighboring pairs by $\succ_{\mathcal{M}_s}$.
To prove that, let us suppose to the contrary that $u_i^t=0$.
Since $q_{ij}<x_i$, there is a pair $(i,j')$ such that $q_{ij'}>0$.
Following Theorem \ref{cor:Neighboring-priority-strong}, we have $b_j^t=0$.
According to the definition of $b_j^t$, we have $v_j^t\le b_j^t=0$, implying that $v_j^t=0$.
However, $v_j^t=0$ but $q_{ij}<y_j$ implies that there is some $(i',j)$ such that $q_{i'j}>0$.
Consequently, $b_j^t\ge v_j^t+q_{i'j}>0$, which contradicts $b_j^t=0$.
Similarly, we can show that $v_i^t=0$ also leads to contradiction.
Thus, both $u_i^t$ and $v_j^t$ are positive.

We now show that increasing the matching quantity between $i$ and $j$ by $\varepsilon:=\min\left\{ u_i^t,v_j^t \right\}$ does not hurt the optimality of $\mathbf{Q}$.
Increasing the matching quantity between $i$ and $j$ by $\varepsilon$ will increase the matching reward in period $t$ by $r_{ij}^t \varepsilon$, but decrease both the post-matching levels of $i$ and $j$ by $\varepsilon$. In other words, we have
\begin{align*}
   & H_t(\mathbf{Q}+\varepsilon\mathbf{e}_{ij}^{m\times n},\mathbf{x},\mathbf{y})
    -H_t(\mathbf{Q}+\varepsilon\mathbf{e}_{ij}^{m\times n},\mathbf{x},\mathbf{y})\\
    = & r_{ij}^t \varepsilon + 
    EV_{t+1}(\alpha\mathbf{u}^t -\alpha \varepsilon\mathbf{e}_i^m + \mathbf{D}^{t+1},\beta \mathbf{v}^t -\beta\varepsilon\mathbf{e}_j^n  +\mathbf{S}^{t+1})
    -EV_{t+1}(\alpha\mathbf{u}^t + \mathbf{D}^{t+1},\beta \mathbf{v}^t+\mathbf{S}^{t+1}).
\end{align*}

Let us consider the case $\beta\ge \alpha$ without loss of generality.
We have
\begin{align*}
    V_{t+1}(\alpha\mathbf{u}^t+\mathbf{D}^{t+1}-\alpha\varepsilon\mathbf{e}_i^m,\beta \mathbf{v}^t+\mathbf{S}^{t+1}-\beta \varepsilon\mathbf{e}_j^n) =& V_{t+1}(\alpha\mathbf{u}^t+\mathbf{D}^{t+1} + (\beta-\alpha)\varepsilon\mathbf{e}_i^m,\beta \mathbf{v}^t+\mathbf{S}^{t+1} ) - \beta\varepsilon r_{ij}^{t+1} \\ 
    \ge & V_{t+1}(\alpha\mathbf{u}^t+\mathbf{D}^{t+1},\beta \mathbf{v}^t+\mathbf{S}^{t+1} ) - \beta\varepsilon r_{ij}^{t+1}
\end{align*}
 where the equality is because of the greedy matching of pair $(i,j)$ for the subsequent periods,
 and the inequality holds because $V_{t+1}$ is increasing in the state vector.
 Therefore, 
 \begin{align*}
     H_t(\mathbf{Q}+\varepsilon\mathbf{e}_{ij}^{m\times n},\mathbf{x},\mathbf{y})
    -H_t(\mathbf{Q}+\varepsilon\mathbf{e}_{ij}^{m\times n},\mathbf{x},\mathbf{y})
    \ge & (r_{ij}^t - \beta r_{ij}^{t+1})\varepsilon\ge 0.
 \end{align*}
 Therefore, we can always weakly improve $\mathbf{Q}$ if it does not greedily match $i$ with $j$.
 \Halmos
\endproof

\bigskip
\proof{Proof of Corollary \ref{cor:euclidean-dist}.}
$R_t$ decreasing in $t$ implies that $r_{ii}^t=R^t$ is decreasing in $t$.
By the triangle inequality, $\text{dist}_{i''j'}- \text{dist}_{i''i}\le \text{dist}_{ij'}$ for any location $i$, $i''$ and $j'$.
Thus $r_{i''i}^{t+1}-r_{i''j'}^{t+1}=\gamma_{t+1}(\dist_{i''j'}-\dist_{i''i})\le \gamma_{t+1} \text{dist}_{ij'}\le \gamma_t \dist_{ij'} = R_t-r_{ij'}^t=r_{ii}^t-r_{ij'}^t $, for any $i''\in\mathcal{D}$ and $j'\in\mathcal{S}$. This shows that $(i,i)\succ_{\mathcal{M}} (i,j')$ for any $j'\in\mathcal{S}$.
By symmetry, $(i,i)\succ_{\mathcal{M}} (i',j)$ for any $i'\in\mathcal{D}$.
Moreover, for all $i'\in\mathcal{D}$ and $j'\in\mathcal{S}$,
$r_{i'i}^{t}-r_{i'j'}^{t}=\gamma_{t+1}(\dist_{i'j'}-\dist_{i'i})\le \gamma_{t} \text{dist}_{ij'}=r_{ii}^t-r_{ij'}^t$. This is equivalent to $r_{ii}^t+r_{i'j'}^{t}\ge r_{i'i}^{t}+r_{ij'}^t$ for all $i'\in\mathcal{D}$ and $j'\in\mathcal{S}$.
This shows that $\succ_{\mathcal{M}}$ can be defined as $\mathcal{M}_s$ according to Definition \ref{def:Neighbor-dominance-strong}. \Halmos
\endproof
\bigskip

\proof{Proof of Proposition \ref{prop:opt-horizontal22}.}
    We focus on the matching in round 2, and only consider the case with $z_1\ge 0$ and $z_1\ge 0$ (the case with $z_1<0$ and $z_2<0$ is symmetric).
    
    Using the formulation (\ref{eqn:horizontal22-alt1})--(\ref{eqn:horizontal22-alt2}) in Appendix \ref{app:horizontal22}, the optimal matching quantity solves $\max_{q\in M(\mathbf{z})}J_t(q,\mathbf{z})$.
    Let us use $p_d:=z_1-q$ and $p_s=z_2-q$ as decision variables in place of $q$.
    Then, $p_d=p_s+z_1-z_2=p_s+I\!B$.
    Since both $p_d$ and $p_s$ need to be nonnegative, the feasible range of $p_s$ is $I\!B^-\le p_s\le z_2$.

    We rewrite $J_t(q,\mathbf{z})$ as a function of $p_d$ and $p_s$.
    Since $0\le q\le \min\{z_1,z_2\}$, we have
    \begin{align*}
        J_t(q,\mathbf{z})=& r_{12}^t q + r_{11}^{t+1} E\min\left\{ \alpha(z_1-q)+D_1^{t+1},S_1^{t+1} \right\} 
        + r_{22}^{t+1} E\min\left\{D_2^{t+1}, \beta (z_2-q)+S_2^{t+1}  \right\} \\
        &+ EU_{t+1}(\alpha(z_1-q)+D_1^{t+1}-S_1^{t+1},\beta(z_2-q) + S_2^{t+1}-D_2^{t+1})\\
        =& r_{12}^t (z_2-p_s) + r_{11}^{t+1} E\min\left\{ \alpha(p_s+I\!B)+D_1^{t+1},S_1^{t+1} \right\} 
        + r_{22}^{t+1} E\min\left\{D_2^{t+1}, \beta p_s+S_2^{t+1}  \right\} \\
        &+ EU_{t+1}(\alpha(p_s+I\!B)+D_1^{t+1}-S_1^{t+1},\beta p_s + S_2^{t+1}-D_2^{t+1}),
    \end{align*}
    which depends on $I\!B$, $p_s$ and also linearly on $z_2$.
    We write $J_t(q,\mathbf{z})=r_{12}^t z_2 + \check{J}_t(p_s,I\!B)$.

    It is easy to see that $\check{J}_t$ is concave in $p_s$ (by Lemma \ref{lem:concavity22} $U_{t+1}$ is concave).
    Let $p_{s,+}^{t,I\!B}\in\arg\max_{p_s\ge I\!B^-} \check{J}_t(p_s,I\!B)$.
    Given the constraint $I\!B^-\le p_s\le z_2$, the optimal decision in terms of $p_s$ is $\min\left\{ z_2, p_{s,+}^{t,I\!B}\right\}$.

    Let us denote $p_{d,+}^{t,I\!B}:=I\!B+p_{s,+}^{t,I\!B}$. The optimal decision in terms of $p_d$ is 
    $I\!B+\min\left\{ z_2, p_{s,+}^{t,I\!B}\right\}=\min\left\{ I\!B+z_2,I\!B+ p_{s,+}^{t,I\!B}\right\}=\min\left\{ z_1,p_{d,+}^{t,I\!B} \right\}$. \Halmos
\endproof

\bigskip

\proof{Proof of Proposition \ref{prop:protection-monotone}.}
Let $J_t$ and $U_t$ be defined as in (\ref{eqn:horizontal22-alt1})--(\ref{eqn:horizontal22-alt1}). To prove the proposition, we present two lemmas.

\begin{lemma}
   Suppose that $(r_{22}^t-r_{12}^t)-(r_{22}^{t+1}-r_{12}^{t+1})$ decreases in $t$. For $0\le \varepsilon'\le \varepsilon$ and $t=1,\ldots,T$, $(r_{11}^{t+1}+r_{22}^{t+1}-r_{12}^t)\varepsilon' + (r_{22}^t-r_{22}^{t+1})\varepsilon + \tilde{U}_{t+1}(z_1+\varepsilon',z_2+\varepsilon')-\tilde{U}_{t+1}(z_1,z_2)\ge 0$, where $r_{ij}^{T+1}$ is defined as zero for all $i$ and $j$.
    \label{lem:diffU-hor22}
\end{lemma}
\proof{Proof of Lemma \ref{lem:diffU-hor22}.}
We prove this lemma by induction.

Since $U_{T+1}(\mathbf{z})\equiv0$ and $r_{ij}^{T+1}\equiv0$ for all $i,j$, we have $\tilde{U}_{T+1}(\mathbf{z})=0$. Then,
\begin{align*}
(r_{11}^{T+1}+r_{22}^{T+1}-r_{12}^T)\varepsilon' + (r_{22}^T-r_{22}^{T+1})\varepsilon
=& -r_{12}^T \varepsilon' + r_{22}^T \varepsilon\ge 0.
\end{align*}

Suppose that $(r_{11}^{t+2}+r_{22}^{t+2}-r_{12}^{t+1})\varepsilon' + (r_{22}^{t+1}-r_{22}^{t+2})\varepsilon + \tilde{U}_{t+2}(z_1+\varepsilon',z_2+\varepsilon')-\tilde{U}_{t+2}(z_1,z_2)\ge 0$ for any $\mathbf{z}=(z_1,z_2)$ and $0\le \varepsilon'\le \varepsilon$. 
We will show that the same inequality holds for $t+1$, i.e., $(r_{11}^{t+1}+r_{22}^{t+1}-r_{12}^t)\varepsilon' + (r_{22}^t-r_{22}^{t+1})\varepsilon + \tilde{U}_{t+1}(z_1+\varepsilon',z_2+\varepsilon')-\tilde{U}_{t+1}(z_1,z_2)\ge 0$ for any $\mathbf{z}=(z_1,z_2)$ and $0\le \varepsilon'\le \varepsilon$.

let $\hat{q}\in\arg\max_{q\ge 0} \tilde{J}_{t+1}(q,\mathbf{z})$.
We consider the following cases.

\emph{Case 1:} $\mathbf{z}\in\mathbb{R}_{++}^2$.
It is easy to see that $\hat{q}+\varepsilon'$ is a feasible matching quantity between type 1 demand and type 2 supply under the state $(z_1+\varepsilon',z_2+\varepsilon')$.
Thus,
\begin{align*}
    \tilde{U}_{t+1}(z_1+\varepsilon',z_2+\varepsilon')-\tilde{U}_{t+1}(z_1,z_2)
    \ge \tilde{J}_{t+1}(\hat{q}+\varepsilon',z_1+\varepsilon',z_2+\varepsilon') - \tilde{J}_{t+1}(\hat{q},z_1,z_2)=(-r_{11}^{t+1}-r_{22}^{t+1} + r_{12}^{t+1})\varepsilon'.
\end{align*}
It follows that
\begin{align*}
    & (r_{11}^{t+1}+r_{22}^{t+1}-r_{12}^t)\varepsilon' + (r_{22}^t-r_{22}^{t+1})\varepsilon + \tilde{U}_{t+1}(z_1+\varepsilon',z_2+\varepsilon')-\tilde{U}_{t+1}(z_1,z_2)\\
    \ge & (r_{11}^{t+1}+r_{22}^{t+1}-r_{12}^t)\varepsilon' + (r_{22}^t-r_{22}^{t+1})\varepsilon
     + (-r_{11}^{t+1}-r_{22}^{t+1} + r_{12}^{t+1})\varepsilon'\\
     =& -(r_{12}^t-r_{12}^{t+1}) \varepsilon' + (r_{22}^t-r_{22}^{t+1})\varepsilon\\
     \ge & [ (r_{22}^t-r_{22}^{t+1}) -(r_{12}^t-r_{12}^{t+1})] \varepsilon\\
     =& [(r_{22}^t-r_{12}^t) - (r_{22}^{t+1}-r_{12}^{t+1})]\ge 0.
\end{align*}

\emph{Case 2:} $\mathbf{z}\in\mathbb{R}_{--}^2$.
Let $\varepsilon''=\min\left\{ \hat{q},\varepsilon \right\}$.
It is easy to see that $\hat{q}+\varepsilon''$ is a feasible decision under the state $\mathbf{z}+\varepsilon' \mathbf{1}_2=(z_1+\varepsilon',z_2+\varepsilon')$. Then,
\begin{align*}
    &\tilde{U}_{t+1}(\mathbf{z}+\varepsilon' \mathbf{1}_2) - \tilde{U}_{t+1}(\mathbf{z})\\
    \ge & \tilde{J}_{t+1}(\hat{q}+\varepsilon'',\mathbf{z}+\varepsilon' \mathbf{1}) - \tilde{J}_{t+1}(\hat{q},\mathbf{z})\\
    =&  -r_{21}^{t+1} \varepsilon'' + E\tilde{U}_{t+2}(z_1-\hat{q}+\varepsilon'-\varepsilon''+D_1^{t+1}-S_1^{t+1},z_2-\hat{q}+\varepsilon'-\varepsilon''+D_2^{t+1}-S_2^{t+1})\\
    &- E\tilde{U}_{t+2}(z_1-\hat{q}+D_1^{t+2}-S_1^{t+2},z_2-\hat{q}+D_2^{t+2}-S_2^{t+2})\\
    \ge & -r_{21}^{t+1} \varepsilon'' -(r_{11}^{t+2}+r_{22}^{t+2}-r_{12}^{t+1})(\varepsilon'-\varepsilon'') - (r_{22}^{t+1}-r_{22}^{t+2})\varepsilon_0,
\end{align*}
for any $\varepsilon_0\ge \varepsilon'-\varepsilon''$.
We set $\varepsilon_0=\varepsilon'-\varepsilon''$.

It follows that 
\begin{align*}
  & (r_{11}^{t+1}+r_{22}^{t+1}-r_{12}^t)\varepsilon' + (r_{22}^t-r_{22}^{t+1})\varepsilon + \tilde{U}_{t+1}(z_1+\varepsilon',z_2+\varepsilon')-\tilde{U}_{t+1}(z_1,z_2)\\
  \ge &  (r_{11}^{t+1}+r_{22}^{t+1}-r_{12}^t)\varepsilon' + (r_{22}^t-r_{22}^{t+1})\varepsilon
       -r_{21}^{t+1} \varepsilon'' -(r_{11}^{t+2}+r_{22}^{t+2}-r_{12}^{t+1})(\varepsilon'-\varepsilon'') - (r_{22}^{t+1}-r_{22}^{t+2})\varepsilon_0\\
       =& (r_{11}^{t+1}+r_{22}^{t+1}-r_{12}^t)\varepsilon'+ (r_{22}^t-r_{22}^{t+1})\varepsilon
       -r_{21}^{t+1} \varepsilon'' -(r_{11}^{t+2}+r_{22}^{t+2}-r_{12}^{t+1})(\varepsilon'-\varepsilon'') - (r_{22}^{t+1}-r_{22}^{t+2})(\varepsilon'-\varepsilon'')\\
       =& (r_{11}^{t+1}+r_{22}^{t+1})\varepsilon' - (r_{12}^t-r_{12}^{t+1}) \varepsilon' +  (r_{22}^t-r_{22}^{t+1})\varepsilon -(r_{11}^{t+2}+r_{22}^{t+2})(\varepsilon'-\varepsilon'')-r_{21}^{t+1} \varepsilon''-r_{12}^{t+1} \varepsilon''- (r_{22}^{t+1}-r_{22}^{t+2})(\varepsilon'-\varepsilon'')\\
       =& (r_{11}^{t+1}+r_{22}^{t+1})\varepsilon' - (r_{12}^t-r_{12}^{t+1}) \varepsilon' +  (r_{22}^t-r_{22}^{t+1})\varepsilon -r_{21}^{t+1} \varepsilon''-r_{12}^{t+1} \varepsilon''- (r_{11}^{t+2}+r_{22}^{t+1})(\varepsilon'-\varepsilon'')\\
       =& (r_{11}^{t+1}-r_{11}^{t+2})(\varepsilon'-\varepsilon'')- (r_{12}^t-r_{12}^{t+1}) \varepsilon' +  (r_{22}^t-r_{22}^{t+1})\varepsilon  -r_{21}^{t+1} \varepsilon''-r_{12}^{t+1} \varepsilon''  + (r_{11}^{t+1}+r_{22}^{t+1})\varepsilon''\\
       \ge& [(r_{22}^t-r_{12}^t)-(r_{22}^{t+1}-r_{12}^{t+1})]\varepsilon' + (r_{11}^{t+1}+r_{22}^{t+1}-r_{21}^{t+1}-r_{12}^{t+1})\varepsilon''.
\end{align*}

It follows from Assumption \ref{as:horizontal22b} that $(r_{22}^t-r_{12}^t)-(r_{22}^{t+1}-r_{12}^{t+1})\ge 0$, and from Assumption \ref{as:reward-22} that $r_{11}^{t+1}+r_{22}^{t+1}-r_{21}^{t+1}-r_{12}^{t+1}\ge 0$.
Therefore, we have $(r_{11}^{t+1}+r_{22}^{t+1}-r_{12}^t)\varepsilon' + (r_{22}^t-r_{22}^{t+1})\varepsilon + \tilde{U}_{t+1}(z_1+\varepsilon',z_2+\varepsilon')-\tilde{U}_{t+1}(z_1,z_2)\ge 0$.

\emph{Case 3:} $\mathbf{z}\in\mathbb{R}_{+-}^2$.
In this case, we have
\begin{align*}
    &\tilde{U}_{t+1}(\mathbf{z}+\varepsilon' \mathbf{1}) - \tilde{U}_{t+1}(\mathbf{z})\\
    =& \tilde{J}_{t+1}(0,\mathbf{z}+\varepsilon' \mathbf{1}) - \tilde{J}_{t+1}(0,\mathbf{z})\\
    =& -r_{11}^{t+1} \varepsilon' + r_{11}^{t+2} \varepsilon' + E\tilde{U}_{t+2}(z_1+\varepsilon'+D_1^t-S_1^t,z_2+\varepsilon'+D_2^t-S_2^t) - E\tilde{U}_{t+2}(z_1+D_1^t-S_1^t,z_2+D_2^t-S_2^t)
\end{align*}
It follows that 
\begin{align*}
     &(r_{11}^{t+1}+r_{22}^{t+1}-r_{12}^t)\varepsilon' + (r_{22}^t-r_{22}^{t+1})\varepsilon + \tilde{U}_{t+1}(z_1+\varepsilon',z_2+\varepsilon')-\tilde{U}_{t+1}(z_1,z_2)\\
 =& (r_{11}^{t+2}+r_{22}^{t+1}-r_{12}^t)\varepsilon' + (r_{22}^t-r_{22}^{t+1})\varepsilon \\
 & -(r_{11}^{t+2}+r_{22}^{t+2}-r_{12}^{t+1})\varepsilon' - (r_{22}^{t+1}-r_{22}^{t+2})\varepsilon\\
&+(r_{11}^{t+2}+r_{22}^{t+2}-r_{12}^{t+1})\varepsilon' + (r_{22}^{t+1}-r_{22}^{t+2})\varepsilon\\
 &+ E\tilde{U}_{t+2}(z_1+\varepsilon'+D_1^t-S_1^t,z_2+\varepsilon'+D_2^t-S_2^t) - E\tilde{U}_{t+2}(z_1+D_1^t-S_1^t,z_2+D_2^t-S_2^t)\\
 =& [(r_{22}^{t+1}-r_{22}^{t+2}) -(r_{12}^t - r_{12}^{t+1})]\varepsilon' 
 + [(r_{22}^t-r_{22}^{t+1}) - (r_{22}^{t+1}-r_{22}^{t+2}) ]\varepsilon\\
 \ge& \left[ (r_{22}^{t+1}-r_{22}^{t+2}) -(r_{12}^t - r_{12}^{t+1}) + (r_{22}^t-r_{22}^{t+1}) - (r_{22}^{t+1}-r_{22}^{t+2})\right] \varepsilon'\\
 =& \left[ (r_{22}^t-r_{22}^{t+1}) - (r_{12}^t - r_{12}^{t+1}) \right]\varepsilon'\ge 0,
\end{align*}
where the first inequality holds because $\varepsilon'\le \varepsilon$ and $(r_{22}^t-r_{22}^{t+1}) - (r_{22}^{t+1}-r_{22}^{t+2})\ge 0$ (By Assumption \ref{as:horizontal22b}),
and the second inequality holds because of Assumption \ref{as:horizontal22b}.

\emph{Case 4:} $\mathbf{z}\in\mathbb{R}_{-+}^+$.
In this case, we have
\begin{align*}
	&\tilde{U}_{t+1}(\mathbf{z}+\varepsilon' \mathbf{1}) - \tilde{U}_t(\mathbf{z})\\
	=& \tilde{J}_{t+1}(0,\mathbf{z}+\varepsilon' \mathbf{1}) - \tilde{J}_{t+1}(0,\mathbf{z})\\
    =& -r_{22}^{t+1} \varepsilon' + r_{22}^{t+2} \varepsilon' + E\tilde{U}_{t+2}(z_1+\varepsilon'+D_1^{t+2}-S_1^{t+2},z_2+\varepsilon'+D_2^{t+2}-S_2^{t+2}) \\
	 &- E\tilde{U}_{t+2}(z_1+D_1^{t+2}-S_1^{t+2},z_2+D_2^{t+2}-S_2^{t+2}).
\end{align*}

It follows that 
\begin{align*}
	 &(r_{11}^{t+1}+r_{22}^{t+1}-r_{12}^t)\varepsilon' + (r_{22}^t-r_{22}^{t+1})\varepsilon + \tilde{U}_{t+1}(z_1+\varepsilon',z_2+\varepsilon')-\tilde{U}_{t+1}(z_1,z_2)\\
	 =&  (r_{11}^{t+1}+r_{22}^{t+1}-r_{12}^t)\varepsilon' + (r_{22}^t-r_{22}^{t+1})\varepsilon  
	      -r_{22}^{t+1} \varepsilon' + r_{22}^{t+2} \varepsilon' \\
				& + E\tilde{U}_{t+2}(z_1+\varepsilon'+D_1^t-S_1^{t+2},z_2+\varepsilon'+D_2^{t+1}-S_2^{t+1}) 
	 - E\tilde{U}_{t+2}(z_1+D_1^{t+1}-S_1^{t+1},z_2+D_2^{t+1}-S_2^{t+1})\\
	= & (r_{11}^{t+1}+r_{22}^{t+2}-r_{12}^t)\varepsilon' + (r_{22}^t-r_{22}^{t+1})\varepsilon\\ 
    & + E\tilde{U}_{t+2}(z_1+\varepsilon'+D_1^{t+2}-S_1^{t+2},z_2+\varepsilon'+D_2^{t+2}-S_2^{t+2}) 
	 - E\tilde{U}_{t+2}(z_1+D_1^{t+2}-S_1^{t+2},z_2+D_2^{t+2}-S_2^{t+2})\\
	= & (r_{11}^{t+1}+r_{22}^{t+2}-r_{12}^t)\varepsilon' + (r_{22}^t-r_{22}^{t+1})\varepsilon\\ 
	    & -(r_{11}^{t+2}+r_{22}^{t+2}-r_{12}^{t+1})\varepsilon' - (r_{22}^{t+1}-r_{22}^{t+2})\varepsilon\\
&+(r_{11}^{t+2}+r_{22}^{t+2}-r_{12}^{t+1})\varepsilon' + (r_{22}^{t+1}-r_{22}^{t+2})\varepsilon\\
& + E\tilde{U}_{t+2}(z_1+\varepsilon'+D_1^{t+2}-S_1^{t+2},z_2+\varepsilon'+D_2^{t+2}-S_2^{t+2}) 
	 - E\tilde{U}_{t+2}(z_1+D_1^{t+2}-S_1^{t+2},z_2+D_2^{t+2}-S_2^{t+2})\\
	 \ge&  (r_{11}^{t+1}+r_{22}^{t+2}-r_{12}^t)\varepsilon' + (r_{22}^t-r_{22}^{t+1})\varepsilon 
	     -(r_{11}^{t+2}+r_{22}^{t+2}-r_{12}^{t+1})\varepsilon' - (r_{22}^{t+1}-r_{22}^{t+2})\varepsilon\\
			 =& [(r_{11}^{t+1}-r_{11}^{t+2})-(r_{12}^t-r_{12}^{t+1}] \varepsilon' + [(r_{22}^t-r_{22}^{t+1})-(r_{22}^{t+1}-r_{22}^{t+2})]\varepsilon\\
			 \ge & [(r_{12}^{t+1}-r_{12}^{t+2})-(r_{12}^t-r_{12}^{t+1})] \varepsilon' + [(r_{22}^t-r_{22}^{t+1})-(r_{22}^{t+1}-r_{22}^{t+2})]\varepsilon\\
			 \ge & -[(r_{12}^t-r_{12}^{t+1})-(r_{12}^{t+1}-r_{12}^{t+2})] \varepsilon + [(r_{22}^t-r_{22}^{t+1})-(r_{22}^{t+1}-r_{22}^{t+2})]\varepsilon\\
             \ge& \left\{ [(r_{22}^t-r_{12}^t)-(r_{22}^{t+1}-r_{12}^{t+1})]
             -[(r_{22}^{t+1}-r_{12}^{t+1})-(r_{22}^{t+2}-r_{12}^{t+2})] \right\} \varepsilon.
\end{align*}
where the first inequality follows from the induction hypothesis, the second one from Assumption \ref{as:horizontal22b}, and the third one holds because $(r_{22}^t-r_{12}^t)-(r_{22}^{t+1}-r_{12}^{t+1})$ decreases in $t$.
The induction is completed.\Halmos
\endproof

\begin{lemma}
    Suppose that $\alpha=\beta=1$. The functions $\tilde{U}_t(\mathbf{z}):=-r_{11}^tz_1^+-r_{22}^tz_2^+ + U_t(\mathbf{z})$
    and $\tilde{J}_t(q,\mathbf{z}):=-r_{11}^tz_1^+-r_{22}^tz_2^++J_t(q,\mathbf{z})$ are $L^\natural$-concave with respect to all variables.
    \label{lem:L-concave22}
\end{lemma}
\proof{Proof of Lemma \ref{lem:L-concave22}.}
By using the equality $\min\left\{ a,b \right\}=a-(a-b)^+$, we can rewrite $J_t$ as follows.
\begin{align*}
   J_t(q,\mathbf{z}) = & r_{12}^t q^+ + r_{21}^t q^- + ED_1^t + ES_2^t \\
   &+  r_{11}^{t+1} E(z_1-q)^+ 
   -r_{11}^{t+1} E(z_1-q+D_1^{t+1}-S_1^{t+1})^+\\
   & + r_{22}^{t+1} E(z_2-q)^+ - r_{22}^{t+1} E(z_2-q+S_2^{t+1}-D_2^{t+1})^+\\
   & + EU_{t+1}(z_1-q+D_1^{t+1}-S_1^{t+1},z_2-q+S_2^{t+1}-D_2^{t+1}).
\end{align*}
Then, by definition, we have
\begin{align}
	\tilde{J}_t(q,\mathbf{z}) = & ED_1^t + ES_2^t -r_{11}^t z_1^+ - r_{22}^t z_2^+ +	r_{12}^t q^+ + r_{21}^t q^-  
    + r_{11}^{t+1}E(z_1-q)^+ + r_{22}^{t+1} E(z_2-q)^+ \nonumber \\
    & + E\tilde{U}_{t+1}(z_1-q+D_1^{t+1}-S_1^{t+1},z_2-q+S_2^{t+1}-D_2^{t+1}). \label{eqn:J-tilde-exp0}
\end{align}

We further rewrite $\tilde{J}_t(q,\mathbf{z})$ as follows.
\begin{align*}
	\tilde{J}_t(q,\mathbf{z}) =& ED_1^t + ES_2^t  -r_{11}^t z_1^+ - r_{22}^t z_2^+ +	(r_{12}^t+r_{21}^t) q^+  - r_{21}^t q   + r_{11}^{t+1}E(z_1-q)^+ + r_{22}^{t+1} E(z_2-q)^+ \\
    & + E\tilde{U}_{t+1}(z_1-q+D_1^{t+1}-S_1^{t+1},z_2-q+S_2^{t+1}-D_2^{t+1})\\
    =&  ED_1^t + ES_2^t -r_{11}^t [z_1^+-(z_1-q)^+] - r_{22}^t [z_2^+-(z_2-q)^+] + (r_{12}^t+r_{21}^t) q^+  - r_{21}^t q \\
    & - (r_{11}^t-r_{11}^{t+1}) (z_1-q)^+ - (r_{22}^t-r_{22}^{t+1}) (z_2-q)^+ \\
    &+ E\tilde{U}_{t+1}(z_1-q+D_1^{t+1}-S_1^{t+1},z_2-q+S_2^{t+1}-D_2^{t+1}) \\
    =& ED_1^t + ES_2^t -r_{11}^t q^+  - r_{22}^t q^+ + (r_{12}^t+r_{21}^t) q^+  - r_{21}^t q \\
    & - (r_{11}^t-r_{11}^{t+1}) (z_1-q)^+ - (r_{22}^t-r_{22}^{t+1}) (z_2-q)^+ \\
    & + E\tilde{U}_{t+1}(z_1-q+D_1^{t+1}-S_1^{t+1},z_2-q+S_2^{t+1}-D_2^{t+1})                          \\
    =& ED_1^t + ES_2^t -(r_{11}^t + r_{22}^t - r_{12}^t - r_{21}^t ) q^+ -r_{21}^t q \\
    & - (r_{11}^t-r_{11}^{t+1}) (z_1-q)^+ - (r_{22}^t-r_{22}^{t+1}) (z_2-q)^+\\
    & + E\tilde{U}_{t+1}(z_1-q+D_1^{t+1}-S_1^{t+1},z_2-q+S_2^{t+1}-D_2^{t+1}),
\end{align*}
where the first inequality is due to the fact $q^-=q^+-q$, and the third equality holds because for $q\in M(\mathbf{z})$, $z_1^+-q^+ = (z_1-q)^+$ and $z_2^+-q^+=(z_2-q)^+$.

We now prove the proposition by induction.
$\tilde{U}_T(\mathbf{z})\equiv 0$ is $L^\natural$-concave.
Suppose that $\tilde{U}_{t+1}(\mathbf{z})$ is $L^\natural$-concave.
To show that $\tilde{U}_{t}(\mathbf{z})$ is $L^\natural$-concave, we need to prove that $\tilde{U}_{t}(\mathbf{z}-\eta\mathbf{1}_2)$ is supermodular in $(\eta,\mathbf{z})$.

Given the conditions $r_{11}^t\ge r_{11}^{t+1}$ and $r_{22}^t\ge r_{22}^{t+1}$, the induction hypothesis and the concavity of $-(\cdot)^+$, it is easy to see that $\tilde{J}_t(q-\eta,\mathbf{z}-\eta\mathbf{1}_2)$ is supermodular in $(\eta,q,\mathbf{z})$.
This implies that $\tilde{J}_t(q,\mathbf{z})$ is $L^\natural$-concave in $(q,\mathbf{z})$, and thus it is also supermodular in $(q,\mathbf{z})$.
Since the set $\{(q,\mathbf{z})\mid q\in M(\mathbf{z})\}$ is a lattice (See Appendix \ref{app:horizontal22}), $\tilde{U}_t(\mathbf{z})=\max_{q\in M(\mathbf{z})}\tilde{J}_t(q,\mathbf{z})$ is supermodular.
As a result, the function $\tilde{U}_t(\mathbf{z}-\eta\mathbf{1}_2)$ is supermodular in $\mathbf{z}$.
To show that $\tilde{U}_t$ is $L^\natural$-concave, it suffices to show that $\tilde{U}(\mathbf{z}-\eta\mathbf{1}_2)$ has increasing difference in $(\eta,z_1)$ and in $(\eta,z_2)$ within the feasible region $\{(\mathbf{z},\eta)\mid \mathbf{z}-\eta\mathbf{1}_2\ge 0\}$.
In the followings, we show that this is true within four regions (i.e., $\mathbf{z}\in \mathbb{R}_{++}$, $\mathbf{z}\in \mathbb{R}_{--}$, $\mathbf{z}\in \mathbb{R}_{+-}$ and $\mathbf{z}\in \mathbb{R}_{-+}$), as well as across the four regions.

For $\mathbf{z}\in \mathbb{R}_{++}$, $M(\mathbf{z})=\left\{ \mathbf{z}\mid 0\le q\le \min\left\{ z_1,z_2 \right\} \right\}$,
and $\left\{ (q,\eta,\mathbf{z})\mid \mathbf{z}-\eta\mathbf{1}_2 \ge 0, q-\eta\in M(\mathbf{z}-\eta\mathbf{1}_2) \right\}=\left\{ (q,\eta,\mathbf{z})\mid \eta\le q \le \min\left\{ z_1,z_2 \right\}  \right\}$.
The latter is a lattice. 
Then, $\tilde{U}_t(\mathbf{z}-\eta\mathbf{1}_2)=\max_{q-\eta\in M(\mathbf{z}-\eta\mathbf{1}_2)}\tilde{J}_t(q-\eta,\mathbf{z}-\eta\mathbf{1}_2)$ is supermodular in $(\eta,\mathbf{z})$ for $\eta\le \min\left\{ z_1,z_2 \right\}$.
This implies that $\tilde{U}_t(\mathbf{z}-\eta\mathbf{1})$ has increasing differences in $(\eta,z_1)$ and in $(\eta,z_2)$ for $\eta\le \min\left\{ z_1,z_2 \right\}$.

For $\mathbf{z}\in \mathbb{R}_{--}$, $M(\mathbf{z})=\left\{ \mathbf{z}\mid \max\left\{ z_1,z_2 \right\} \le q\le 0 \right\}$.
$\tilde{U}_t(\mathbf{z}-\eta\mathbf{1})=\max_{\max\left\{ z_1,z_2 \right\}\le q\le \eta} \tilde{J}_t(q-\eta,\mathbf{z}-\eta\mathbf{1})$ is supermodular in $(\eta,\mathbf{z})$ for $\eta\ge \max\left\{ z_1,z_2 \right\}$ because $\left\{ (q,\eta,\mathbf{z}) \mid \max\left\{ z_1,z_2 \right\}\le q\le \eta \right\}$ is a lattice.
Thus, $\tilde{U}_t(\mathbf{z}-\eta\mathbf{1}_2)$ has increasing differences in $(\eta,z_1)$ and in $(\eta,z_2)$ for $\eta\le \max\left\{ z_1,z_2 \right\}$.

For $\mathbf{z}\in \mathbb{R}_{+-}$ or $\mathbf{z}\in \mathbb{R}_{-+}$, we have $M(\mathbf{z})=\left\{ 0 \right\}$.
It is easy to verify that $\tilde{U}_t(\mathbf{z}-\eta\mathbf{1})=\tilde{J}(0,\mathbf{z}-\eta\mathbf{1})$ is supermodular in $\{(\eta,\mathbf{z})\mid \mathbf{z}-\eta \mathbf{1}_2\in \mathbb{R}_{+-}\}$, and in $\{(\eta,\mathbf{z})\mid \mathbf{z}-\eta \mathbf{1}_2\in \mathbb{R}_{-+}\}$.
Thus, $\tilde{U}_t(\mathbf{z}-\eta\mathbf{1}_2)$ has increasing differences in $(\eta,z_1)$ and in $(\eta,z_2)$ for $\mathbf{z}-\eta \mathbf{1}_2\in \mathbb{R}_{+-}$ and for $\mathbf{z}-\eta \mathbf{1}_2\in \mathbb{R}_{-+}$.

It remains to show that $\tilde{U}(\mathbf{z}-\eta\mathbf{1}_2)$ has increasing differences in $(\eta,z_1)$ and in $(\eta,z_2)$ across the 4 regions.
In the followings, we focus on the difference $\tilde{U}_t(\mathbf{z}+\varepsilon\mathbf{1}_2)-\tilde{U}_t(\mathbf{z})$ across the boundary between $\mathbb{R}^2_{++}$ and $\mathbb{R}^2_{+-}$.
The same property across the other boundaries can be proved similarly.
More specifically, we will prove the following inequality holds for sufficiently small $\varepsilon>0$.
\begin{align*}
    \tilde{U}_t(z_1,0) - \tilde{U}_t(z_1-\varepsilon,-\varepsilon)\ge 
     - \tilde{U}_t(z_1,\varepsilon) - \tilde{U}_t(z_1-\varepsilon,0),
\end{align*}
which implies that $\tilde{U}(\mathbf{z}-\eta\mathbf{1}_2)$ has increasing differences in $(\eta,z_2)$ across across the boundary between $\mathbb{R}^2_{++}$ and $\mathbb{R}^2_{+-}$. (The increasing difference property with respect to $(\eta,z_1)$ can be proved similarly.)

Let $z_1>0$ and $\hat{q}\in\arg\max_{q\in M(z_1,\varepsilon)}$. Also, let $\varepsilon'=\varepsilon-\hat{q}$. Then, by using the expression of $\tilde{J}_t$ given in (\ref{eqn:J-tilde-exp0}), we have
\begin{align*}
    & \tilde{U}_t(z_1,\varepsilon) - \tilde{U}_t(z_1-\varepsilon,0)\\
    =& \tilde{J}_t(\hat{q},z_1,\varepsilon) - \tilde{J}_t(0,z_1-\varepsilon,0)\\
    =& -r_{11}^t z_1 - r_{22}^t \varepsilon +r_{12}^t \hat{q} + r_{11}^{t+1} (z_1-\hat{q}) + r_{22}^{t+1}(\varepsilon-\hat{q})+ E\tilde{U}_{t+1}(z_1-\hat{q}+D_1^{t+1}-S_1^{t+1},\varepsilon-\hat{q} + S_2^{t+1}-D_2^{t+1})\\
    & -[-r_{11}^t (z_1-\varepsilon) + r_{11}^{t+1} (z_1-\varepsilon) + E\tilde{U}_{t+1}(z_1-\varepsilon+D_1^{t+1}-S_1^{t+1},S_2^{t+1}-D_2^{t+1})]\\
    =& -r_{11}^t\varepsilon -r_{22}^t \varepsilon +r_{12}^t \hat{q} + r_{11}^{t+1} (\varepsilon-\hat{q}) + r_{22}^{t+1} (\varepsilon-\hat{q}) \\
    & + E\tilde{U}_{t+1}(z_1-\hat{q}+D_1^{t+1}-S_1^{t+1},\varepsilon -\hat{q} + S_2^{t+1}-D_2^{t+1})-E\tilde{U}_{t+1}(z_1-\varepsilon+D_1^{t+1}-S_1^{t+1},S_2^{t+1}-D_2^{t+1})\\
    =&-r_{11}^t\varepsilon -r_{22}^t \varepsilon +r_{12}^t(\varepsilon-\varepsilon') + r_{11}^{t+1} \varepsilon' + r_{22}^{t+1} \varepsilon' 
    +E\tilde{U}_{t+1}(z_1-\varepsilon+\varepsilon'+D_1^{t+1}-S_1^{t+1},\varepsilon' + S_2^{t+1}-D_2^{t+1})\\
    & -E\tilde{U}_{t+1}(z_1-\varepsilon+D_1^{t+1}-S_1^{t+1},S_2^{t+1}-D_2^{t+1}),
\end{align*}


Also, for $z_1>0$, 
\begin{align*}
    &\tilde{U}_t(z_1,0) - \tilde{U}_t(z_1-\varepsilon,-\varepsilon)\\
    =&
    \tilde{J}_t(0,z_1,0) - \tilde{J}_t(0,z_1-\varepsilon,-\varepsilon)\\
    =& -r_{11}^t z_1 + r_{11}^{t+1} z_1 + E\tilde{U}_{t+1}(z_1+D_1^{t+1}-S_1^{t+1},S_2^{t+1}-D_2^{t+1})\\
    & -[ -r_{11}^t (z_1-\varepsilon) + r_{11}^{t+1} (z_1-\varepsilon)
    + E\tilde{U}_{t+1}(z_1-\varepsilon+D_1^{t+1}-S_1^{t+1},-\varepsilon+S_2^{t+1}-D_2^{t+1})]\\
    =& -(r_{11}^t-r_{11}^{t+1})\varepsilon\\
    &+ E\tilde{U}_{t+1}(z_1+D_1^{t+1}-S_1^{t+1},S_2^{t+1}-D_2^{t+1}) - E\tilde{U}_{t+1}(z_1-\varepsilon+D_1^{t+1}-S_1^{t+1},-\varepsilon+S_2^{t+1}-D_2^{t+1})\\
    =& -(r_{11}^t-r_{11}^{t+1})\varepsilon\\
    &+ E\tilde{U}_{t+1}(z_1+D_1^{t+1}-S_1^{t+1},S_2^{t+1}-D_2^{t+1}) - E\tilde{U}_{t+1}(z_1-\varepsilon +\varepsilon'+D_1^{t+1}-S_1^{t+1},-\varepsilon + \varepsilon' +S_2^{t+1}-D_2^{t+1})\\
    & + E\tilde{U}_{t+1}(z_1-\varepsilon +\varepsilon'+D_1^{t+1}-S_1^{t+1},-\varepsilon + \varepsilon' +S_2^{t+1}-D_2^{t+1})\\
    &- E\tilde{U}_{t+1}(z_1-\varepsilon+D_1^{t+1}-S_1^{t+1},-\varepsilon+S_2^{t+1}-D_2^{t+1})\\
    \ge &  -(r_{11}^t-r_{11}^{t+1})\varepsilon \\
    &+ E\tilde{U}_{t+1}(z_1+D_1^{t+1}-S_1^{t+1},S_2^{t+1}-D_2^{t+1}) - E\tilde{U}_{t+1}(z_1-\varepsilon +\varepsilon'+D_1^{t+1}-S_1^{t+1},-\varepsilon + \varepsilon' +S_2^{t+1}-D_2^{t+1})\\
    & + E\tilde{U}_{t+1}(z_1-\varepsilon +\varepsilon'+D_1^{t+1}-S_1^{t+1},  \varepsilon' +S_2^{t+1}-D_2^{t+1})- E\tilde{U}_{t+1}(z_1-\varepsilon+D_1^{t+1}-S_1^{t+1},S_2^{t+1}-D_2^{t+1}),
\end{align*}
where we obtain the inequality due to the increasing difference property of $\tilde{U}_{t+1}$ by induction.

It then follows that 
 \begin{align*}
     &[\tilde{U}(z_1,0) - \tilde{U}(z_1-\varepsilon,-\varepsilon)]
     - [\tilde{U}_t(z_1,\varepsilon) - \tilde{U}_t(z_1-\varepsilon,0)]\\
     \ge & r_{11}^{t+1} (\varepsilon-\varepsilon') + r_{22}^t \varepsilon - r_{12}^t (\varepsilon-\varepsilon') -r_{22}^{t+1} \varepsilon' \\
     &+ E\tilde{U}_{t+1}(z_1+D_1^t-S_1^t,S_2^t-D_2^t) - E\tilde{U}_{t+1}(z_1-\varepsilon +\varepsilon'+D_1^t-S_1^t,-\varepsilon + \varepsilon' +S_2^t-D_2^t)\\
     \ge& 0.
 \end{align*}
 By denoting $\tilde{\varepsilon}:=\varepsilon-\varepsilon'$, one can verify that the last inequality follows from Lemma \ref{lem:diffU-hor22}.
\Halmos
\endproof

We proceed to complete the proof of Proposition \ref{prop:protection-monotone}.
Let us focus on the case with $z_1\ge0$ and $z_2\ge 0$ to show that $p_{d,+}^{t,I\!B}$ is increasing in $I\!B$ and $p_{s,+}^{t,I\!B}$ is decreasing in $I\!B$.

Given the protection levels $p_{d,+}^{t,I\!B}$ and $p_{s,+}^{t,I\!B}$, the optimal matching quantity (between type 1 demand and type 2 supply) in round 2 is $q_t^*=(z_1-p_{d,+}^{t,I\!B})^+=(z_2-p_{s,+}^{t,I\!B})^+$.
The optimal matching quantity $q_t^*$ solves $\max_{0\le q\le \min\left\{ z_1,z_2 \right\}}\tilde{J}_t(q,\mathbf{z})$.
By Lemma \ref{lem:diffU-hor22}, $q_t^*$ is increasing in $z_1$ and $z_1$, with the increasing rate bounded from above by 1.

For a fixed value of $I\!B=z_1-z_2$,  we have 
 $z_2>p_{d,+}^{t,I\!B}$ for sufficiently large $z_2$ (in the mean time, sufficiently large $z_1$ and fixed $I\!B$). 
 Therefore $q_t^*=z_2-p_{s,+}^{t,I\!B}$.
If we increase $I\!B$ by further increasing $z_1$ while holding the value of $z_2$, $q_t^*=z_2-p_{s,+}^{t,I\!B}$ increases at a rate no greater than 1 (since $q_t^*$ increases in $z_1$ at a rate no greater than 1).
This is possible only when $p_{s,+}^{t,I\!B}$ is decreasing in $I\!B$ at a rate no greater than 1.

Similarly, with $I\!B$ fixed and both $z_1$ and $z_2$ are sufficiently large, we have $q_t^*=z_1-p_{d,+}^{t,I\!B}$.
We can then increase $I\!B$ be decreasing $z_2$ while holding $z_1$ fixed.
This will decrease $q_t^*=z_1-p_{d,+}^{t,I\!B}$ at a rate no greater than 1, which implies that $p_{d,+}^{t,I\!B}$ increases in $I\!B$ at a rate no greater than 1.
\Halmos
\endproof

\proof{Proof of Proposition \ref{prop:protection-lost-demand}.}
Given that $\alpha=0$ and $\beta=1$, we can rewrite equation
(\ref{eqn:horizontal22-alt2}) as follows.
\begin{align*}
    J_t(q,\mathbf{z}) = & r_{12}^t q^+ + r_{21}^t q^- + r_{11}^{t+1} E\min \left\{  D_1^{t+1}, (z_1-q)^-+S_1^{t+1}\right\}
     + r_{22}^{t+1}E\min\left\{ D_2^{t+1}, (z_2-q)^++S_2^{t+1} \right\}\nonumber \\
    & + EU_{t+1}(D_1^{t+1} -  (z_1-q)^--S_1^{t+1}, (z_2-q)^++S_2^{t+1}-D_2^{t+1}).
\end{align*}
Let us focus on the case with $z_1>0$ and $z_2>0$ (the case with negative $z_1$ and $z_2$ can be dealt with analogously).
In this case, the feasible set becomes $M(q)=\left\{ q\mid 0\le q\le \min\left\{ z_1,z_2 \right\} \right\}$.
Thus, for $\mathbf{z}\in M(q)\bigcap \mathbb{R}_{++}^2$,
\begin{align*}
    J_t(q,\mathbf{z}) = & r_{12}^t q +  r_{11}^{t+1} E\min \left\{  D_1^{t+1}, S_1^{t+1}\right\}
     + r_{22}^{t+1}E\min\left\{ D_2^{t+1}, z_2-q+S_2^{t+1} \right\}\nonumber \\
    & + EU_{t+1}(D_1^{t+1} -S_1^{t+1}, z_2-q+S_2^{t+1}-D_2^{t+1}),
\end{align*}
which is independent of $z_1$.
Let $\bar{J}_t(p):=-r_{12}^t p +  r_{22}^{t+1}E\min\left\{ D_2^{t+1}, p+S_2^{t+1} \right\}+EU_{t+1}(D_1^{t+1} -S_1^{t+1}, p+S_2^{t+1}-D_2^{t+1})$.
We have $J_t(q,\mathbf{z})=r_{12}^t z_2+ \bar{J}_t(z_2-q)$.
Given that $0\le q\le \min\left\{ z_1,z_2 \right\}$, the feasible range of $p$ is $p\in[(z_2-z_1)^+,z_2]$.
Let $p_{s,+}^t\in\arg\max_{p\ge 0} \bar{J}_t(p)$.
Then, the optimal solution to $\max_{p\in[(z_2-z_1)^+,z_2]} \bar{J}_t(p)$ is $p^{t*}=\max\left\{ (z_2-z_1)^+,z_2\wedge p_{s,+}^t \right\}=\max\left\{ z_2-z_1, z_2\wedge p_{s,+}^t\right\}$.
Thus, it is optimal to reduce type 2 supply to $\max\left\{ z_2-z_1, z_2\wedge p_{s,+}^t\right\}$ by matching demand type 1 and supply type 2 in round 2 of period $t$.

Similarly, we can prove that in the case with $z_1<0$ and $z_2<0$, it is optimal to reduce type 1 supply to $z_2-z_1,(-z_1)\wedge p_{s,-}^t$ by matching type 2 demand with type 1 supply in round 2 of period $t$. \Halmos
\endproof

\bigskip

\proof{Proof of Lemma \ref{lem:horizontal-dominance}.}
We prove part (i), and part (ii) follows symmetrically.

Suppose that along the direction $o\to d$, $i$ is closer to $j$ than $i'$ is.
We show that $(i,j)\succ_{\mathcal{M}_s} (i',j)$ in this case.
It is easy to see that $\dist_{i\leftarrow j}\le \dist_{i'\leftarrow j}$, which implies that $r_{ij}^t\ge r_{ij}^{t'}$ for all $t$.
It remains to verify that $r_{ij}^t-r_{i'j}^t\ge \alpha(r_{ij''}^{t+1}-r_{i'j''}^{t+1})$.
To see that, we have
\begin{align*}
    r_{ij}^t-r_{i'j}^t=(R_t-\dist_{i\leftarrow j}) - (R_t-\dist_{i'\leftarrow j})
    = \text{dist}_{i'\leftarrow j}-\text{dist}_{i\leftarrow j}=\text{dist}_{i'\leftarrow i}.
\end{align*}

Consider $j''\in\mathcal{S}$.
We consider the following two possibilities.

If $j''$ is located between $i$ and endpoint $d$, we have $r_{i\leftarrow j''}^{t+1}=0$ since $i$ is not accessible from $j''$.
Then, $\alpha(r_{ij''}^{t+1}-r_{i'j''}^{t+1})=-\alpha r_{i'j''}^{t+1}\le 0\le r_{ij}^t-r_{i'j}^t$.

If $j''$ is located between endpoint $o$ and $i$, 
then 
\begin{align*}
    &\alpha(r_{ij''}^{t+1}-r_{i'j}^{t+1}) 
    = \alpha[(R^{t+1}-\dist_{i\leftarrow j''}) - (R^{t+1}-\dist_{i'\leftarrow j''})]
    = \alpha (\dist_{i'\leftarrow j''} - \dist_{i\leftarrow j''})=\alpha\cdot \dist_{i'\leftarrow i}\\
    \le & \dist_{i'\leftarrow i}= r_{ij}^t-r_{i'j}^t.
\end{align*}

Therefore $(i,j)\succ_{\mathcal{M}_s} (i',j)$.

The above arguments also imply that $(i',j)\succ_{\mathcal{M}_s} (i,j)$ if along the direction $o\to d$, $i'$ is closer to $j$ than $i$ is.
This proves part (i)  \Halmos
\endproof

\bigskip
\proof{Proof of Proposition \ref{prop:horizontal-structure}.}
Part (i) follows directly from Lemma \ref{lem:horizontal-dominance} and Theorem \ref{cor:Neighboring-priority-strong}.

The condition $R^t$ decreasing with respect to $t$ in part (ii) ensures that $r_{ij}^t\ge r_{ij}^{t+1}$ for any period $t$.
Then, part (ii) follows from Proposition \ref{prop:greedy-optimal}. \Halmos
\endproof

\bigskip

\proof{Proof of Proposition \ref{prop:one-level}.}
We prove the following statement by induction, which would imply the proposition.

For any feasible matching policy $P=\left\{ \mathbf{Q}^t(\mathbf{x},\mathbf{y}) \right\}_{t=1,\ldots,T}$, we can construct another feasible matching policy $\bar{P}$ such that the property in the proposition is satisfied (i.e., type $k$ demand is always matched with type $k$ greedily).

First we consider the single-period problem (i.e., $t=1$).
Let the current matching decision in that single period be $\mathbf{Q}$ under the state $(\mathbf{x},\mathbf{y})$.
We construct a feasible decision $\bar{\mathbf{Q}}$ by repeatedly applying the following modification steps until the desired property is satisfied.

\bigskip

\noindent\emph{Modifying the matching decision for the single-period probem.}

\begin{enumerate}[{Modification Step} 1.]
     \item If $q_{k+1,k}>0$ and $q_{k,k-1}>0$ at the same time, let $\mathbf{Q}\leftarrow \mathbf{Q}+\varepsilon\mathbf{e}_{kk}^{m\times n}- \varepsilon\mathbf{e}_{k+1,k}^{m\times n}-\varepsilon\mathbf{e}_{k,k-1}$, where $\varepsilon=\min\left\{ q_{k+1,k},q_{k,k-1} \right\}$.

    \item If $q_{k+1,k}>0$ and $u_k>0$, let $\mathbf{Q}\leftarrow \mathbf{Q}+\varepsilon\mathbf{e}_{kk}^{m\times n}- \varepsilon\mathbf{e}_{k+1,k}^{m\times n}$, where $\varepsilon=\min\left\{ q_{k+1,k},u_k \right\}$.

    \item If $q_{k,k-1}>0$ and $v_k>0$, let $\mathbf{Q}\leftarrow \mathbf{Q}+\varepsilon\mathbf{e}_{kk}^{m\times n}- \varepsilon\mathbf{e}_{k,k-1}^{m\times n}$, where $\varepsilon=\min\left\{ q_{k,k-1},v_k \right\}$.
\end{enumerate}

The above steps transfer matching quantities from a low-priority pair to a high-priority (i.e., from $(k-1,k)$ to $(k,k)$ or from $(k,k+1)$ to $(k,k)$; priority is defined by $\succ$). 
Note that any quantity is transferred at most once (no quantity would be transferred out from $(k,k)$ since it is already a highest priority pair under $\succ$).
In the followings, we verify that for any quantity transferred, the matching reward associated with this quantity after the transfer is at least $\frac{1}{2}$ of its reward before the transfer.

In step 1 of the above procedure, before the transfer, a total reward $r_{k+1,k}^t \varepsilon+r_{k,k-1}^t \varepsilon$ is received by matching type $k-1$ demand with type $k$ supply for the quantity $\varepsilon$, and matching type $k$ demand with type $k+1$ supply for the same quantity.
After applying step 1, a total reward $r_{kk}^t \varepsilon$ is received. 
Thus, the ratio of the post-transfer reward to the pre-transfer reward for the affected matching quantity is $r_{kk}^t/(r_{k+1,k}^t+r_{k,k-1}^t)\ge \frac{1}{2}$, where the inequality holds because $r_{k+1,k}^t\le r_{kk}^t$ and $r_{k,k-1}^t\le r_{kk}^t$.

In step 2 of the above procedure, before the transfer, a reward $r_{k+1,k}^t\varepsilon$ is received from the matching quantity $\varepsilon$ between type $k+1$ demand and type $k$ supply.
After the transfer, a reward $r_{kk}^t$ is received by redirecting the quantity to $(k,k)$.
Thus, the ratio of the post-transfer reward to the pre-transfer reward for the affected matching quantity is $r_{kk}^t/r_{k+1,k}^t\ge 1$.

Step 3 is analogous to step 2. The ratio of the post-transfer reward to the pre-transfer reward for the affected matching quantity is $r_{kk}^t/r_{k,k-1}^t\ge 1$.

For any matching quantity that is not transferred, the ratio of the post-transfer reward to the pre-transfer reward for the affected matching quantity is exactly 1.

Let $\varepsilon_1,\ldots,\varepsilon_L$ be the matching quantities transferred at some point of the procedure, and $\varepsilon_{L+1}$ be the untransferred quantity.
Also we let $r_1^{\varepsilon_1},\ldots, r_{L+1}^{\varepsilon_{L+1}}$ be the matching rewards corresponding to $\varepsilon_1,\ldots,\varepsilon_{K+1}$ before applying the procedure, and $\bar{r}_1^{\varepsilon_1},\ldots, \bar{r}_{L+1}^{\varepsilon_{L+1}}$ the rewards after applying the procedure.
Then, the ratio of the total post-transfer reward to the total pre-transfer reward is
\begin{align*}
    \frac{\sum_{\ell=\ell}^{L+1} \bar{r}_\ell^{\varepsilon_\ell}}{ \sum_{\ell=\ell}^{L+1}r_\ell^{\varepsilon_\ell}}\ge \min_{k=1,\ldots,n} \frac{r_{kk}^t}{r_{k-1,k}^t+r_{k,k+1}^t} \ge \frac{1}{2}.
\end{align*}

Now suppose that the statement holds for any $T$-period problem.
We will show that it is also true for a $(T+1)$-period problem.

Let $\Omega_{T+1}$ be the set of all possible realizations of demand and supply realization over the $T+1$ periods.
For a realization $\omega\in \Omega_{T+1}$, let $\mathbf{Q}^t(\omega)$ be the matching decision in period $t$ under policy $P$, and $(\mathbf{x}^t(\omega|P),\mathbf{y}^t(\omega|P))$ the state in period $t$ under policy $P$ for the realization $\omega$.

As in the single-period problem, for any given realization $\omega$, we modify the matching decision in period 1 so that it satisfies the desired properties. 

\emph{Modification 1.}
If  $q^1_{k,k-1}(\omega)>0$ and $q^1_{k+1,k}(\omega)>0$ at the same time, let $\mathbf{Q}^1(\omega)\leftarrow \mathbf{Q}^1(\omega)+\varepsilon\mathbf{e}_{kk}^{m\times n}- \varepsilon\mathbf{e}_{k,k-1}^{m\times n}-\varepsilon\mathbf{e}_{k+1,k}^{m\times n}$, where $\varepsilon=\min\left\{ q_{k,k-1}^1(\omega),q_{k+1,k}^1(\omega) \right\}$. 
This leads to a reduction $r_{k,k-1}^1\varepsilon+r_{k+1,k}^1\varepsilon$ in rewards,
in exchange for an extra reward $r_{kk}^1 \varepsilon$ between type $k$ demand and type $k$ supply in period 1. 
The ratio of post-modification reward $r_{kk}^t \varepsilon$ to the pre-modification reward  $r_{k,k-1}^1\varepsilon+r_{k+1,k}^1\varepsilon$
is therefore $r_{kk}^1 /(r_{k,k-1}^1+r_{k+1,k}^1)\ge \frac{1}{2}$ since $r_{kk}^1\ge \max\left\{ r_{k,k-1}^1,r_{k+1,k}^1 \right\}$.

\emph{Modification Step 2.}
If $q_{k,k-1}^1(\omega)>0$ and $v_k^1(\omega)>0$, we let $\mathbf{Q}^1(\omega)\leftarrow \mathbf{Q}^1(\omega)+\varepsilon\mathbf{e}_{kk}^{m\times n}- \varepsilon\mathbf{e}_{k,k-1}^{m\times n}$, where $\varepsilon=\min\left\{ q_{k,k-1}^1(\omega),v_k^1(\omega) \right\}$.
This modification in period 1 reduces the the post-matching level $v_k^1(\omega)$ by $\varepsilon$,
and thus the capacity of type $k$ supply in the beginning of period 2, $y_k^2(\omega\mid P)$, by $\beta\varepsilon$.
We will need to modify the matching decision $\mathbf{Q}^2(\omega)$ in period 2 to keep it feasible.
Due to the reduction in $y_k^2(\omega\mid P)$, there exist $\varepsilon_{kk}^2\ge 0$ and $\varepsilon_{k+1,k}^2\ge 0$ such that $\varepsilon_{kk}^2+\varepsilon_{k+1,k}^2\le \beta \varepsilon$ and $\mathbf{Q}^2(\omega)-\varepsilon_{kk}^2\mathbf{e}_{kk}^{m\times n} - \varepsilon_{k+1,k}^2 \mathbf{e}_{k+1,k}^{m\times n}$ is feasible after we modify the matching decision in period 1.
The post-matching level of type $k$ supply is reduced by $\beta\varepsilon - (\varepsilon_{kk}^2+\varepsilon_{k+1,k}^2)$, compared with the post-matching level $v_k^2(\omega)$ before the modification.
Thus, the capacity of type $k$ supply in period 3 is reduced by $\beta[\beta\varepsilon - (\varepsilon_{kk}^2+\varepsilon_{k+1,k}^2)]\le \beta\varepsilon - (\varepsilon_{kk}^2+\varepsilon_{k+1,k}^2)$.
There exist $\varepsilon_{kk}^3\ge 0$ and $\varepsilon_{k+1,k}^3\ge 0$ such that $\varepsilon_{kk}^3+\varepsilon_{k+1,k}^3\le \beta^2\varepsilon - \beta\varepsilon_{kk}^2+\varepsilon_{k+1,k}^2$ and $\mathbf{Q}^3(\omega)-\varepsilon_{kk}^3\mathbf{e}_{kk}^{m\times n} - \varepsilon_{k+1,k}^3 \mathbf{e}_{k+1,k}^{m\times n}$ is feasible after we modify the decisions in periods 1 and 2.
Recursively, we can show that there exists $\varepsilon_{kk}^t\ge 0$ and $\varepsilon_{k+1,k}^t\ge 0$ for $t=2,\ldots,T$ such that $\sum_{t=2}^{T}\beta^{T-t} (\varepsilon_{kk}^t+\varepsilon_{k+1,k}^t)\le \beta^{T-1} \varepsilon$ and that $\mathbf{Q}^t(\omega)-\varepsilon_{kk}^t\mathbf{e}_{kk}^{m\times n} - \varepsilon_{k+1,k}^t \mathbf{e}_{k+1,k}^{m\times n}$ is feasible in period $t$ for the realization $\omega$, after we have modified the decisions in previous periods.
We will let $\mathbf{Q}^t(\omega)\leftarrow\mathbf{Q}^t(\omega)-\varepsilon_{kk}^t\mathbf{e}_{kk}^{m\times n} - \varepsilon_{k+1,k}^t \mathbf{e}_{k+1,k}^{m\times n}$
In this step of modification, there is a total reduction $r_{k,k-1}^t\varepsilon + \sum_{t=2}^{T}(r_{kk}^t\varepsilon_{kk}^t+ r_{k+1,k}^t\varepsilon_{k+1,k}^t)$ in rewards from period 1 to period $T$. 
In the mean time, in period 1, an extra reward $r_{kk}^1\varepsilon$ is received in period $t$.
The ratio of the latter to the former is 
$r_{kk}^1 \varepsilon/[r_{k,k-1}^1\varepsilon + \sum_{t=2}^{T}(r_{kk}^t\varepsilon_{kk}^t+ r_{k+1,k}^t\varepsilon_{k+1,k}^t)]\ge \frac{1}{2} $, where the inequality holds because  $r_{k,k-1}^1\le r_{kk}^1$ and
$\sum_{t=2}^{T}(r_{kk}^t\varepsilon_{kk}^t+ r_{k+1,k}^t\varepsilon_{k+1,k}^t)\le \sum_{t=2}^{T}r_{kk}^t(\varepsilon_{kk}^t+\varepsilon_{k+1,k}^t)\le \sum_{t=2}^{T}\beta^{-(t-1)}r_{kk}^1(\varepsilon_{kk}^t+\varepsilon_{k+1,k}^t)=r_{kk}^1\beta^{-(T-1)}\sum_{t=2}^{T}\beta^{(T-t)}(\varepsilon_{kk}^t+\varepsilon_{k+1,k}^t)\le r_{kk}^1\beta^{-(T-1)} \cdot r_{kk}^1\beta^{T-1}\varepsilon=r_{kk}^1 \varepsilon$.

\emph{Modification Step 3.}
If $q_{k+1,k}^1(\omega)>0$ and $u_k^1(\omega)>0$, we let $\mathbf{Q}^1(\omega)\leftarrow \mathbf{Q}^1(\omega)+\epsilon\mathbf{e}_{kk}^{m\times n}- \epsilon\mathbf{e}_{k+1,k}^{m\times n}$. 
Analogous to Modification Step 2, there exist $\epsilon_{kk}^t\ge 0$ and $\epsilon_{k,k-1}^t\ge 0$ such that $\sum_{t=2}^{T}\beta^{T-t} (\epsilon_{kk}^t+\epsilon_{k,k-1}^t)\le \beta^{T-1} \epsilon$.
We modifify the matching deicision as $\mathbf{Q}^t\leftarrow \mathbf{Q}^t(\omega)-\epsilon_{kk}^t\mathbf{e}_{kk}^{m\times n} - \epsilon_{k,k-1}^t \mathbf{e}_{k+1,k}^{m\times n}$.

After repeatedly applying the above modification steps, we arrive at a matching policy $\tilde{P}=\left\{ \tilde{\mathbf{Q}}^t \right\}_{t=1,\ldots,T}$ that satisfies the statement in period $t$.
For the affected matching quantities in those modification steps, the post-modification rewards are at least $\frac{1}{2}$ of the pre-modification rewards.

Now for any $t=2,\ldots,T$, $\tilde{\mathbf{Q}}^t(\omega)$ is a feasible decision in period $t$ for the realization $\omega$.
Let the policy $\tilde{P}\mid_{t=2,\ldots,T}$ be policy $\tilde{P}$ confined to periods $t=2,\ldots,T$.
By induction, we can construct another policy $\bar{P}\mid_{t=2,\ldots,T}$ for periods $t=2,\ldots,T$ such that the total reward received under $\bar{P}\mid_{t=2,\ldots,T}$ is at least $\frac{1}{2}$ of the total reward under $\tilde{P}\mid_{t=2,\ldots,T}$.

Finally, let us define policy $\bar{P}:=\left\{ \tilde{\mathbf{Q}}^1 \right\}\bigoplus \{\tilde{P}\mid_{t=2,\ldots,T}\}$, i.e., applying the matching decision  $\tilde{\mathbf{Q}}^1$ in period 1 and policy $\tilde{P}\mid_{t=2,\ldots,T}$ in the remaining periods. 
According to our analysis, the total reward under policy $\bar{P}$ is at least $\frac{1}{2}$ of that under the original policy $P$. \Halmos
\endproof

\bigskip

\proof{Proof of Lemma \ref{lem:vert-comp}.}
We prove that $(i,j)\succ_{\mathcal{M}_s} (i',j)$ if $i>i'$.
The second statement in the lemma would then follow symmetrically.

We have $r_{ij}^t=r_{id}^t+r_{js}^t\ge r_{i'd}^t+r_{js}^t=r_{i'j}^t$ given that $i'>i$.

For any $j''\in\mathcal{S}$, we have
$r_{ij}^t-r_{i'j}^t=r_{id}^t-r_{i'd}^t\ge r_{id}^t-r_{i'd}^{t+1}=r_{ij''}^{t+1}-r_{i'j''}^{t+1}$.

It remains to show that $r_{ij}^t+r_{i'j'}^t\ge r_{ij'}^t+r_{i'j}a^t$ for $i'>i$ and $j'>i$.
This holds trivially since 
It remains to show that $r_{ij}^t+r_{i'j'}^t= r_{ij'}^t+r_{i'j}a^t$ under the additive reward structure. \Halmos
\endproof

\bigskip

\proof{Proof of Proposition \ref{prop:vert-top-down}.}
The proposition follows directly from Lemma \ref{lem:vert-comp}.
More details are provided in Appendix \ref{app:vert-alt} on how the total matching quantity in a period uniquelly determines the matching quantities. \Halmos
\endproof

\bigskip

\proof{Proof of Proposition \ref{prop:1sa}.}
Let $P^{\text{OSA}[1,t],\text{Greedy}[t+1,T]}$ be the policy that applies the one-step-ahead policy up to period $t$, and uses greedy matching from period $t+1$ to period $T$. 
In the followings, we show that the policy $P^{\text{OSA}[1,t],\text{Greedy}[t+1,T]}$ achieves a higher total reward than $P^{\text{OSA}[1,t-1],\text{Greedy}[t,T]}$.

The two policies, $P^{\text{OSA}[1,t],\text{Greedy}[t+1,T]}$ and $P^{\text{OSA}[1,t-1],\text{Greedy}[t,T]}$ coincide with each other in periods $1,\ldots,t-1$, and therefore have the same expected rewards in those periods.

For any state in the beginning of period $t$, the policy $P^{\text{OSA}[1,t],\text{Greedy}[t+1,T]}$ uses the one-step-ahead policy in that period, which is optimal (for maximizing the total expected reward from period $t$ to period $T$) given that $P^{\text{OSA}[1,t],\text{Greedy}[t+1,T]}$ will use greedy matching from the next period on.
In contrast, the policy $P^{\text{OSA}[1,t-1],\text{Greedy}[t,T]}$ uses greedy matching in period $t$, which is suboptimal in response to the greedy matching it enforces from period $t+1$ to period $T$. Consequently, $P^{\text{OSA}[1,t],\text{Greedy}[t+1,T]}$ leads to a higher total expected reward from period $t$ to period $T$ than $P^{\text{OSA}[1,t-1],\text{Greedy}[t,T]}$.
The overall total expected matching reward from period $1$ to period $T$ is higher under $P^{\text{OSA}[1,t],\text{Greedy}[t+1,T]}$ than under $P^{\text{OSA}[1,t-1],\text{Greedy}[t,T]}$.

The one-step-ahead policy coincides with $P^{\text{OSA}[1,T-1],\text{Greedy}[T,T]}$, and the greedy matching policy coincides with $P^{\text{OSA}[1,0],\text{Greedy}[1,T]}$.
Thus, the former leads to a  higher total expected reward than the latter. \Halmos
\endproof

\bigskip

\proof{Proof of Proposition \ref{prop:vertical-mono}.}
We use the formulation (\ref{eqn:vertical-alt11})--(\ref{eqn:vertical-alt12}) in Appendix \ref{app:vert-alt} to prove the following lemma.

\begin{lemma}
  \label{lem:L-alpha=beta} Suppose $\alpha=\beta=1$. 
 $\tilde{V}_t(\tilde{\mathbf{x}},\tilde{\mathbf{y}})$ is $L^\natural$-concave in $(\tilde{\mathbf{x}},\tilde{\mathbf{y}})$ for $t=1,\ldots, T+1$,  and $\tilde{G}_t(Q,\tilde{\mathbf{x}},\tilde{\mathbf{y}})$ is $L^\natural$-concave in $(Q,\tilde{\mathbf{x}},\tilde{\mathbf{y}})$ for $t=1,\ldots,T$.
\end{lemma}
\proof{Proof of Lemma \ref{lem:L-alpha=beta}.}
The proof is by induction on $t$. 
Clearly, $\tilde{V}_{T+1}(\tilde{\mathbf{x}},\tilde{\mathbf{y}})=-\tilde{\mathbf{x}}\mathbf{U}_m^{-1}(\mathbf{r}_d^t)^{\tt{T}} - \tilde{\mathbf{y}}\mathbf{U}_n^{-1}(\mathbf{r}_s^t)^{\tt{T}}$ is $L^\natural$-concave in $(\tilde{\mathbf{x}},\tilde{\mathbf{y}})$. 
We suppose that $\tilde{V}_{t+1}(\tilde{\mathbf{x}},\tilde{\mathbf{y}})$ is $L^\natural$-concave in $(\tilde{\mathbf{x}},\tilde{\mathbf{y}})$.
Then by definition of $L^\natural$-concavity and submodularity,  for any given $\tilde{\mathbf{D}}^{t+1}$ and $\tilde{\mathbf{S}}^{t+1}$, $\tilde{V}_{t+1}(\tilde{\mathbf{x}}+\tilde{\mathbf{D}}^{t+1},\alpha\tilde{\mathbf{y}}+\tilde{\mathbf{S}}^{t+1})$ is $L^\natural$-concave in $(\tilde{\mathbf{x}},\tilde{\mathbf{y}})$.
 Now consider period $t$. 
Since $Q\le\min\{ \tilde{x}_n, \tilde{y}_m\}$.
\begin{align*}
    &   \tilde{V}_{t+1}((\tilde{\mathbf{x}}-Q\bm{1}^n)^++\tilde{\mathbf{D}}^{t+1}, (\tilde{\mathbf{y}}-Q\bm{1}^m)^+ +\tilde{\mathbf{S}}^{t+1})\\
  &=  \tilde{V}_{t+1}((\tilde{\mathbf{x}}_{[1,n-1]}-Q\bm{1}^{n-1})^++\tilde{\mathbf{D}}_{[1,n-1]}^{t+1},(\tilde{x}_n-Q)+\tilde{D}_n^{t+1}, (\tilde{\mathbf{y}}_{[1,m-1]}-Q\bm{1}^{m-1})^++\tilde{\mathbf{S}}_{[1,m-1]}^{t+1},(\tilde{y}_m-Q)+\tilde{S}_m^{t+1} ),
\end{align*} 
which is $L^\natural$-concave in $(Q, \tilde{\mathbf{x}},\tilde{\mathbf{y}})$
by applying \citet[Lemma 4]{Chen:2014uw} and noting the monotonicity proved in Lemma \ref{lem:vert-alt-mono}. 
By \citet[Proposition 2.3.4(c)]{SimchiLevi:2013tj}, $E_{\tilde{\mathbf{D}}^{t+1},\tilde{\mathbf{S}}^{t+1}}[\tilde{V}_{t+1}((\tilde{\mathbf{x}}-Q\bm{1}_m)^++\tilde{\mathbf{D}}^{t+1}, (\tilde{\mathbf{y}}-Q\bm{1}_n)^+ +\tilde{\mathbf{S}}^{t+1})]$ is $L^\natural$-concave in $(Q, \tilde{\mathbf{x}},\tilde{\mathbf{y}})$, thus the last term in (\ref{eqn:vertical-alt12}) is $L^\natural$-concave in $(Q, \tilde{\mathbf{x}},\tilde{\mathbf{y}})$.
The first two terms in (\ref{eqn:vertical-alt12}) are $L^\natural$-concave in $(Q, \tilde{\mathbf{x}},\tilde{\mathbf{y}})$, because $-(\tilde{x}_{i'}-Q)^+$ is supermodular in $(Q,\tilde{x}_{i'})$, $-(\tilde{y}_{j'}-Q)^+$ is supermodular in $(Q,\tilde{y}_{j'})$ and $L^\natural$-concavity is preserved under any nonnegative linear combination. Since the other terms are linear, 
$\tilde{G}_t(Q,\tilde{\mathbf{x}},\tilde{\mathbf{y}})$ is $L^\natural$-concave in $(Q,\tilde{\mathbf{x}},\tilde{\mathbf{y}})$. 
By \citet[Proposition 2.3.4(e)]{SimchiLevi:2013tj}, 
$\tilde{V}_t(\tilde{\mathbf{x}},\tilde{\mathbf{y}})$ is $L^\natural$-concave in $(\tilde{\mathbf{x}},\tilde{\mathbf{y}})$. 
This completes the induction.\Halmos
\endproof
\endproof

We now proceed to prove the proposition.

Since $L^\natural$-concavity implies supermodularity, by Lemma \ref{lem:L-alpha=beta}, 
$\tilde{G}_t(Q,\tilde{\mathbf{x}},\tilde{\mathbf{y}})$ is $L^\natural$-concave, a fortiori, supermodular in $(Q,\tilde{\mathbf{x}},\tilde{\mathbf{y}})$. 
By \citet[Theorem 2.2.8]{SimchiLevi:2013tj}, the optimal solution to (\ref{eqn:vertical-alt11}), denoted by $\hat{Q}_t(\tilde{\mathbf{x}},\tilde{\mathbf{y}})$, is nondecreasing in $(\tilde{\mathbf{x}},\tilde{\mathbf{y}})$. Since the higher the original state $(\mathbf{x},\mathbf{y})$, the higher the transformed state $(\tilde{\mathbf{x}},\tilde{\mathbf{y}})$, the optimal solution $Q^*_t(\mathbf{x},\mathbf{y})$, expressed in terms of the original state, is nondecreasing in  $(\mathbf{x},\mathbf{y})$.
Thus, $\frac{\partial Q^{t*}(\mathbf{x},\mathbf{y})}{\partial x_i}\ge 0$ for any $i\in\mathcal{D}$.

By the definition of $L^\natural$-concavity, $\tilde{G}_t(Q-\xi,\tilde{\mathbf{x}}-\xi\bm{1}_m,\tilde{\mathbf{y}}-\xi\bm{1}_n)$ is supermodular in $(Q,\tilde{\mathbf{x}},\tilde{\mathbf{y}},\xi)$.
Then, for $Q>\hat{Q}_t(\tilde{\mathbf{x}},\tilde{\mathbf{y}})+\epsilon$, we have
\begin{align*}
    \tilde{G}_t(Q,\tilde{\mathbf{x}}+\epsilon\bm{1}_m,\tilde{\mathbf{y}}+\epsilon\bm{1}_n)
    - \tilde{G}_t(\hat{Q}_t(\tilde{\mathbf{x}},\tilde{\mathbf{y}})+\epsilon,\tilde{\mathbf{x}}+\epsilon\bm{1}_m,\tilde{\mathbf{y}}+\epsilon\bm{1}_n)
    \le \tilde{G}_t(Q-\epsilon,\tilde{\mathbf{x}},\tilde{\mathbf{y}})
    - \tilde{G}_t(\hat{Q}_t(\tilde{\mathbf{x}},\tilde{\mathbf{y}}),\tilde{\mathbf{x}},\tilde{\mathbf{y}})\le 0,
\end{align*}
where the first inequality is  derived by definition of supermodularity and the second inequality is due to the optimality of $\hat{Q}_t$. 
This implies that any matching quantity $Q>\hat{Q}_t(\tilde{\mathbf{x}},\tilde{\mathbf{y}})+\epsilon$ is no better than $\hat{Q}_t(\tilde{\mathbf{x}},\tilde{\mathbf{y}})+\epsilon$ for the state $ (\tilde{\mathbf{x}}+\epsilon\bm{1}_m,\tilde{\mathbf{y}}+\epsilon\bm{1}_n)$.
Therefore, $\hat{Q}_t(\tilde{\mathbf{x}}+\epsilon\bm{1}_m,\tilde{\mathbf{y}}+\epsilon\bm{1}_n)\le \hat{Q}_t(\tilde{\mathbf{x}},\tilde{\mathbf{y}})+\epsilon$. 
By the monotonicity of $\hat{Q}_t(\tilde{\mathbf{x}},\tilde{\mathbf{y}})$, $\hat{Q}_t(\tilde{\mathbf{x}}+\epsilon\bm{1}_m,\tilde{\mathbf{y}})\le \hat{Q}_t(\tilde{\mathbf{x}}+\epsilon\bm{1}_m,\tilde{\mathbf{y}}+\epsilon\bm{1}_n)\le \hat{Q}_t(\tilde{\mathbf{x}},\tilde{\mathbf{y}})+\epsilon$. Expressed in the original state, $Q^*_t(\mathbf{x}+\epsilon \mathbf{e}_1^n, \mathbf{y})\le Q^*_t(\mathbf{x},\mathbf{y})+\epsilon$.
This implies that $\frac{\partial Q^{t*}(\mathbf{x},\mathbf{y})}{\partial x_1}\le 1$.

For any two original states $(\mathbf{x}+\epsilon\mathbf{e}_k^n,\mathbf{y})$ and $(\mathbf{x}+\epsilon\mathbf{e}_{k+1}^n,\mathbf{y})$, $k=1,\ldots,n-1$, their transformed states can be ordered as $(\tilde{\mathbf{x}}+\epsilon \mathbf{1}_{[k,m]},\tilde{\mathbf{y}})\ge (\tilde{\mathbf{x}}+\epsilon \mathbf{1}_{[k+1,m]},\tilde{\mathbf{y}})$, where $\mathbf{1}_{[k,m]}$ is an $n$-dimensional vector with the $k$-th up to $m$-th entry being one and the rest of the entries being all zeros. By the monotonicity of $\hat{Q}_t(\tilde{\mathbf{x}},\tilde{\mathbf{y}})$, $\hat{Q}_t(\tilde{\mathbf{x}}+\epsilon \mathbf{1}_{[k,m]},\tilde{\mathbf{y}})\ge\hat{Q}_t(\tilde{\mathbf{x}}+\epsilon \mathbf{1}_{[k+1,m]},\tilde{\mathbf{y}})$. 
This implies that $Q^{t*}(\mathbf{x}+\epsilon \mathbf{e}_k^m,\mathbf{y})\ge Q^{t*}(\mathbf{x}+\epsilon \mathbf{e}_{k+1}^m,\mathbf{y})$.
Thus we have $Q^{t*}(\mathbf{x}+\epsilon \mathbf{e}_k^m,\mathbf{y}) - Q^{t*}(\mathbf{x},\mathbf{y})\ge Q^{t*}(\mathbf{x}+\epsilon \mathbf{e}_{k+1}^m,\mathbf{y})- Q^{t*}(\mathbf{x},\mathbf{y})$,
which implies that $\frac{\partial Q^{t*}(\mathbf{x},\mathbf{y})}{\partial x_{k+1}}\le \frac{\partial Q^{t*}(\mathbf{x},\mathbf{y})}{\partial x_{k}}$. 

We have proved that $0\le \frac{\partial Q^{t*}(\mathbf{x},\mathbf{y})}{\partial x_i}\le 1$ and $\frac{\partial Q^{t*}(\mathbf{x},\mathbf{y})}{\partial x_{i}}\ge \frac{\partial Q^{t*}(\mathbf{x},\mathbf{y})}{\partial x_{i+1}}$.
Analogously we can show that 
 $0\le \frac{\partial Q^{t*}(\mathbf{x},\mathbf{y})}{\partial y_j}\le 1$ for all $j\in\mathcal{S}$, 
    and  $\frac{\partial Q^{t*}(\mathbf{x},\mathbf{y})}{\partial y_j} \ge \frac{\partial Q^{t*}(\mathbf{x},\mathbf{y})}{\partial y_{j+1}}$.
\Halmos
\endproof

\bigskip

\proof{Proof of Proposition \ref{prop:vert-1sa-structure}.}
As part of the one-step-ahead heuristic, the greedy matching policy is implemented from period $t+1$.
We first prove two lemmas on the greedy matching policy. 
With vertically differentiated types, the greedy matching policy also follows the top-down structure but does not reserve demand or supply. 
Let $V_t^g(\mathbf{x},\mathbf{y})$ be the expected total discounted surplus under the greedy matching policy from the current period $t$ to the end of the horizon, given the current state $(\mathbf{x},\mathbf{y})$.

We first prove the following lemma.

\begin{lemma}
    For $\epsilon>0$ and $\delta_k\ge 0$, $k=1,\ldots,j$ such that $\sum_{k=1}^{j}\delta_k=\epsilon$,
     the difference $V_t^g(\mathbf{x},\mathbf{y}-\sum_{k=1}^j\delta_k\mathbf{e}_k+\epsilon\mathbf{e}_j)-V_t^g(\mathbf{x},\mathbf{y})$ 
     depends only on $\delta_k$, $k=1,\ldots,j$, $\tilde{x}_m$ and $\mathbf{y}_{[1,j-1]}$.
     Symmetrically, $V_t^g(\mathbf{x}-\sum_{k=1}^i\delta_k\mathbf{e}_k+\epsilon\mathbf{e}_i,\mathbf{y})-V_t^g(\mathbf{x},\mathbf{y})$ depends only on $\delta_k$, $k=1,\ldots,i$, $\tilde{y}_n$ and $\mathbf{x}_{[1,i-1]}$, for $\epsilon>0$, $\delta_k\ge 0$ and $\sum_{k=1}^{i}\delta_k=\epsilon$.
    \label{lem:v-g-dep}
\end{lemma}

\proof{Proof of Lemma \ref{lem:v-g-dep}.}
We will focus on the difference $V_t^g(\mathbf{x},\mathbf{y}-\sum_{k=1}^j\delta_k\mathbf{e}_k+\epsilon\mathbf{e}_j)-V_t^g(\mathbf{x},\mathbf{y})$ and the other difference satisfies the desired property by symmetry.
If we define $\delta_{j+1}=\cdots =\delta_{m}=0$ and $\bm{\delta}=(\delta_1,\ldots,\delta_m)$, the difference can be rewritten as 
$V_t^g(\mathbf{x},\mathbf{y}-\bm{\delta}+\epsilon\mathbf{e}_j)-V_t^g(\mathbf{x},\mathbf{y})$.

We prove the lemma by induction. Suppose the desired property holds for $t+1$.

If $\tilde{x}_n<\tilde{y}_{j-1}$, 
there exists $1\le j'\le j-1$ such that $\tilde{y}_{j'-1}\le\tilde{x}_n<\tilde{y}_{j'}$. 
Under the greedy matching policy, all the demand and types $1,\ldots,j'-1$ supply is matched, a quantity $\tilde{x}_n-\tilde{y}_{j'-1}$ in type $j'$ supply is matched, and types $j'+1,\ldots,m$ supply will not be matched.
This leads to the post-matching levels $\mathbf{Y}=(\mathbf{0}_{[1,j'-1]},\tilde{y}_{j'}-\tilde{x}_n,\mathbf{y}_{[j'+1,m]})$ for the supply types.
Then $V_t^g(\mathbf{x},\mathbf{y})=\mathbf{r}_{d}^t\mathbf{x}^{\tt{T}}+\mathbf{r}_{[1,j'-1],s}^t\mathbf{y}_{[1,j'-1]}^{\tt{T}}+r_{j's}^t(\tilde{x}_n-\tilde{y}_{j'-1})+\gamma EV_{t+1}^g(\mathbf{D}^{t+1},\mathbf{Y}+\mathbf{S}^{t+1})$.

On the other hand, under the state $(\mathbf{x},\mathbf{y}-\bm{\delta}+\epsilon\mathbf{e}_j)$, all the demand will again be fully matched, and the total amounts of demand and supply do not change compared to the state $(\mathbf{x},\mathbf{y})$.
There exists $j'\le j''\le j-1$ such that types $1,\ldots,j''-1$ supply are fully matched and types $j''+1,\ldots,m$ supply are not matched.
This leads to the post-matching levels $\mathbf{Y}'=(\mathbf{0}_{[1,j''-1]},\tilde{y}_{j''}-\tilde{\delta}_{j''}-\tilde{x}_n,\mathbf{y}_{[j''+1,m]}-\bm{\delta}_{[j''+1,m]})+\epsilon \mathbf{e}_j$, where 
$\tilde{\delta}_k\defeq\sum_{\ell=1}^{k}\delta_\ell$ for $1\le k\le m$.
Then $V_t^g(\mathbf{x},\mathbf{y})=\mathbf{r}_d^t\mathbf{x}^{\tt{T}}+ \mathbf{r}^t_{[1,j''-1],s} (\mathbf{y}_{[1,j''-1]}-\bm{\delta}_{[1,j''-1]})^{\tt{T}} + r^t_{j''s} (\tilde{x}_n-\tilde{y}_{j''-1}+\tilde{\delta}_{j''-1}) + \gamma EV_{t+1}^g(\mathbf{D},\alpha\mathbf{Y}'+\mathbf{S})$.
Let
\begin{align*}
    \bm{\Delta}=\left\{
    \begin{aligned}
        & (\mathbf{0}_{[1,j']},\tilde{y}_{j'}-\tilde{x}_n, \mathbf{y}_{[j'+1,j''-1]},y_{j''}-(\tilde{y}_{j''}-\tilde{\delta}_{j''}-\tilde{x}_n),  \bm{\delta}_{[j''+1,m]} ) & \text{ if } j''>j',\\
        & (\mathbf{0}_{[1,j'-1]},\tilde{\delta}_{j'},\bm{\delta}_{[j'+1,m]})  & \text{ if }j''=j'.
\end{aligned}
\right.
\end{align*}

It is easy to verify that $\bm{\Delta}\ge 0$, $\sum_{k=1}^{m}\Delta_k=\sum_{k=1}^{m}\delta_k=\epsilon$, $\Delta_k=0$ for $k=j+1,\ldots,m$, and $\bm{\Delta}$ depends only on $\epsilon$, $\bm{\delta}$ and $(\tilde{x}_{n},\mathbf{y}_{[1,j-1]})$.
In addition, $\mathbf{Y}'=\mathbf{Y}-\bm{\Delta}+\epsilon\mathbf{e}_j$. 
Then $V_{t+1}^g(\mathbf{D}^{t+1},\mathbf{Y}+\mathbf{S}^{t+1}) - V_{t+1}^g(\mathbf{D}^{t+1},\mathbf{Y} -\bm{\Delta}+\epsilon\mathbf{e}_j+\mathbf{S}^{t+1})$ depends only on $\bm{\delta}$, $\tilde{x}_n$, $\tilde{D}_n^{t+1}$ and $\mathbf{y}_{[1,j-1]}$. (Note that $\mathbf{Y}_{[1,j-1]}$ is uniquely determined by $\mathbf{y}_{[1,j-1]}$.)
Then the difference
\begin{align*}
    &V_t^g(\mathbf{x},\mathbf{y}-\bm{\delta}+\epsilon\mathbf{e}_j)-V_t^g(\mathbf{x},\mathbf{y})\\
    =& \mathbf{r}_{[1,j'-1],s}^t\mathbf{y}_{[1,j'-1]}^{\tt{T}}+r_{j's}^t(\tilde{x}_n-\tilde{y}_{j'-1})
    - \mathbf{r}^t_{[1,j''-1],s} (\mathbf{y}_{[1,j''-1]}-\bm{\delta}_{[1,j''-1]})^{\tt{T}} - r_{j''s}^t (\tilde{x}_n-\tilde{y}_{j''-1}+\tilde{\delta}_{j''-1})\\
    &+  E \left[ V_{t+1}^g(\mathbf{D}^{t+1},\mathbf{Y}+\mathbf{S}^{t+1}) - V_{t+1}^g(\mathbf{D}^{t+1},\mathbf{Y} -\bm{\Delta}+\epsilon\mathbf{e}_j+\mathbf{S}) \right]
\end{align*}
depends only on $\bm{\delta}$, $\tilde{x}_n$ and $\mathbf{y}_{[1,j-1]}$.

If $\tilde{x}_n\ge \tilde{y}_{j-1}$, the greedy matching policy leads to the same post-matching levels under the two states $(\mathbf{x},\mathbf{y})$ and $(\mathbf{x},\mathbf{y}-\bm{\delta}+\epsilon\mathbf{e}_j)$. 
We see that $V_t^g(\mathbf{x},\mathbf{y}-\bm{\delta}+\epsilon\mathbf{e}_j)-V_t^g(\mathbf{x},\mathbf{y})=\sum_{k=1}^{j}\delta_k (r_j^s- r_k^s )$, which is independent of $(\mathbf{x},\mathbf{y})$.

Combining the above analysis, we see that the difference $V_t^g(\mathbf{x},\mathbf{y}-\bm{\delta}+\epsilon\mathbf{e}_j)-V_t^g(\mathbf{x},\mathbf{y})$ depends only on $\delta_k$, $k=1,\ldots,j$, $\tilde{x}_n$ and $\mathbf{y}_{[1,j-1]}$.\Halmos
\endproof

\begin{lemma}     \label{lem:V-g-dep2}
    The difference $V_t^g(\mathbf{x}+\epsilon\mathbf{e}_i^n,\mathbf{y}+\epsilon\mathbf{e}_j^m)-V_t^g(\mathbf{x},\mathbf{y})$
depends only on $\epsilon$ and $(\mathbf{x}_{[1,i-1]},\mathbf{y}_{[1,j-1]},\tilde{x}_n,\tilde{y}_m)$.
\end{lemma}
\proof{Proof of Lemma \ref{lem:V-g-dep2}.}
First, consider the case $\tilde{x}_n\le \tilde{y}_m$.
If $\tilde{x}_n\ge \tilde{y}_{j-1}$, then $V_t^g(\mathbf{x}+\epsilon\mathbf{e}_i^n,\mathbf{y}+\epsilon\mathbf{e}_j^m)-V_t^g(\mathbf{x},\mathbf{y})=(r_i^d+r_j^s)\epsilon$, which is independent of $(\mathbf{x},\mathbf{y})$.

If $\tilde{y}_{j'-1}\le \tilde{x}_n<\tilde{y}_{j'}$ for some $1\le j'\le j-1$, 
Under state $(\mathbf{x},\mathbf{y})$, 
the post-matching levels for the supply types are $\mathbf{v}=(\mathbf{0}_{[1,j'-1]},\tilde{y}_{j'}-\tilde{x}_n,\mathbf{y}_{[j'+1,m]})$.
The additional amount, $\epsilon$, of type $i$ demand and the extra quantities $\delta_{j'}\defeq\min\left\{ \epsilon,\tilde{y}_{j'}-\tilde{x}_n \right\}$, $\delta_{j'+1}\defeq\min\left\{ [\epsilon-(\tilde{y}_{j'}-\tilde{x}_n)]^+,y_{j'+1} \right\}$, $\ldots, \delta_{j-1}\defeq\min\left\{ [\epsilon-(\tilde{y}_{j-2}-\tilde{x}_n)]^+,y_{j-1}\right\}$,  $\delta_j\defeq[\epsilon-(\tilde{y}_{j-1}-\tilde{x}_n)]^+$ are matched for supply types $j',j'+1,\ldots,j-1,j$,
respectively, under the new state $(\mathbf{x}+\epsilon\mathbf{e}_i^n,\mathbf{y}+\epsilon\mathbf{e}_j^m)$.
Let $\bm{\delta}=(0,\ldots,0,\delta_{j'},\ldots,\delta_j,0,\ldots,0)\in\mathbb{R}_+^m$.
Note that $\bm{\delta}$ is a function of $(\tilde{x}_n,\mathbf{y}_{[1,j-1]})$.
Then we have
\begin{align*}
    &V_t^g(\mathbf{x}+\epsilon\mathbf{e}_i^n,\mathbf{y}+\epsilon\mathbf{e}_j^m)-V_t^g(\mathbf{x},\mathbf{y})\\
    =& r_{id}^t \epsilon + \sum_{k=1}^{m}r_k^s \delta_k^s 
    +  E\left[ V_{t+1}(\mathbf{D}^{t+1}, (\mathbf{v} -\bm{\delta} +\epsilon \mathbf{e}_j^m ) + \mathbf{S}^{t+1}) - V_{t+1}(\mathbf{D}^{t+1}, \mathbf{v}+\mathbf{S}^{t+1}) \right].
\end{align*}
Since $ V_{t+1}(\mathbf{D}^{t+1}, (\mathbf{v} -\bm{\delta} +\epsilon \mathbf{e}_j^m ) + \mathbf{S}^{t+1}) - V_{t+1}(\mathbf{D}^{t+1}, \mathbf{v}+\mathbf{S}^{t+1})$ depends only on $\epsilon$, $\bm{\delta}$, $\tilde{D}_n^{t+1}$ and $\mathbf{v}_{[1,j-1]}$ by Lemma \ref{lem:v-g-dep},
the difference $V_t^g(\mathbf{x}+\epsilon\mathbf{e}_i^n,\mathbf{y}+\epsilon\mathbf{e}_j^m)-V_t^g(\mathbf{x},\mathbf{y})
$ depends only on $\epsilon$, $\tilde{x}_n$ and $\mathbf{y}_{[1,j-1]}$, because $\bm{\delta}$ is defined in terms of $\epsilon$, $\tilde{x}_n$ and $\mathbf{y}_{[1,j-1]}$.

Now consider the case $\tilde{x}_n>\tilde{y}_m$. By symmetry, we can show that  $V_t^g(\mathbf{x}+\epsilon\mathbf{e}_i^n,\mathbf{y}+\epsilon\mathbf{e}_j^m)-V_t^g(\mathbf{x},\mathbf{y})
$ depends only on $\epsilon$, $\tilde{y}_m$ and $\mathbf{x}_{[1,i-1]}$.
Combining those two cases, the difference depends only on $\epsilon$ and $(\mathbf{x}_{[1,i-1]},\mathbf{y}_{[1,j-1]},\tilde{x}_n,\tilde{y}_m)$.\Halmos
\endproof

\endproof

\bigskip

We proceed to prove Proposition \ref{prop:vert-1sa-structure}.
We focus on the matching between type $i$ demand and type $j$ supply, which happens under topdown-matching only when $\tilde{x}_i>\tilde{y}_{j-1}$ and $\tilde{y}_j>\tilde{x}_{i-1}$.
Note that if $\tilde{x}_i\le \tilde{y}_{j-1}$, there would be no type $i$ demand left when we start to use type $j$ supply. 
If $\tilde{y}_j\le \tilde{x}_{i-1}$, there would be no type $j$ supply left when we start to use type $j$ demand.
This proves part (i)

Let $F_t(Q,\tilde{\mathbf{x}},\tilde{\mathbf{y}})$ be defined as in (\ref{eqn:osa-obj}) of Appendix \ref{sec_app:1sa}, which is the total expected reward to be received from period $t$ to period $T$ for using the total matching quantity $Q$ in period $t$ and greedy matching from period $t+1$ on.
Following the top-down structure, we start to match $i$ with $j$ when the matching quantity $Q$ reaches $\max\left\{ \tilde{x_{i-1}},\tilde{y}_{j-1} \right\}$, and complete the matching between $i$ and $j$ when the quantity reaches $\min\left\{ \tilde{x}_i,\tilde{y}_j \right\}$.
Within this range of $Q$, we have $(\tilde{\mathbf{x}}-Q\mathbf{1}_m)^+=(\mathbf{0}_{i-1},\tilde{\mathbf{x}}_{[i,m]}-Q\mathbf{1}_{m-i+1})$ and $(\tilde{\mathbf{y}}-Q\mathbf{1}_m)^+=(\mathbf{0}_{j-1},\tilde{\mathbf{y}}_{[j,n]}-Q\mathbf{1}_{n-j+1})$. 
For sufficiently small $\varepsilon>0$, we have
\begin{align*}
    &(\tilde{\mathbf{x}}-(Q+\varepsilon)\mathbf{1}_m)^+\mathbf{U}_m^{-1} (\mathbf{r}_{d}^t-\alpha \mathbf{r}_{d}^{t+1})^{\tt{T}} 
    - (\tilde{\mathbf{x}}-Q\mathbf{1}_m)^+\mathbf{U}_m^{-1} (\mathbf{r}_{d}^t-\alpha \mathbf{r}_{d}^{t+1})^{\tt{T}}\\
    =& \sum_{i''=i}^{m} (\tilde{x}_{i''}-Q-\varepsilon) [ (r_{i''d}^t-r_{i''+1,d}^{t})-\alpha (r_{i''d}^{t+1}-r_{i''+1,d}^{t+1})] - \sum_{i''=i}^{m} (\tilde{x}_{i''}-Q) [ (r_{i''d}^t-r_{i''+1,d}^{t})-\alpha (r_{i''d}^{t+1}-r_{i''+1,d}^{t+1})]\\
        =&  -\varepsilon\sum_{i''=i}^{m}[ (r_{i''d}^t-r_{i''+1,d}^{t})-\alpha (r_{i''d}^{t+1}-r_{i''+1,d}^{t+1})]\\
        =& -\varepsilon (r_{id}^t-\alpha r_{id}^{t+1}).
\end{align*}

Similarly, we have
\begin{align*}
    (\tilde{\mathbf{y}}-(Q+\varepsilon)\mathbf{1}_m)^+ \mathbf{V}_n^{-1}(\mathbf{r}_{s}^t- \beta \mathbf{r}_{s}^{t+1})^{\tt{T}}-(\tilde{\mathbf{y}}-Q\mathbf{1}_m)^+ \mathbf{V}_n^{-1}(\mathbf{r}_{s}^t- \beta \mathbf{r}_{s}^{t+1})^{\tt{T}}
    =-\varepsilon (r_{js}^t-\alpha r_{js}^{t+1}).
\end{align*}

By (\ref{eqn:osa-obj}),
\begin{align*}
         F_t(Q,\tilde{\mathbf{x}},\tilde{\mathbf{y}})=& \tilde{\mathbf{x}}\mathbf{U}_m^{-1} (\mathbf{r}_{d}^t)^{\tt{T}}
+ \tilde{\mathbf{y}}\mathbf{V}_n^{-1} (\mathbf{r}_{s}^t)^{\tt{T}} +E\tilde{\mathbf{D}}^{t+1}\mathbf{U}_m^{-1}+E\tilde{\mathbf{S}}^{t+1}\mathbf{V}_n^{-1}\nonumber\\
      & - (\tilde{\mathbf{x}}-Q\mathbf{1}_m)^+\mathbf{U}_m^{-1} (\mathbf{r}_{d}^t-\alpha \mathbf{r}_{d}^{t+1})^{\tt{T}}
       - (\tilde{\mathbf{y}}-Q\mathbf{1}_m)^+ \mathbf{V}_n^{-1}(\mathbf{r}_{s}^t- \beta \mathbf{r}_{s}^{t+1})^{\tt{T}}\nonumber \\
       &+ E\tilde{V}_{t+1}^g(\alpha(\tilde{\mathbf{x}}-Q\mathbf{1}_m)^+ +\tilde{\mathbf{D}}^{t+1},\beta (\tilde{\mathbf{y}}-Q\mathbf{1}_n)^++\tilde{\mathbf{S}}^{t+1}).
\end{align*}

Then, for $\max\left\{ \tilde{x}_{i-1},\tilde{y}_{j-1} \right\}\le Q\le \min\left\{ \tilde{x}_i,\tilde{y}_j \right\}$, 
\begin{align*}
   & F_t(Q+\varepsilon,\tilde{\mathbf{x}},\tilde{\mathbf{y}}) - F_t(Q,\tilde{\mathbf{x}},\tilde{\mathbf{y}}) \\
   =& \varepsilon (r_{id}^t-\alpha r_{id}^{t+1} + r_{js}^t-\alpha r_{js}^{t+1})\\
   &+ E\tilde{V}_{t+1}^g(\alpha(\tilde{\mathbf{x}}-(Q+\varepsilon)\mathbf{1}_m)^+ +\tilde{\mathbf{D}}^{t+1},\beta (\tilde{\mathbf{y}}-(Q+\varepsilon)\mathbf{1}_n)^++\tilde{\mathbf{S}}^{t+1})\\
  & -E\tilde{V}_{t+1}^g(\alpha(\tilde{\mathbf{x}}-Q\mathbf{1}_m)^+ +\tilde{\mathbf{D}}^{t+1},\beta (\tilde{\mathbf{y}}-Q\mathbf{1}_n)^++\tilde{\mathbf{S}}^{t+1})\\
  =& \varepsilon (r_{id}^t-\alpha r_{id}^{t+1} + r_{js}^t-\alpha r_{js}^{t+1})\\
  &+ E\tilde{V}_{t+1}^g(\alpha(\tilde{\mathbf{x}}-Q\mathbf{1}_m)^+ +\tilde{\mathbf{D}}^{t+1} - \alpha\varepsilon\sum_{i''=i}^{m}\mathbf{e}_{i''}^m,\beta (\tilde{\mathbf{y}}-Q\mathbf{1}_n)^++\tilde{\mathbf{S}}^{t+1}- \beta\varepsilon\sum_{j''=j}^{m}\mathbf{e}_{j''}^n)\\
  & -E\tilde{V}_{t+1}^g(\alpha(\tilde{\mathbf{x}}-Q\mathbf{1}_m)^+ +\tilde{\mathbf{D}}^{t+1},\beta (\tilde{\mathbf{y}}-Q\mathbf{1}_n)^++\tilde{\mathbf{S}}^{t+1})\\
=& \varepsilon (r_{id}^t + r_{js}^t)\\
&+ EV_{t+1}^g(\alpha(\tilde{\mathbf{x}}-Q\mathbf{1}_m)^+\mathbf{U}_m^{-1} +\mathbf{D}^{t+1} - \alpha\varepsilon\mathbf{e}_{i}^n,\beta (\tilde{\mathbf{y}}-Q\mathbf{1}_n)^+\mathbf{V}_n^{-1}+\mathbf{S}^{t+1}- \beta\varepsilon\mathbf{e}_{j}^m)\\
  & -EV_{t+1}^g(\alpha(\tilde{\mathbf{x}}-Q\mathbf{1}_m)^+ +\tilde{\mathbf{D}}^{t+1},\beta (\tilde{\mathbf{y}}-Q\mathbf{1}_n)^++\tilde{\mathbf{S}}^{t+1}).
\end{align*}

By Lemma \ref{lem:v-g-dep}, $F_t(Q+\varepsilon,\tilde{\mathbf{x}},\tilde{\mathbf{y}}) - F_t(Q,\tilde{\mathbf{x}}$ depends only on the first $i-1$ entries of $(\tilde{\mathbf{x}}-Q\mathbf{1}_m)^+$, the first $j-1$ entries of $(\tilde{\mathbf{y}}-Q\mathbf{1}_n)^+$, $\tilde{x}_m-Q$ and $\tilde{y}_n-Q$.
However, the first $i-1$ entries of $(\tilde{\mathbf{x}}-Q\mathbf{1}_m)^+$ and the first $j-1$ entries of $(\tilde{\mathbf{y}}-Q\mathbf{1}_n)^+$ are equal to zero for $\max\left\{ \tilde{x}_{i-1},\tilde{y}_{j-1} \right\}\le Q\le \min\left\{ \tilde{x}_i,\tilde{y}_j \right\}$.
Thus, the difference only depends on $\tilde{x}_m-Q$ and $\tilde{y}_n-Q$.
Then, the partial derivative $\frac{\partial}{\partial Q} F_t(Q,\tilde{\mathbf{x}},\tilde{\mathbf{y}})$ only depends on $\tilde{x}_m-Q$ and $\tilde{y}_n-Q$ for $\max\left\{ \tilde{x}_{i-1},\tilde{y}_{j-1} \right\}\le Q\le \min\left\{ \tilde{x}_i,\tilde{y}_j \right\}$.
(If $F_t$ is not differentiable in $Q$, we define $\frac{\partial}{\partial Q} F_t(Q+\varepsilon,\tilde{\mathbf{x}},\tilde{\mathbf{y}})$ as $\frac{\partial}{\partial Q} F_t(Q,\tilde{\mathbf{x}},\tilde{\mathbf{y}}):= \lim_{\varepsilon\to0+} [F_t(Q+\varepsilon,\tilde{\mathbf{x}},\tilde{\mathbf{y}})-F_t(Q,\tilde{\mathbf{x}},\tilde{\mathbf{y}})]/\varepsilon$.)
Consequently, there exists a function $\bar{F}_t(\tilde{x}_m-Q,\tilde{y}_n-Q)$ such that it is identical to the function $F_t(Q,\tilde{\mathbf{x}},\tilde{\mathbf{y}})$ up to a factor that is independent of $Q$.
To find the optimal solution to $\max_{\max\left\{ \tilde{x}_{i-1},\tilde{y}_{j-1} \right\}\le Q\le \min\left\{ \tilde{x}_i,\tilde{y}_j \right\}} {F}_t(Q,\tilde{\mathbf{x}}, \tilde{\mathbf{y}})$, it is equivalent to solve $\max_{\max\left\{ \tilde{x}_{i-1},\tilde{y}_{j-1} \right\}\le Q\le \min\left\{ \tilde{x}_i,\tilde{y}_j \right\}} \bar{F}_t(\tilde{x}_m-Q,\tilde{y}_n-Q)$.

Let $p_d:=\tilde{x}_m-Q$ and $p_s:=\tilde{y}_n-Q$ be the target levels for aggregate demand and supply, respectively. 
Then, $p_s=p_d-(\tilde{x}_m-\tilde{y}_n)=p_d-I\!B$.
To ensure that $\max\left\{ \tilde{x}_{i-1},\tilde{y}_{j-1} \right\}\le Q\le \min\left\{ \tilde{x}_i,\tilde{y}_j \right\}$, we require $\tilde{x}_m-\tilde{x}_i\wedge \tilde{y}_j\le p_d\le \tilde{x}_m-\tilde{x}_{i-1}\vee \tilde{y}_{j-1}$. 
The problem $\max_{\max\left\{ \tilde{x}_{i-1},\tilde{y}_{j-1} \right\}\le Q\le \min\left\{ \tilde{x}_i,\tilde{y}_j \right\}} \bar{F}_t(\tilde{x}_m-Q,\tilde{y}_n-Q)$ is further equivalent to the following problem.
\begin{align*}
    \max &\quad\bar{F}_t(p_d,p_d-I\!B)\\
    S.t. &\quad \tilde{x}_m-\tilde{x}_i\wedge \tilde{y}_j\le p_d\le \tilde{x}_m-\tilde{x}_{i-1}\vee \tilde{y}_{j-1}.
\end{align*}

Let $p_{ij,d}^{t,I\!B}:=\arg\max_{p_{ij,d}\ge I\!B^+} \bar{F}_t(p_d,p_d-I\!B)$ (recall that $\bar{F}_t$ is dependent on $i$ and $j$).
Then, it is optimal (under the one-step-ahead policy) to reduce the aggregate demand to the target level $p_{ij,d}^{t,I\!B}$, if $p_{ij,d}^{t,I\!B}$ lies within the range $[\tilde{x}_m-\tilde{x}_i\wedge \tilde{y}_j, \tilde{x}_m-\tilde{x}_{i-1}\vee \tilde{y}_{j-1}]$.
If $p_{ij,d}^{t,I\!B}> \tilde{x}_m-\tilde{x}_{i-1}\vee \tilde{y}_{j-1}$, then the aggregate demand is already below $p_{ij,d}^{t,I\!B}$ prior to matching $i$ with $j$. In that case, $i$ and $j$ are not matched at all.
If $p_{ij,d}^{t,I\!B}<\tilde{x}_m-\tilde{x}_i\wedge \tilde{y}_j$, then the aggregate demand would be still above  $p_{ij,d}^{t,I\!B}$ even if we fully match $i$ with $j$.
In that case, we would match $i$ with $j$ to the full extent.

By defining $p_{ij,s}^{t,I\!B}:=p_{ij,d}^{t,I\!B}-I\!B$, we can show the one-step-ahead policy aims to let the aggregate supply reach $p_{ij,s}^{t,I\!B}$ or as much as possible.

Part (ii) is thus proved. 

Similar to the proof of Lemma \ref{lem:L-alpha=beta}, we can show that $\tilde{V}^g_t$ is $L^\natural$-concave.
By utilizing the $L^\natural$-concavity, part (ii) follows similar analysis as in the proof of Proposition \ref{prop:vertical-mono}.
\Halmos
\endproof

\section{Connection between weak and strong compatibility}\label{sec_appendix:weak--strong}

We now discuss a connection between policies weakly respecting and those (strongly) respecting the partial order  $\succ_{\mathcal{M}}$. 
For a policy $P_W$ that weakly respects $\succ_{\mathcal{M}}$, let us construct another policy $P_S$ that (strongly) respects $\succ_{\mathcal{M}}$, such that $P_S$ leads to the same post-matching levels as $P_W$. 
To do this, let us consider $(i,j)\succ_{\mathcal{M}} (i',j)$.
In period $t$ under the policy $P_W$, the matching quantities and post-matching levels satisfy that $u_i^t=0$ if $q_{i'j}^t>0$. 
Suppose that $a_i^{t,\succ_{\mathcal{M}}}>0$.
Then by definition, $0<a_i^{t,\succ_{\mathcal{M}}}=x_i-\sum_{j'': (i,j'')\notin \mathcal{B}_{ij,L}^{\succ_{\mathcal{M}}}}q_{ij''}^t=x_i- \sum_{j'': (i,j'')\in \mathcal{A}}q_{ij''}^t +\sum_{j'': (i,j'')\in \mathcal{B}_{ij,L}^{\succ_{\mathcal{M}}}}q_{ij''}^t=u_i^t+\sum_{j'': (i,j'')\in \mathcal{B}_{ij,L}^{\succ_{\mathcal{M}}}}q_{ij''}^t= \sum_{j'': (i,j'')\in \mathcal{B}_{ij,L}^{\succ_{\mathcal{M}}}}q_{ij''}^t$.
Thus, there exists $(i,j'')\in \mathcal{B}_{ij,L}^{\succ_{\mathcal{M}}}$ (i.e., $(i,j)\succ_{\mathcal{M}} (i,j'')$) such that $q_{ij''}^t>0$.
We can reduce the matching quantities $q_{i'j}^t$ and $q_{ij'}^t$ by the same amount $\varepsilon:=\min\left\{ q_{i'j}^t,  q_{ij'}^t\right\}$, and increase the matching quantities $q_{ij}^t$ and $q_{i'j'}^t$ by $\varepsilon$ at the same time.
By doing so, we reduce both the matching quantity $q_{i'j}^t$ and $a_i^{t,\succ_{\mathcal{M}}}=u_i^t+\sum_{j'': (i,j'')\in \mathcal{B}_{ij,L}^{\succ_{\mathcal{M}}}}q_{ij''}^t$.
By repeatedly modifying the matching quantities in period $t$ as described above, we will eventually arrive at a policy $P_S$, which (strongly) respects $\mathcal{M}$ (i.e., either $q_{i'j}^t=0$ or $a_i^{t,\succ_{\mathcal{M}}}=0$).

Now consider an optimal policy $P_W^*$ that weakly respects $\succ_{\mathcal{M}}$.
Following the above logic, there exists a policy $P_S$ with the same post-matching levels $(\mathbf{u}^t,\mathbf{v}^t)$ as the policy $P_W^*$ in any period $t$.
Under both policies, in period $t$ with state $(\mathbf{x},\mathbf{y})$, the consumption of demand and supply is given by $(\mathbf{x}-\mathbf{u}^t,\mathbf{y}-\mathbf{v}^t)$.
Thus, if we have obtained the post-matching levels $(\mathbf{u}^t,\mathbf{v}^t)$ under policy $P_S$, by solving the single-period transportation problem with demand and supply consumption $(\mathbf{x}-\mathbf{u}^t,\mathbf{y}-\mathbf{v}^t)$ in each period $t$, we can arrive at the optimal policy $P_W^*$.
This implies that the optimal policy $P_W^*$ is only different from a certain policy $P_S$ (strongly) respecting $\mathcal{M}$ by solving a single-period transportation problem. But it would be challenging to identify such a policy $P_S$ by restricting the search within the class of policies that (strongly) respect $\mathcal{M}$. 

\section{Accounting for waiting costs}\label{sec_appendix:waiting}
{\label{marker:wait-cost} Even though we did not consider waiting costs of those demand and supply types that are not immediately matched, those costs can be easily incorporated into the matching rewards.
To see this, we consider the following two problems. In first problem, demand type $i$ (resp., supply type $j$) incurs a per-unit waiting cost $c_i^t$ (resp., $h_j^t$) in period $t$ if unmatched, and the unit matching reward between $i\in\mathcal{D}$ and $j\in\mathcal{S}$ in period $t$ is $r_{ij}^t$.
In the second problem, all waiting costs are equal to zero, and the unit matching reward between $i\in\mathcal{D}$ and $j\in\mathcal{S}$ in period $t$ is $\bar{r}_{ij}:=r_{ij}^t+ \sum_{\tau=t}^{T} \alpha^{\tau-t}c_j^\tau+\sum_{\tau=t}^{T} \beta^{\tau-t}h_j^\tau$.
Let us refer to the first problem as Problem (W) and the second one as Problem (NW).

\begin{proposition}
    Problem (W) is equivalent to Problem (NW), in the sense that they share the same optimal matching policy.\label{prop:W-NW}
\end{proposition}
\proof{Proof of Proposition \ref{prop:W-NW}.}
Let us focus on the problem (W) and its total reward less waiting costs.    

  For Problem (W), consider an arbitrary sample path of demand and supply realizations $\left\{ \mathbf{d}^t,\mathbf{s}^t \right\}_{t=1,\ldots,T}$, where $\mathbf{d}^t=(d_1^t,\ldots,d_m^t)$ and $\mathbf{s}^t=(s_1^t,\ldots,s_n^t)$ are the realizations of demand and supply in period $t$, respectively, and a set of matching decisions $\left\{ \mathbf{Q}^t \right\}_{t=1,\ldots,T}$ feasible under the aforementioned sample path. 
  For a quantity $q$ of type $i$ demand that arrives in period $t_1$ and matched in period $t_2\ge t_1$, it incurs waiting cost in periods $t_1,t_1+1,\ldots,t_2-1$.
  In each period $\tau$ ($t_1\le \tau\le t_2-1$), only a fraction $\alpha^{\tau-t_1}$ of the original amount is retained, due to the partial carry-over.
  Therefore, the corresponding original quantity in period $t_1$ (to the matched quantity in period $t_2$) is $q \alpha^{-(t_2-t_1)}$.
  In period $\tau$, the amount $q \alpha^{-(t_2-t_1)}\times \alpha^{\tau-t_1}=q \alpha^{-(t_2-\tau)}$ is retained and incurs the waiting cost $ q \alpha^{-(t_2-\tau) } c_i^\tau$. 
  Thus, the total waiting cost incurred by the original quantity $q\alpha^{-(t_2-t_1)}$ of type $i$ demand is $\sum_{\tau=t_1}^{t_2-1} \alpha^{-(t_2-\tau)} q c_i^\tau= \sum_{\tau=t_1}^{T} \alpha^{-(t_2-\tau)} q c_i^\tau - \sum_{\tau=t_2}^{T} \alpha^{-(t_2-\tau)} q c_i^\tau$.
  If the quantity $q$ waits until the end of period $T$ but is never matched, the total waiting cost is simply $\sum_{\tau=t_1}^{T} \alpha^{-(t_2-\tau)} c_i^\tau q $.
   The demand $d_i^t$ that arrives in period $t$, is either matched in some period $\tau$ ($t\le \tau\le T$), lost, or never matched until the end. 
   We let $\xi_i^{[t_1,t_2]}$ be the amount of $d_i^{t_1}$ that is matched in period $t_2$ ($t_1\le t_2\le T$), and $\xi_i^{[t_1,T+1]}$ be the amount of $d_i^{t_1}$ that waits until the end but is never matched.
   Then, $d_i^{t_1}=\sum_{t_2=t_1}^{T+1} \alpha^{-(t_2-t_1)}\xi_i^{[t_1,t_2]}$,
   and the total waiting cost incurred by $d_i^{t_1}$ is
   \begin{align*}
       &\sum_{t_2=t_1}^{T+1}[ \sum_{\tau=t_1}^{T} \alpha^{-(t_2-\tau)}  c_i^\tau \xi_i^{[t_1, t_2]} - \sum_{\tau=t_2}^{T} \alpha^{-(t_2-\tau)}  c_i^\tau \xi_i^{[t_1, t_2]}]\\
       =& \sum_{t_2=t_1}^{T+1} \xi_i^{[t_1,t_2]} [ \sum_{\tau=t_1}^{T} \alpha^{-(t_2-\tau)}  c_i^\tau  - \sum_{\tau=t_2}^{T} \alpha^{-(t_2-\tau)}  c_i^\tau ]\\
       =& \sum_{t_2=t_1}^{T+1} \xi_i^{[t_1,t_2]}  \sum_{\tau=t_1}^{T} \alpha^{-(t_2-\tau)}  c_i^\tau  - \sum_{t_2=t_1}^{T+1} \xi_i^{[t_1,t_2]}\sum_{\tau=t_2}^{T} \alpha^{-(t_2-\tau)}  c_i^\tau\\
       =&  \sum_{t_2=t_1}^{T+1} \xi_i^{[t_1,t_2]} \alpha^{-(t_2-t_1)} \sum_{\tau=t_1}^{T} \alpha^{-(\tau-t_1)}  c_i^\tau - \sum_{t_2=t_1}^{T+1} \xi_i^{[t_1,t_2]}\sum_{\tau=t_2}^{T} \alpha^{-(t_2-\tau)}  c_i^\tau\\
       =& d_i^{t_1} \sum_{\tau=t_1}^{T} \alpha^{-(\tau-t_1)}  c_i^\tau - \sum_{t_2=t_1}^{T} \xi_i^{[t_1,t_2]}\sum_{\tau=t_2}^{T} \alpha^{-(t_2-\tau)}  c_i^\tau.
   \end{align*}

   It follows that the total waiting cost incurred by type $i$ demand in all periods is,
\begin{align*}
    & \sum_{t_1=1}^{T}\left[ d_i^{t_1} \sum_{\tau=t_1}^{T} \alpha^{-(\tau-t_1)}  c_i^\tau - \sum_{t_2=t_1}^{T} \xi_i^{[t_1,t_2]}\sum_{\tau=t_2}^{T} \alpha^{-(t_2-\tau)}  c_i^\tau \right]\\
    =& \sum_{t_1=1}^{T} d_i^{t_1} \sum_{\tau=t_1}^{T} \alpha^{-(\tau-t_1)}  c_i^\tau
    - \sum_{t_1=1}^{T}\sum_{t_2=t_1}^{T} \xi_i^{[t_1,t_2]}\sum_{\tau=t_2}^{T} \alpha^{-(t_2-\tau)}  c_i^\tau\\
    =& \sum_{t_1=1}^{T} d_i^{t_1} \sum_{\tau=t_1}^{T} \alpha^{-(\tau-t_1)}  c_i^\tau
    - \sum_{t_2=1}^{T}\sum_{t_1=1}^{t_2} \xi_i^{[t_1,t_2]}\sum_{\tau=t_2}^{T} \alpha^{-(t_2-\tau)}  c_i^\tau.
   \end{align*}

   Let $\xi_i^{t_2}=\sum_{t_1=1}^{t_2} \xi_i^{[t_1,t_2]}$ be the total quantity of type $i$ demand matched in period $t_2$.
   Then, it is equal to the total quantity of type $i$ demand matched will all types of supply in period $t$, i.e., $\xi_i^{t_2}=\sum_{j=1}^{n} q_{ij}^{t_2}$.
   Thus, we can rewrite the total waiting cost incurred by type $i$ demand in all periods as,
   \begin{align*}
       &\sum_{t_1=1}^{T} d_i^{t_1} \sum_{\tau=t_1}^{T} \alpha^{-(\tau-t_1)}  c_i^\tau
    - \sum_{t_2=1}^{T}\sum_{t_1=1}^{t_2} \xi_i^{[t_1,t_2]}\sum_{\tau=t_2}^{T} \alpha^{-(t_2-\tau)}  c_i^\tau\\
    =& \sum_{t_1=1}^{T} d_i^{t_1} \sum_{\tau=t_1}^{T} \alpha^{-(\tau-t_1)}  c_i^\tau
    - \sum_{t_2=1}^{T} \xi_i^{t_2}\sum_{\tau=t_2}^{T} \alpha^{-(t_2-\tau)}  c_i^\tau\\
    =& \sum_{t_1=1}^{T} d_i^{t_1} \sum_{\tau=t_1}^{T} \alpha^{-(\tau-t_1)}  c_i^\tau
    - \sum_{t_2=1}^{T} \sum_{j=1}^{n} q_{ij}^{t_2}\sum_{\tau=t_2}^{T} \alpha^{-(t_2-\tau)}  c_i^\tau.
   \end{align*}

   Similarly, we can show that, the total waiting cost incurred by type $j$ supply is 
   \begin{align*}
       \sum_{t_1=1}^{T} s_j^{t_1} \sum_{\tau=t_1}^{T} \beta^{-(\tau-t_1)}  h_j^\tau
    - \sum_{t_2=1}^{T} \sum_{i=1}^{m} q_{ij}^{t_2}\sum_{\tau=t_2}^{T} \beta^{-(t_2-\tau)}  h_j^\tau.
   \end{align*}

   Then, the total matching reward less the waiting costs in all periods is
    \begin{align*}
        &\sum_{t=1}^{T} \sum_{i\in\mathcal{D},j\in\mathcal{S}}r_{ij}^tq_{ij}^t
        -\sum_{i\in\mathcal{D}}\left[  \sum_{t_1=1}^{T} d_i^{t_1} \sum_{\tau=t_1}^{T} \alpha^{-(\tau-t_1)}  c_i^\tau
         - \sum_{t_2=1}^{T} \sum_{j=1}^{n} q_{ij}^{t_2}\sum_{\tau=t_2}^{T} \alpha^{-(t_2-\tau)}  c_i^\tau\right]\\
         &- \sum_{j\in\mathcal{S}}\left[  \sum_{t_1=1}^{T} s_j^{t_1} \sum_{\tau=t_1}^{T} \beta^{-(\tau-t_1)}  h_j^\tau
    - \sum_{t_2=1}^{T} \sum_{i=1}^{m} q_{ij}^{t_2}\sum_{\tau=t_2}^{T} \beta^{-(t_2-\tau)}  h_j^\tau\right]\\
    =& -\sum_{t_1=1}^{T} \sum_{i\in\mathcal{D}}d_i^{t_1} \sum_{\tau=t_1}^{T} \alpha^{-(\tau-t_1)}  c_i^\tau - \sum_{t_1=1}^{T} \sum_{j\in\mathcal{S}}s_j^{t_1} \sum_{\tau=t_1}^{T} \beta^{-(\tau-t_1)}  h_j^\tau\\
    &+ \sum_{t=1}^{T} \sum_{i\in\mathcal{D},j\in\mathcal{S}}r_{ij}^tq_{ij}^t
    + \sum_{t=1}^{T} \sum_{i\in\mathcal{D}} \sum_{j=1}^{n} q_{ij}^{t}\sum_{\tau=t}^{T} \alpha^{-(t-\tau)}  c_i^\tau + \sum_{t=1}^{T}\sum_{j\in\mathcal{S}}  \sum_{i=1}^{m} q_{ij}^{t}\sum_{\tau=t}^{T} \beta^{-(t-\tau)}  h_j^\tau\\
    =& -\sum_{t=1}^{T} \sum_{i\in\mathcal{D}}d_i^{t} \sum_{\tau=t}^{T} \alpha^{-(\tau-t)}  c_i^\tau - \sum_{t=1}^{T} \sum_{j\in\mathcal{S}}s_j^{t} \sum_{\tau=t}^{T} \beta^{-(\tau-t)}  h_j^\tau\\
    &+ \sum_{t=1}^{T} \sum_{i\in\mathcal{D},j\in\mathcal{S}}
    \left[ r_{ij}^t + \sum_{\tau=t}^{T} \alpha^{\tau-t}c_j^\tau+\sum_{\tau=t}^{T} \beta^{\tau-t}h_j^\tau \right] q_{ij}^t\\
    =& -\sum_{t=1}^{T} \sum_{i\in\mathcal{D}}d_i^{t} \sum_{\tau=t}^{T} \alpha^{-(\tau-t)}  c_i^\tau - \sum_{t=1}^{T} \sum_{j\in\mathcal{S}}s_j^{t} \sum_{\tau=t}^{T} \beta^{-(\tau-t)}  h_j^\tau
    + \sum_{t=1}^{T} \sum_{i\in\mathcal{D},j\in\mathcal{S}} \bar{r}_{ij}^t q_{ij}^t.
    \end{align*}

    Note that the term $\sum_{t=1}^{T} \sum_{i\in\mathcal{D},j\in\mathcal{S}} \bar{r}_{ij}^t q_{ij}^t$ is the total reward for Problem (NW), if the same matching decisions are used. Consequently, under the same feasible matching policy, for any realization of demand and supply, the total reward for Problem (NW) is different from the total reward less waiting costs for Problem (W) 
    by $-\sum_{t=1}^{T} \sum_{i\in\mathcal{D}}d_i^{t} \sum_{\tau=t}^{T} \alpha^{-(\tau-t)}  c_i^\tau - \sum_{t=1}^{T} \sum_{j\in\mathcal{S}}s_j^{t} \sum_{\tau=t}^{T} \beta^{-(\tau-t)}  h_j^\tau$.
    The expected different between the two is $-\sum_{t=1}^{T} \sum_{i\in\mathcal{D}}ED_i^{t} \sum_{\tau=t}^{T} \alpha^{-(\tau-t)}  c_i^\tau - \sum_{t=1}^{T} \sum_{j\in\mathcal{S}}ES_j^{t} \sum_{\tau=t}^{T} \beta^{-(\tau-t)}  h_j^\tau$, which is a constant.
    Thus, the two problems are equivalent to each other.
    \Halmos
\endproof

\section{An alternative formulation of the $2\times 2$ horizontal model}
\label{app:horizontal22}

We reformulate the problem in terms of the new state $\mathbf{z}=(z_1,z_2)$ immediately prior to round 2 of a period.

Matching between an imperfect pair happens only if $z_1$ and $z_2$ have the same sign (i.e., $z_1z_2\ge 0$).
More specifically, we consider the following cases.

\underline{Case 1: $z_1\ge 0$ and $z_2\ge 0$.}
After round 1 matching, a quantity $z_1$ of type 1 demand is available to be matched with a quantity $z_2$ of type 2 supply.
 Let $q$ be the matching quantity in round 2 between type 1 demand and type 2 supply.
 We have $0\le q\le \min\left\{ z_1,z_2 \right\}$.
 After round 2 matching, the remaining quantity of type 1 demand is $z_1$ and that of type 2 supply is $z_2-q$.
 The post-matching state is therefore $(z_1-q,z_2-q)$.

 \underline{Case 2: $z_1<0$ and $z_2<0$.}
 After round 1 matching, a quantity $-z_1$ of type 1 supply is available to be matched with a quantity $-z_2$ of type 2 demand. 
  Let $-q$ be the matching quantity in the round 2 between type 2 demand and type 1 supply.
  We have $0\le -q\le \min\left\{ -z_1,-z_2 \right\}$, or equivalently, $\max\left\{ z_1,z_2 \right\}\le q\le 0$.
  After round 2 matching, the remaining quantity of type 1 supply is $-z_1+q$ and that of type 2 demand is $-z_2+q$.
  In other words, the post-matching state is $(z_1-q,z_2-q)$.

  \underline{Case 3: $z_1z_2<0$.}
  After round 1 matching, either there is only demand available or only supply available.  
  The matching quantity in round 2 is $q=0$.
    The post-matching state is $(z_1-q,z_2-q)=(z_1,z_2)$ (it is identical to the pre-matching state since there is no matching in round 2).

In any of the above three cases, the feasible space of matching decision in round 2 of a period $t$ is given by:
\begin{align}
    M(\mathbf{z})=\left\{ q\mid 0\le q\le \min\left( z_1,z_2 \right) \text{ or } \max\left( z_1,z_2 \right)\le q\le 0 \text{ or } q=0 \right\}.
\end{align}

One can easily verify that $M(\mathbf{z})$ is a lattice.

To reformulate the problem, we consider the total expected reward received from round 2 matching in period $t$ to the end of period $T$. 

In period $t$, the matching quantity between type 1 demand and type 2 supply is $q^+$, and that between type 2 demand and type 1 supply is $q^-$.
Thus, a total reward $r_{12}^t q^+ + r_{21}^t q^-$ is received in round 2 of period $t$.

Given that the post-matching state in period $t$ is $(z_1-q,z_2-q)$ after round 2, in the beginning of period $t+1$ the available quantity of type 1 demand is $\alpha(z_1-q)^++D_1^{t+1}$, that of type 2 demand is $\alpha(z_2-q)^-+D_2^{t+1}$, that of type 1 supply is $\beta(z_1-q)^-+S_1^{t+1}$, and that of type 2 supply is $\beta(z_2-q)^+ + S_2^{t+1}$.
In round 1 of period $t+1$, type 1 demand and type 1 supply will be matched  greedily, and so will type 2 demand and type 2 supply.
This results in the total expected reward $r_{11}^{t+1} E\min \left\{  \alpha (z_1-q)^++D_1^{t+1}, \beta (z_1-q)^-+S_1^{t+1}\right\}+r_{22}^{t+1}E\min\left\{ \alpha (z_2-q)^-+D_2^{t+1}, \beta (z_2-q)^++S_2^{t+1} \right\}$ in round 1 of period $t+1$. 
The state immediately prior to round 2 of period $t+1$ is $(\alpha(z_1-q)^++D_1^{t+1} - \beta (z_1-q)^--S_1^{t+1}, \beta(z_2-q)^++S_2^{t+1}-\alpha(z_2-q)^--D_2^{t+1})$.

Let us define $J_t(q,\mathbf{z})$ as the total expected reward received from round 2 of period $t$ until the end of period $T$ if the round 2 matching decision in period $t$ is $q$.
We also define $U_t(\mathbf{z})$ as the optimal total expected reward achievable (by using the optimal $q$) from round 2 of period $t$ until the end of period $T$.
We are now ready to present the reformulation.
\begin{align}
    U_t(\mathbf{z}) = & \max_{q\in M(\mathbf{z})} J_t(q,\mathbf{z}) \label{eqn:horizontal22-alt1}\\
    J_t(q,\mathbf{z}) = & r_{12}^t q^+ + r_{21}^t q^- + r_{11}^{t+1} E\min \left\{  \alpha(z_1-q)^++D_1^{t+1}, \beta(z_1-q)^-+S_1^{t+1}\right\}\nonumber \\
    & + r_{22}^{t+1}E\min\left\{ \alpha(z_2-q)^-+D_2^{t+1}, \beta(z_2-q)^++S_2^{t+1} \right\}\nonumber \\
    & + EU_{t+1}(\alpha(z_1-q)^++D_1^{t+1} - \beta (z_1-q)^--S_1^{t+1}, \beta(z_2-q)^++S_2^{t+1}-\alpha(z_2-q)^--D_2^{t+1}).\label{eqn:horizontal22-alt2}
\end{align}

We show the concavity of $U_t$ and $J_t$ in the following lemma.

\begin{lemma}\label{lem:concavity22}
$U_t(\mathbf{z})$ is concave in $\mathbf{z}$ and $J_t(q,\mathbf{z})$ is concave in $q$ for any given $\mathbf{z}$.
\end{lemma}
\proof{Proof of Lemma \ref{lem:concavity22}.}
Suppose that the (original) state in the beginning of period $t$ is given as $x_1=z_1$, $x_2=0$, $y_1=0$ and $y_2=z_2$.
The matching quantity in round 1 is zero since there is no type 2 demand or type 1 supply available.
By definition, we have $U_t(\mathbf{z})=V_t(z_1,0,0,z_2)$.
Since $V_t$ is concave, $U_t(\mathbf{z})$ is concave in $\mathbf{z}$.

To show that $J_t(q,\mathbf{z})$ is concave in $q$ for any given $\mathbf{z}$, we assume $z_1\ge 0$ and $z_2\ge 0$ without loss of generality.
The concavity of $J_t$ with respect to $q$ can be proved analogously.
\Halmos
\endproof

\section{An alternative formulation of the vertical model}
\label{app:vert-alt}

We reformulate the vertical model with a \emph{transformed} system state and the total matching quantity $Q$ as the decision variable in each period.

We define $\tilde{x}_i:=\sum_{k=1}^{i} x_k$ for $i=1,\ldots, m$ and $\tilde{y}_j:=\sum_{k=1}^{j} y_k$ ($\tilde{x}_0$ and $\tilde{y}_0$ are defined as zero) as the \emph{transformed} system state, $\tilde{u}_i=\sum_{k=1}^{i} u_k$ and $\tilde{v}_j=\sum_{k=1}^i v_k$ as the transformed post-matching levels.
In addition, let $\tilde{D}_i^t=\sum_{k=1}^{i} D_k^t$ and $\tilde{S}_j^t=\sum_{k=1}^{j} S_k^t$ be the transformed random variables that represents new arrivals of demand and supply in period $t$.

Let us define $\mathbf{U}_k$ as the $k\times k$ upper triangular matrix with all the entries on or above the diagonal equal to one.
Then the state transformation can be written in a matrix form: $\mathbf{x}\mathbf{U}_n=\tilde{\mathbf{x}}$ and $\mathbf{y}\mathbf{U}_m=\tilde{\mathbf{y}}$.
Equivalently, we can write $\tilde{\mathbf{x}}=\mathbf{x}\mathbf{U}_m^{-1}$ and $\tilde{\mathbf{y}}=\mathbf{y}\mathbf{V}_n^{-1}$.
Here $\mathbf{U}_m^{-1}$ and $\mathbf{V}_n^{-1}$ are the inverse matrices of $\mathbf{U}_m$ and $\mathbf{V}_n$, respectively.
One can easily verify that both $\mathbf{U}_m^{-1}$ and $\mathbf{V}_n^{-1}$ have all their diagonal entries equal to 1 and each off-diagonal entry right above a diagonal entry equal to $-1$.

The decision variable, total matching quantity $Q$, satisfies that $Q\le \min\left\{ \tilde{x}_m,\tilde{y}_n \right\}$.
Under top-down matching, a total quantity $\min\left\{ \tilde{x}_i,Q \right\}$ is consumed for types $1,\ldots,i$ demand combined, for any $i\in\mathcal{D}$.
Thus, the quantity of type $i$ demand being consumed is 
\begin{align*}
    \min\left\{ \tilde{x}_i,Q \right\}-\min\left\{ \tilde{x}_{i-1},Q \right\}
    =& \tilde{x}_i - \tilde{x}_{i-1} - (\tilde{x}_i-Q)^+ + (\tilde{x}_{i-1}-Q)^+,
\end{align*}
i.e., the total consumed quantity of types $1,\ldots,i$ demand less the total consumed quantity of types $1,\ldots,i-1$ demand.
This contributes the reward $r_{id}^t [\tilde{x}_i - \tilde{x}_{i-1} - (\tilde{x}_i-Q)^+ + (\tilde{x}_{i-1}-Q)^+]$ from demand type $i$.
Likewise, supply type $j$ contributes the reward $r_{js}^t [\tilde{y}_j-\tilde{y}_{j-1}-(\tilde{y}_j-Q)^+ + (\tilde{y}_{j-1}-Q)^+]$ in period $t$.
Consequently, the total reward received in period $t$ is
\begin{align*}
    & \sum_{i=1}^{m} r_{id}^t [\tilde{x}_i - \tilde{x}_{i-1} - (\tilde{x}_i-Q)^+ + (\tilde{x}_{i-1}-Q)^+]
    + \sum_{j=1}^{n} r_{js}^t [\tilde{y}_j-\tilde{y}_{j-1}-(\tilde{y}_j-Q)^+ + (\tilde{y}_{j-1}-Q)^+]\\
    =& \sum_{i=1}^{m} ( r_{id}^t -r_{i+1,d}^t) \tilde{x}_i
    + \sum_{i=1}^n ( r_{js}^t -r_{j+1,d}^t) \tilde{y}_j
    - \sum_{i=1}^{m} ( r_{id}^t -r_{i+1,d}^t) (\tilde{x}_i-Q)^+
    - \sum_{i=1}^n ( r_{js}^t -r_{j+1,d}^t) (\tilde{y}_j-Q)^+\\
    =& \tilde{\mathbf{x}}\mathbf{U}_m^{-1} (\mathbf{r}_{d}^t)^{\tt{T}}
       + \tilde{\mathbf{y}}\mathbf{V}_n^{-1} (\mathbf{r}_{s}^t)^{\tt{T}}
       - (\tilde{\mathbf{x}}-Q\mathbf{1}_m)^+\mathbf{U}_m^{-1} (\mathbf{r}_{d}^t)^{\tt{T}}
       - (\tilde{\mathbf{y}}-Q\mathbf{1}_m)^+ \mathbf{V}_n^{-1}(\mathbf{r}_{s}^t)^{\tt{T}}
\end{align*}
where $r_{m+1,d}^t=r_{n+1,s}^t:=0$.

In the end of period $t$, the remaining quantity of types $1,\ldots,i$ demand combined is $(\tilde{x}_i-Q)^+$ and the remaining quantity of types $1,\ldots,j$ supply combined is $(\tilde{y}_j-Q)^+$.
Thus, the transformed post-matching levels are given by $(\tilde{\mathbf{u}},\tilde{\mathbf{v}})=((\tilde{\mathbf{x}}-Q\mathbf{1}_m)^+,(\tilde{\mathbf{y}}-Q\mathbf{1}_n)^+)$.
The transformed state in period $t+1$ is $(\tilde{\mathbf{x}}_{t+1},\tilde{\mathbf{y}}_{t+1})=(\alpha(\tilde{\mathbf{x}}-Q\mathbf{1}_m)^++\tilde{\mathbf{D}}^{t+1},\beta (\tilde{\mathbf{y}}-Q\mathbf{1}_n)^++\tilde{\mathbf{S}}^{t+1})$,
which can be converted back to the original state as $(\mathbf{x}_{t+1},\mathbf{y}_{t+1})=(\tilde{\mathbf{x}}_{t+1}\mathbf{U}_m^{-1},\tilde{\mathbf{y}}_{t+1}\mathbf{V}_n^{-1})$.

If the total matching quantity in period is $Q$, the maximun total expected reward achievable from period $t$ to period $T$ is
\begin{align}
    G_t(Q,\tilde{\mathbf{x}},\tilde{\mathbf{y}}) 
    =& \tilde{\mathbf{x}}\mathbf{U}_m^{-1} (\mathbf{r}_{d}^t)^{\tt{T}}
       + \tilde{\mathbf{y}}\mathbf{V}_n^{-1} (\mathbf{r}_{s}^t)^{\tt{T}}
       - (\tilde{\mathbf{x}}-Q\mathbf{1}_m)^+\mathbf{U}_m^{-1} (\mathbf{r}_{d}^t)^{\tt{T}}
       - (\tilde{\mathbf{y}}-Q\mathbf{1}_m)^+ \mathbf{V}_n^{-1}(\mathbf{r}_{s}^t)^{\tt{T}}\nonumber \\
       &+ EV_{t+1}(\alpha(\tilde{\mathbf{x}}-Q\mathbf{1}_m)^+\mathbf{U}_m^{-1}+\tilde{\mathbf{D}}^{t+1}\mathbf{U}_m^{-1},\beta (\tilde{\mathbf{y}}-Q\mathbf{1}_n)^+\mathbf{V}_n^{-1}+\tilde{\mathbf{S}}^{t+1}\mathbf{V}_n^{-1})\label{eqn:vertical-alt02}
\end{align}
given the transformed state $(\tilde{\mathbf{x}},\tilde{\mathbf{y}})$ in period $t$.
The optimal total expected reward from period $t$ to period $T$ is thus
\begin{align}
    V_t(\tilde{\mathbf{x}}\mathbf{U}_m^{-1},\tilde{\mathbf{y}}\mathbf{V}_n^{-1})
    =\max_{0\le Q\le \tilde{x}_m\wedge \tilde{y}_n } G_t(Q,\tilde{\mathbf{x}},\tilde{\mathbf{y}}).\label{eqn:vertical-alt01}
\end{align}

Let $\tilde{V}_t(\tilde{\mathbf{x}},\tilde{\mathbf{y}}):=V_t(\tilde{\mathbf{x}}\mathbf{U}_m^{-1},\tilde{\mathbf{y}}\mathbf{U}_n^{-1})-\tilde{\mathbf{x}}\mathbf{U}_m^{-1}(\mathbf{r}_d^t)^{\tt{T}} - \tilde{\mathbf{y}}\mathbf{U}_n^{-1}(\mathbf{r}_s^t)^{\tt{T}}$.
Equations (\ref{eqn:vertical-alt02}) and (\ref{eqn:vertical-alt01}) are equivalent to:
\begin{align}
    & \tilde{V}_t(\tilde{\mathbf{x}},\tilde{\mathbf{y}})
    =\max_{0\le Q\le \tilde{x}_m\wedge \tilde{y}_n} \tilde{G}_t(\tilde{\mathbf{x}},\tilde{\mathbf{y}}).\label{eqn:vertical-alt11}\\
    & \tilde{G}_t(Q,\tilde{\mathbf{x}},\tilde{\mathbf{y}}) 
    =
    - (\tilde{\mathbf{x}}-Q\mathbf{1}_m)^+\mathbf{U}_m^{-1} (\mathbf{r}_{d}^t-\alpha\mathbf{r}_d^{t+1})^{\tt{T}}
    - (\tilde{\mathbf{y}}-Q\mathbf{1}_m)^+ \mathbf{V}_n^{-1}(\mathbf{r}_{s}^t - \beta \mathbf{r}_{s}^{t+1})^{\tt{T}}\nonumber \\
       &\quad\quad\quad\quad+ E\tilde{V}_{t+1}(\alpha(\tilde{\mathbf{x}}-Q\mathbf{1}_m)^++\tilde{\mathbf{D}}^{t+1},\beta (\tilde{\mathbf{y}}-Q\mathbf{1}_n)^++\tilde{\mathbf{S}}^{t+1})\label{eqn:vertical-alt12}
\end{align}

Since $V_{T+1}(\mathbf{x},\mathbf{y})\equiv0$, we have $\tilde{V}_{T+1}(\tilde{\mathbf{x}},\tilde{\mathbf{y}}) \equiv -\tilde{\mathbf{x}}\mathbf{U}_m^{-1}(\mathbf{r}_d^t)^{\tt{T}} - \tilde{\mathbf{y}}\mathbf{U}_n^{-1}(\mathbf{r}_s^t)^{\tt{T}}$.

Finally, we have the following property for the function $\tilde{V}_t(\tilde{\mathbf{x}},\tilde{\mathbf{y}})$.

\begin{lemma}
    Suppose that $r_{id}^t-r_{i+1,d}^t\ge \alpha(r_{id}^{t+1}-r_{i+1,d}^{t+1}) $  and $r_{js}^t-r_{j+1,s}^t\ge \beta(r_{js}^{t+1}-r_{j+1,s}^{t+1}) $ for $i=1,\ldots,m-1$, $j=1,\ldots,n-1$ and $t=1,\ldots,T-1$.
    Then, for any period $t=1,\ldots,T$, the function $\tilde{V}_t(\tilde{\mathbf{x}},\tilde{\mathbf{y}})$ is decreasing in $x_i$ for $i=1,\ldots,m-1$ and in $y_j$for all $j=1,\ldots,n-1$.
    \label{lem:vert-alt-mono}
\end{lemma}
\proof{Proof of Lemma \ref{lem:vert-alt-mono}.}
By definition of the function $\tilde{V}_t$, for $i=1,\ldots,m-1$, we have
\begin{align}
	&\tilde{V}_t(\tilde{\mathbf{x}}+\varepsilon\mathbf{e}_i^m,\tilde{\mathbf{y}}) -\tilde{V}_t(\tilde{\mathbf{x}},\tilde{\mathbf{y}})\nonumber \\
	=& V_t(\mathbf{x}+\varepsilon\mathbf{e}_i^m-\varepsilon\mathbf{e}_{i+1}^m,\mathbf{y}) -(\mathbf{x}+\varepsilon\mathbf{e}_i^m-\varepsilon\mathbf{e}_{i+1}^m)(\mathbf{r}_d^t)^{\tt{T}} -\mathbf{y} (\mathbf{r}_s^t)^{\tt{T}} - V_t(\mathbf{x},\mathbf{y}) +\mathbf{x}(\mathbf{r}_d^t)^{\tt{T}} -\mathbf{y} (\mathbf{r}_s^t)^{\tt{T}}\nonumber\\
    =&V_t(\mathbf{x}+\varepsilon\mathbf{e}_i^m-\varepsilon\mathbf{e}_{i+1}^m,\mathbf{y})- V_t(\mathbf{x},\mathbf{y}) - (r_{id}^t-r_{i+1,d}^t)\varepsilon, \label{eqn:Vtilde-bound}
\end{align}
where $\mathbf{x}=\tilde{\mathbf{x}}\mathbf{U}_m^{-1}$ and $\mathbf{y}=\tilde{\mathbf{y}}\mathbf{V}_n^{-1}$.

By Lemma \ref{lem:ineq-0}, there exists $\lambda_{j'}^\tau\ge 0$ for $j'=1,\ldots,n$ and $\tau=t,\ldots,T$ such that $\sum_{\tau=t}^{T} \alpha^{-(\tau-t)}\sum_{j'=1}^{n} \lambda_{j'}^\tau \le \varepsilon$ and $V_t(\mathbf{x}+\varepsilon\mathbf{e}_i^m-\varepsilon\mathbf{e}_{i+1}^m,\mathbf{y})- V_t(\mathbf{x},\mathbf{y})\le \sum_{\tau=t}^{T} \sum_{j'=1}^{n} \lambda_{j'}^\tau (r_{ij'}^\tau - r_{i+1,j'}^\tau)$.
Following (\ref{eqn:Vtilde-bound}), we have
\begin{align*}
    & \tilde{V}_t(\tilde{\mathbf{x}}+\varepsilon\mathbf{e}_i^m,\tilde{\mathbf{y}}) -\tilde{V}_t(\tilde{\mathbf{x}},\tilde{\mathbf{y}})\\
    \le& \sum_{\tau=t}^{T} \sum_{j'=1}^{n} \lambda_{j'}^\tau (r_{ij'}^\tau - r_{i+1,j'}^\tau)- (r_{id}^t-r_{i+1,d}^t)\varepsilon\\
    \le & \sum_{\tau=t}^{T} \alpha^{-(\tau-t)}\sum_{j'=1}^{n} \lambda_{j'}^\tau (r_{ij'}^t - r_{i+1,j'}^t)- (r_{id}^t-r_{i+1,d}^t)\varepsilon\\
    =& \sum_{\tau=t}^{T} \alpha^{-(\tau-t)}\sum_{j'=1}^{n} \lambda_{j'}^\tau (r_{id}^t - r_{i+1,d}^t)- (r_{id}^t-r_{i+1,d}^t)\varepsilon\\
    =& [ \sum_{\tau=t}^{T} \alpha^{-(\tau-t)}\sum_{j'=1}^{n} \lambda_{j'}^\tau -\varepsilon]  (r_{id}^t-r_{i+1,d}^t)\\
    \le & 0.
\end{align*}
Therefore, $\tilde{V}_t(\tilde{\mathbf{x}},\tilde{\mathbf{y}})$ is decreasing in $\tilde{x}_i$ for $i=1,\ldots,m-1$.
Similarly, we can show that $\tilde{V}_t(\tilde{\mathbf{x}},\tilde{\mathbf{y}})$ is decreasing in $\tilde{y}_j$ for $j=1,\ldots,n-1$.
\Halmos
\endproof

\section{One-step-ahead policy for the vertical model}

\label{sec_app:1sa}

The one-step-ahead policy assumes greedy matching from the next period to the end of the horizon.
Let $V_t^g(\mathbf{x},\mathbf{y})$ be the total expected reward received under the greedy policy from period $t$ to period $T$, given that the state in period $t$ is $(\mathbf{x},\mathbf{y})$.
The one-step-ahead policy chooses the decision in period $t$ to maximize the sum of the immediate reward in period $t$ and the future expected reward $V_t^g(\mathbf{x},\mathbf{y})$. 
In this appendix, we explore properties of the function $V_t^g(\mathbf{x},\mathbf{y})$ and the one-step-ahead policy.

Let $\tilde{V}_t^g(\mathbf{x},\mathbf{y}):= V_t^g(\mathbf{x},\mathbf{y})-\tilde{\mathbf{x}}\mathbf{U}_m^{-1}(\mathbf{r}_d^t)^{\tt{T}} - \tilde{\mathbf{y}}\mathbf{U}_n^{-1}(\mathbf{r}_s^t)^{\tt{T}}$.
In the following lemma, we present recursive equation satisfied by $\tilde{V}_t^g(\mathbf{x},\mathbf{y})$.

\begin{lemma}
    The function $\tilde{V}_t^g(\mathbf{x},\mathbf{y})$ satisfies the following recursive equations:
\begin{align}
    \tilde{V}_t^g(\tilde{\mathbf{x}},\tilde{\mathbf{y}})
    =& E\tilde{\mathbf{D}}^{t+1}\mathbf{U}_m^{-1}+E\tilde{\mathbf{S}}^{t+1}\mathbf{V}_n^{-1} - (\tilde{\mathbf{x}}-\tilde{y}_n\mathbf{1}_m)^+\mathbf{U}_m^{-1} (\mathbf{r}_{d}^t-\alpha \mathbf{r}_{d}^{t+1})^{\tt{T}}
    - (\tilde{\mathbf{y}}-\tilde{x}_m\mathbf{1}_n)^+ \mathbf{V}_n^{-1}(\mathbf{r}_{s}^t- \beta \mathbf{r}_{s}^{t+1})^{\tt{T}}\nonumber \\
    &+ E\tilde{V}_{t+1}^g(\alpha(\tilde{\mathbf{x}}-\tilde{y}_n\mathbf{1}_m)^+ +\tilde{\mathbf{D}}^{t+1},\beta (\tilde{\mathbf{y}}-\tilde{x}_m\mathbf{1}_n)^++\tilde{\mathbf{S}}^{t+1})
    \label{eqn:V-tilde-g}
\end{align}
    \label{lem:Vg-recursion}
\end{lemma}
\proof{Proof of Lemma \ref{lem:Vg-recursion}.}
Assuming that the greedy policy will be used from period $t+1$ to period $T$,
the total expected reward received by using a top-down matching in period $t$  with the total matching quantity $Q$ is 
\begin{align}
    F_t(Q,\tilde{\mathbf{x}},\tilde{\mathbf{y}})
    =&  \tilde{\mathbf{x}}\mathbf{U}_m^{-1} (\mathbf{r}_{d}^t)^{\tt{T}}
       + \tilde{\mathbf{y}}\mathbf{V}_n^{-1} (\mathbf{r}_{s}^t)^{\tt{T}}
       - (\tilde{\mathbf{x}}-Q\mathbf{1}_m)^+\mathbf{U}_m^{-1} (\mathbf{r}_{d}^t)^{\tt{T}}
       - (\tilde{\mathbf{y}}-Q\mathbf{1}_m)^+ \mathbf{V}_n^{-1}(\mathbf{r}_{s}^t)^{\tt{T}}\nonumber \\
     &+ EV_{t+1}^g(\alpha(\tilde{\mathbf{x}}-Q\mathbf{1}_m)^+\mathbf{U}_m^{-1}+\tilde{\mathbf{D}}^{t+1}\mathbf{U}_m^{-1},\beta (\tilde{\mathbf{y}}-Q\mathbf{1}_n)^+\mathbf{V}_n^{-1}+\tilde{\mathbf{S}}^{t+1}\mathbf{V}_n^{-1}\nonumber
\\
=& \tilde{\mathbf{x}}\mathbf{U}_m^{-1} (\mathbf{r}_{d}^t)^{\tt{T}}
+ \tilde{\mathbf{y}}\mathbf{V}_n^{-1} (\mathbf{r}_{s}^t)^{\tt{T}} +E\tilde{\mathbf{D}}^{t+1}\mathbf{U}_m^{-1}+E\tilde{\mathbf{S}}^{t+1}\mathbf{V}_n^{-1}\nonumber\\
      & - (\tilde{\mathbf{x}}-Q\mathbf{1}_m)^+\mathbf{U}_m^{-1} (\mathbf{r}_{d}^t-\alpha \mathbf{r}_{d}^{t+1})^{\tt{T}}
       - (\tilde{\mathbf{y}}-Q\mathbf{1}_m)^+ \mathbf{V}_n^{-1}(\mathbf{r}_{s}^t- \beta \mathbf{r}_{s}^{t+1})^{\tt{T}}\nonumber \\
       &+ E\tilde{V}_{t+1}^g(\alpha(\tilde{\mathbf{x}}-Q\mathbf{1}_m)^+ +\tilde{\mathbf{D}}^{t+1},\beta (\tilde{\mathbf{y}}-Q\mathbf{1}_n)^++\tilde{\mathbf{S}}^{t+1})  \label{eqn:osa-02}
\end{align}

If greedy matching is used in period $t$, then the matching quantity $Q$ is equal to $\tilde{x}_m\wedge \tilde{y}_n$.
Thus,
\begin{align}
    V_t^g(\tilde{\mathbf{x}},\tilde{\mathbf{y}}) 
    =& F_t(\tilde{x}_m\wedge \tilde{y}_n,\tilde{\mathbf{x}},\tilde{\mathbf{y}})\\
    =& \tilde{\mathbf{x}}\mathbf{U}_m^{-1} (\mathbf{r}_{d}^t)^{\tt{T}}
+ \tilde{\mathbf{y}}\mathbf{V}_n^{-1} (\mathbf{r}_{s}^t)^{\tt{T}} +E\tilde{\mathbf{D}}^{t+1}\mathbf{U}_m^{-1}+E\tilde{\mathbf{S}}^{t+1}\mathbf{V}_n^{-1}\nonumber\\
      & - (\tilde{\mathbf{x}}-\tilde{x}_m\wedge \tilde{y}_n\mathbf{1}_m)^+\mathbf{U}_m^{-1} (\mathbf{r}_{d}^t-\alpha \mathbf{r}_{d}^{t+1})^{\tt{T}}
       - (\tilde{\mathbf{y}}-\tilde{x}_m\wedge \tilde{y}_n\mathbf{1}_m)^+ \mathbf{V}_n^{-1}(\mathbf{r}_{s}^t- \beta \mathbf{r}_{s}^{t+1})^{\tt{T}}\nonumber \\
       &+ E\tilde{V}_{t+1}^g(\alpha(\tilde{\mathbf{x}}-\tilde{x}_m\wedge \tilde{y}_n\mathbf{1}_m)^++\tilde{\mathbf{D}}^{t+1},\beta (\tilde{\mathbf{y}}-\tilde{x}_m\wedge \tilde{y}_n\mathbf{1}_n)^++\tilde{\mathbf{S}}^{t+1}).\label{eqn:V-g}
\end{align}

If $\tilde{x}_i\le \tilde{y}_n$, then $\tilde{x}_i\le \tilde{x}_m\wedge \tilde{y}_n$ and thus $(\tilde{x}_i- \tilde{x}_m\wedge \tilde{y}_n)^+= 0$.
If $\tilde{x}_i>\tilde{y}_n$, then $\tilde{x}_m\ge \tilde{x}_i>\tilde{y}_n$ and thus $(\tilde{x}_i- \tilde{x}_m\wedge \tilde{y}_n)^+=(\tilde{x}_i-\tilde{y}_n)^+=\tilde{x}_i-\tilde{y}_n$.
It follows that $(\tilde{x}_i- \tilde{x}_m\wedge \tilde{y}_n)^+ = (\tilde{x}_i-\tilde{y}_n)^+$ and therefore $(\tilde{\mathbf{x}}-\tilde{x}_m\wedge \tilde{y}_n\mathbf{1}_m)^+=(\tilde{\mathbf{x}}-\tilde{y}_n\mathbf{1}_m)^+$.
Likewise, we can show that $(\tilde{\mathbf{y}}-\tilde{x}_m\wedge \tilde{y}_n\mathbf{1}_n)^+=(\tilde{\mathbf{y}}-\tilde{x}_m\mathbf{1}_n)^+$.

Equation (\ref{eqn:V-tilde-g}) then follows from (\ref{eqn:V-g}). \Halmos
\endproof

The next lemma shows the monotonicity of the function $\tilde{V}_{t}^g(\tilde{\mathbf{x}},\tilde{\mathbf{y}})$.

\begin{lemma}
    Suppose that $r_{id}^t - r_{i+1,d}^t \ge \alpha(r_{id}^{t+1}-r_{i+1,d}^{t+1})$ for all $t=1,\ldots,T-1$.
    Then, for any period $t=1,\ldots,T$, the function $\tilde{V}_t^g(\tilde{\mathbf{x}},\tilde{\mathbf{y}})$ is decreasing in $x_i$ for all $i=1,\ldots,m$ and in $y_j$ for all $j=1\ldots,n$.
    \label{lem:Vg-decrease}
\end{lemma}
\proof{Proof of Lemma \ref{lem:Vg-decrease}.}
We will show by induction that $\tilde{V}_t^g(\tilde{\mathbf{x}}+\varepsilon\mathbf{e}_i^m,\tilde{\mathbf{y}})$ decreases in $\tilde{x}_i$.

It is trivial to prove for $t=T+1$, given that $\tilde{V}_t^g(\tilde{\mathbf{x}},\tilde{\mathbf{y}})\equiv -\tilde{\mathbf{x}}\mathbf{U}_m^{-1}(\mathbf{r}_d^t)^{\tt{T}} - \tilde{\mathbf{y}}\mathbf{U}_n^{-1}(\mathbf{r}_s^t)^{\tt{T}}$.
Let us suppose that $\tilde{V}_{t+1}^g(\tilde{\mathbf{x}},\tilde{\mathbf{y}})$ is decreasing in $\tilde{x}_i$.

To show that $\tilde{V}_{t}^g(\tilde{\mathbf{x}},\tilde{\mathbf{y}})$ is decreasing in $\tilde{x}_i$,
we note that 
$$-\sum_{i=1}^{m}(\tilde{\mathbf{x}}-\tilde{y}_n\mathbf{1}_m)^+\mathbf{U}_m^{-1} (\mathbf{r}_{d}^t-\alpha \mathbf{r}_{d}^{t+1})^{\tt{T}}=(\tilde{x}_i-\tilde{y}_n)^+[(r_{id}^t - r_{i+1,d}^t) - \alpha(r_{id}^{t+1}-r_{i+1,d}^{t+1})]$$
is decreasing in $\tilde{x}_i$.
According to the induction hypothesis, the last term in (\ref{eqn:V-tilde-g}),
$E\tilde{V}_{t+1}^g(\alpha(\tilde{\mathbf{x}}-\tilde{y}_n\mathbf{1}_m)^+\mathbf{U}_m^{-1}+\tilde{\mathbf{D}}^{t+1}\mathbf{U}_m^{-1},\beta (\tilde{\mathbf{y}}-\tilde{x}_m\mathbf{1}_n)^+\mathbf{V}_n^{-1}+\tilde{\mathbf{S}}^{t+1}\mathbf{V}_n^{-1})$ is decreasing in $\tilde{x}_i$.
Thus, all terms in (\ref{eqn:V-tilde-g}) are either constant or decreasing in $\tilde{x}_i$.
This completes the induction and shows that $\tilde{V}_t^g(\tilde{\mathbf{x}},\tilde{\mathbf{y}})$ is decreasing in $\tilde{x}_i$. 
We show that it is also decreasing in $\tilde{y}_j$ similarly.
\Halmos
\endproof

Next, we show that $\tilde{V}_{t}^g(\tilde{\mathbf{x}},\tilde{\mathbf{y}})$ is concave.

\begin{lemma}
    Suppose that $r_{id}^t - r_{i+1,d}^t \ge \alpha(r_{id}^{t+1}-r_{i+1,d}^{t+1})$ for all $t=1,\ldots,T-1$.
    The function $\tilde{V}_{t}^g(\tilde{\mathbf{x}},\tilde{\mathbf{y}})$ is concave in $(\tilde{\mathbf{x}},\tilde{\mathbf{y}})$.
    \label{lem:Vg-concave}
\end{lemma}
\proof{Proof of Lemma \ref{lem:Vg-concave}.}
The proof is again inductive. It is easy to see that $\tilde{V}^g_{T+1}$ is concave (it is actually linear).
Suppose that $\tilde{V}^g_{t+1}$ is concave.

 The terms $- (\tilde{\mathbf{x}}-\tilde{y}_n\mathbf{1}_m)^+\mathbf{U}_m^{-1} (\mathbf{r}_{d}^t-\alpha \mathbf{r}_{d}^{t+1})^{\tt{T}}$ and
    $- (\tilde{\mathbf{y}}-\tilde{x}_m\mathbf{1}_n)^+ \mathbf{V}_n^{-1}(\mathbf{r}_{s}^t- \beta \mathbf{r}_{s}^{t+1})^{\tt{T}}$ are concave due to the concavity of the function $f(x):=-x^+$.
    It remains to show that $E\tilde{V}_{t+1}^g(\alpha(\tilde{\mathbf{x}}-\tilde{y}_n\mathbf{1}_m)^++\tilde{\mathbf{D}}^{t+1},\beta (\tilde{\mathbf{y}}-\tilde{x}_m\mathbf{1}_n)^++\tilde{\mathbf{S}}^{t+1})$ is concave.
For $\lambda_1\ge 0$ and $\lambda_2\ge 0$ such that $\lambda_1+\lambda_2=1$, we have
\begin{align*}
    (\lambda_1\tilde{\mathbf{x}}+\lambda_2\tilde{\mathbf{x}}'-\lambda_1\tilde{y}_n\mathbf{1}_m -\lambda_2 \tilde{y}_n'\mathbf{1}_m)^+
    \le \lambda_1 (\tilde{\mathbf{x}}-\tilde{y}_n\mathbf{1}_m)^+ +\lambda_2 (\tilde{\mathbf{x}}'-\tilde{y}_n'\mathbf{1}_m)^+,
\end{align*}
and
\begin{align*}
    (\lambda_1\tilde{\mathbf{y}}+\lambda_2\tilde{\mathbf{y}}'-\lambda_1\tilde{x}_m\mathbf{1}_n -\lambda_2 \tilde{x}_m'\mathbf{1}_n)^+
    \le \lambda_1 (\tilde{\mathbf{y}}-\tilde{x}_m\mathbf{1}_n)^+ +\lambda_2 (\tilde{\mathbf{y}}'-\tilde{x}_m'\mathbf{1}_n)^+,
\end{align*}
where both inequalities follow from the convexity of the function $g(x):=x^+$.
Since $\tilde{V}_{t+1}^g$ is decreasing in its arguments (Lemma \ref{lem:Vg-decrease}), we have
\begin{align*}
    & \tilde{V}_{t+1}^g(\alpha(\lambda_1\tilde{\mathbf{x}}+\lambda_2\tilde{\mathbf{x}}'-\lambda_1\tilde{y}_n\mathbf{1}_m-\lambda_2\tilde{y}_n'\mathbf{1}_m)^++\tilde{\mathbf{D}}^{t+1},\beta (\lambda_1\tilde{\mathbf{y}}+\lambda_2\tilde{\mathbf{y}}'-\lambda_1\tilde{x}_m\mathbf{1}_n -\lambda_2\tilde{x}_m'\mathbf{1}_n)^++\tilde{\mathbf{S}}^{t+1})\\
    \ge &\tilde{V}_{t+1}^g(\lambda_1\cdot[\alpha(\tilde{\mathbf{x}}-\tilde{y}_n\mathbf{1}_m)^++\tilde{\mathbf{D}}^{t+1}]+\lambda_2\cdot[\alpha(\tilde{\mathbf{x}}'-\tilde{y}_n'\mathbf{1}_m)^++\tilde{\mathbf{D}}^{t+1}],\\
    &\quad\quad\quad\quad \lambda_1\cdot[\beta (\tilde{\mathbf{y}}-\tilde{x}_m\mathbf{1}_n)^++\tilde{\mathbf{S}}^{t+1}] + \lambda_2\cdot [\beta (\tilde{\mathbf{y}}'-\tilde{x}_m'\mathbf{1}_n)^++\tilde{\mathbf{S}}^{t+1}])\\
    \ge&\lambda_1 \tilde{V}_{t+1}^g(\alpha(\tilde{\mathbf{x}}-\tilde{y}_n\mathbf{1}_m)^++\tilde{\mathbf{D}}^{t+1},\beta (\tilde{\mathbf{y}}-\tilde{x}_m\mathbf{1}_n)^++\tilde{\mathbf{S}}^{t+1}) \\
    &+ \lambda_2 \tilde{V}_{t+1}^g(\alpha(\tilde{\mathbf{x}}'-\tilde{y}_n'\mathbf{1}_m)^++\tilde{\mathbf{D}}^{t+1},\beta (\tilde{\mathbf{y}}'-\tilde{x}_m'\mathbf{1}_n)^++\tilde{\mathbf{S}}^{t+1}),
\end{align*}
where the last inequality follows from the induction hypothesis of the concavity of $\tilde{V}^g_{t+1}$.
\Halmos
\endproof

The following proposition shows that the one-step-ahead policy has the top-down structure, i.e., in any period $t$, a lower-quality demand/supply type will not be consumed, unless all higher-quality types are fully used.

\begin{proposition}
    The one-step-ahead policy has the top-down structure in each period $t$.
    \label{prop:ons-top-down}
\end{proposition}
\proof{Proof of Proposition \ref{prop:ons-top-down}.}
Let $\mathbf{Q}$ be the matching decision in period $t$ under the one-step-ahead policy and $(\mathbf{u},\mathbf{v})$ be the post-matching levels. Since we consider linearly additive reward, it is sufficient to show that there is no ``blanks'' for the matching decision in any period $t$, i.e., for $i'>i$, $q_{i'j}>0$ would imply that $u_{i}=0$, and for $j'>j$, $q_{ij'}>0$ would imply that $v_j=0$. That is, due to the linearly additive reward structure, $\succ_{\mathcal{M}}$ is sufficient to ensure $\succ_{\mathcal{M}_s}$. 

We suppose to the contrary that both $q_{i'j}$ and $u_{i}$ are positive for some $i'>i$, in a period $t$ under the one-step-ahead policy.
We will construct a feasible decision that has the top-down structure and is weakly better than the current decision in period $t$.

To construct the new decision, in period $t$ we reduce the consumption of type $i'$ demand by $\varepsilon$ and increase that of type $i$ demand by $\varepsilon$, where $\varepsilon=\min\left\{ q_{ij'}, u_i \right\}$.
By doing so, either the consumption of type $i'$ demand becomes zero or type $i$ demand is fully used.
the matching reward received in period $t$ increases by $(r_{id}^t-r_{i'd}^t)\varepsilon$, and the post matching levels become $(\mathbf{u}-\varepsilon\mathbf{e}_i^m+\varepsilon\mathbf{e}_{i'}^m,\mathbf{v})$.
Under greedy matching, the change in the total expected reward from period $t+1$ to period $T$ is
\begin{align*}
    &EV^g_{t+1}(\alpha \mathbf{u} - \alpha\varepsilon\mathbf{e}_{i}^m+\alpha\varepsilon\mathbf{e}_{i'}^m  +\mathbf{D}^{t+1},\beta\mathbf{v} + \mathbf{S}^{t+1})
    - EV^g_{t+1}(\alpha \mathbf{u} +\mathbf{D}^{t+1},\beta\mathbf{v} + \mathbf{S}^{t+1})\\
    =& E\tilde{V}^g_{t+1}(\alpha \tilde{\mathbf{u}} - \alpha\varepsilon \sum_{i''=i}^{i'-1}\mathbf{e}_{i''}^m  +\tilde{\mathbf{D}}^{t+1},\beta\tilde{\mathbf{v}} + \tilde{\mathbf{S}}^{t+1})\\
   & + (\alpha \tilde{\mathbf{u}} - \alpha\varepsilon \sum_{i''=i}^{i'-1}\mathbf{e}_{i''}^m  +E\tilde{\mathbf{D}}^{t+1})\tilde{\mathbf{U}}_m^{-1}(\mathbf{r}_d^{t+1})^{\tt{T}}
     + (\beta\tilde{\mathbf{v}} + E\tilde{\mathbf{S}}^{t+1}) \mathbf{V}_n^{-1} (\mathbf{r}_s^{t+1})^{\tt{T}}\\
     &- E\tilde{V}^g_{t+1}(\alpha \tilde{\mathbf{u}}   +\tilde{\mathbf{D}}^{t+1},\beta\tilde{\mathbf{v}} + \tilde{\mathbf{S}}^{t+1})
      - (\alpha \tilde{\mathbf{u}} +E\tilde{\mathbf{D}}^{t+1})\tilde{\mathbf{U}}_m^{-1}(\mathbf{r}_d^{t+1})^{\tt{T}}
     - (\beta\tilde{\mathbf{v}} + E\tilde{\mathbf{S}}^{t+1}) \mathbf{V}_n^{-1} (\mathbf{r}_s^{t+1})^{\tt{T}}\\
     =& E\tilde{V}^g_{t+1}(\alpha \tilde{\mathbf{u}} - \alpha\varepsilon \sum_{i''=i}^{i'-1}\mathbf{e}_{i''}^m  +\tilde{\mathbf{D}}^{t+1},\beta\tilde{\mathbf{v}} + \tilde{\mathbf{S}}^{t+1})
     -E\tilde{V}^g_{t+1}(\alpha \tilde{\mathbf{u}}   +\tilde{\mathbf{D}}^{t+1},\beta\tilde{\mathbf{v}} + \tilde{\mathbf{S}}^{t+1})
      - \alpha\varepsilon (\sum_{i''=i}^{i'-1}\mathbf{e}_{i''}^m)\tilde{\mathbf{U}}_m^{-1}(\mathbf{r}_d^{t+1})^{\tt{T}}\\
      \ge& - \alpha\varepsilon (\sum_{i''=i}^{i'-1}\mathbf{e}_{i''}^m)\tilde{\mathbf{U}}_m^{-1}(\mathbf{r}_d^{t+1})^{\tt{T}}\\
      =& -\alpha\varepsilon (r_{id}^{t+1}-r_{i'd}^{t+1}).
\end{align*}

Thus, the change in the total reward from period $t$ to period $T$ is no less than $(r_{id}^t-r_{i'd}^t)\varepsilon-\alpha\varepsilon (r_{id}^{t+1}-r_{i'd}^{t+1})=[(r_{id}^t-r_{i'd}^t)-\alpha (r_{id}^{t+1}-r_{i'd}^{t+1})]\varepsilon\ge 0$.
This implies that the total expected reward increases by reducing the consumption of type $i'$ demand by $\varepsilon$ and increasing that of type $i$ demand by $\varepsilon$.

Analogously, we can show that  total expected reward also increases if we reduce the consumption of type $j'$ supply and increase that of type $j<j'$ supply by the same amount.

We repeatedly transfer quantity from a lower-quality type to a higher-quality type in period $t$, and will eventually arrive at a decision that has the top-down structure.
\Halmos
\endproof

Proposition \ref{prop:ons-top-down} implies that under the one-step-ahead policy, the matching decision in a period $t$ is fully determined by the total matching quantity $Q$. 
Therefore, the one-step-ahead policy reduces to a one-dimension problem for choosing $Q$ to maximize the total expected reward from period $t$ to period $T$, provided that greedy matching is enforced starting from period $t+1$.
We conclude this appendix by presenting a formulation of the optimization problem associated with one-step-ahead policy and showing its concavity.

\begin{proposition}
    The optimal matching quantity $Q^*$ in period $t$ under the one-step-ahead policy solves the following problem.
    \begin{align}
        \max_{0\le Q\le \tilde{x}_m\wedge \tilde{y}_n} \quad F_t(Q,\tilde{\mathbf{x}},\tilde{\mathbf{y}}):=& \tilde{\mathbf{x}}\mathbf{U}_m^{-1} (\mathbf{r}_{d}^t)^{\tt{T}}
+ \tilde{\mathbf{y}}\mathbf{V}_n^{-1} (\mathbf{r}_{s}^t)^{\tt{T}} +E\tilde{\mathbf{D}}^{t+1}\mathbf{U}_m^{-1}+E\tilde{\mathbf{S}}^{t+1}\mathbf{V}_n^{-1}\nonumber\\
      & - (\tilde{\mathbf{x}}-Q\mathbf{1}_m)^+\mathbf{U}_m^{-1} (\mathbf{r}_{d}^t-\alpha \mathbf{r}_{d}^{t+1})^{\tt{T}}
       - (\tilde{\mathbf{y}}-Q\mathbf{1}_m)^+ \mathbf{V}_n^{-1}(\mathbf{r}_{s}^t- \beta \mathbf{r}_{s}^{t+1})^{\tt{T}}\nonumber \\
       &+ E\tilde{V}_{t+1}^g(\alpha(\tilde{\mathbf{x}}-Q\mathbf{1}_m)^+ +\tilde{\mathbf{D}}^{t+1},\beta (\tilde{\mathbf{y}}-Q\mathbf{1}_n)^++\tilde{\mathbf{S}}^{t+1}).\label{eqn:osa-obj}
   \end{align}
   The function $F_t(Q,\tilde{\mathbf{x}},\tilde{\mathbf{y}})$ is concave in $Q$, and represents the total expected reward to be received from period $t$ to period $T$ for using the total matching quantity $Q$ in period $t$ and greedy matching from period $t+1$ on.
   \label{prop:osa-formulation}
\end{proposition}
\proof{Proof of Proposition \ref{prop:osa-formulation}.}
In the proof of Lemma \ref{lem:Vg-recursion}, we have shown that $F_t$ defined in (\ref{eqn:osa-obj}) is the total expected reward to be received from period $t$ to period $T$ for matching a total quanitty $Q$ under the top-down matching structure in period $t$ and applying the greedy policy from period $t$ to period $T+1$.

Thus, it remains to show that $F_t$ is concave in $Q$. We now show that all terms on the right-hand-side of (\ref{eqn:osa-obj}) is concave in $Q$.
It is easy to see that the term $\tilde{\mathbf{x}}\mathbf{U}_m^{-1} (\mathbf{r}_{d}^t)^{\tt{T}}
+ \tilde{\mathbf{y}}\mathbf{V}_n^{-1} (\mathbf{r}_{s}^t)^{\tt{T}}$ is concave, due to its linearity.
The terms $- (\tilde{\mathbf{x}}-Q\mathbf{1}_m)^+\mathbf{U}_m^{-1} (\mathbf{r}_{d}^t-\alpha \mathbf{r}_{d}^{t+1})^{\tt{T}}$ and $ (\tilde{\mathbf{y}}-Q\mathbf{1}_m)^+ \mathbf{V}_n^{-1}(\mathbf{r}_{s}^t- \beta \mathbf{r}_{s}^{t+1})^{\tt{T}}$ are concave because the function $g(x):=-x^+$ is concave.

It remains to prove that $E\tilde{V}_{t+1}^g(\alpha(\tilde{\mathbf{x}}-Q\mathbf{1}_m)^+ +\tilde{\mathbf{D}}^{t+1},\beta (\tilde{\mathbf{y}}-Q\mathbf{1}_n)^++\tilde{\mathbf{S}}^{t+1})$ is concave.
To that end, let $\lambda_1\ge 0$ and $\lambda_2\ge 0$ such that $\lambda_1+\lambda_2=1$.
For two quantities $Q$ and $Q'$, we have
\begin{align*}
    &\tilde{V}_{t+1}^g(\alpha(\tilde{\mathbf{x}}-(\lambda_1 Q + \lambda_2 Q')\mathbf{1}_m)^+ +\tilde{\mathbf{D}}^{t+1},\beta (\tilde{\mathbf{y}}-(\lambda_1 Q + \lambda_2 Q')\mathbf{1}_n)^++\tilde{\mathbf{S}}^{t+1})\\
    =& \tilde{V}_{t+1}^g(\alpha[\lambda_1(\tilde{\mathbf{x}}- Q \mathbf{1}_m) + \lambda_2 (\tilde{\mathbf{x}}- Q' \mathbf{1}_m)]^+ +\tilde{\mathbf{D}}^{t+1}, 
     \beta [\lambda_1(\tilde{\mathbf{y}}- Q \mathbf{1}_n)+\lambda_2 (\tilde{\mathbf{y}}- Q' \mathbf{1}_n)]^++\tilde{\mathbf{S}}^{t+1}  )\\
    \ge &\tilde{V}_{t+1}^g(\lambda_1[\alpha(\tilde{\mathbf{x}}- Q \mathbf{1}_m)^+ +\tilde{\mathbf{D}}^{t+1}] + \lambda_2[\alpha(\tilde{\mathbf{x}}- Q' \mathbf{1}_m)^+ +\tilde{\mathbf{D}}^{t+1}], \\
    &\quad\quad\quad\quad \lambda_1 [\beta (\tilde{\mathbf{y}}- Q \mathbf{1}_n)^++\tilde{\mathbf{S}}^{t+1} ] + \lambda_2 [\beta (\tilde{\mathbf{y}}- Q' \mathbf{1}_n)^++\tilde{\mathbf{S}}^{t+1} ])\\
    \ge& \lambda_1 \tilde{V}_{t+1}^g(\alpha(\tilde{\mathbf{x}}-Q\mathbf{1}_m)^+ +\tilde{\mathbf{D}}^{t+1},\beta (\tilde{\mathbf{y}}-Q\mathbf{1}_n)^++\tilde{\mathbf{S}}^{t+1})
    + \lambda_2 \tilde{V}_{t+1}^g(\alpha(\tilde{\mathbf{x}}-Q'\mathbf{1}_m)^+ +\tilde{\mathbf{D}}^{t+1},\beta (\tilde{\mathbf{y}}-Q'\mathbf{1}_n)^++\tilde{\mathbf{S}}^{t+1}),
\end{align*}
where the first inequality follows from the convexity of the function $g(x):=x^+$ and the fact that $\tilde{V}_{t+1}$ is decreasing in its arguments (Lemma \ref{lem:Vg-decrease}),
and the second inequality follows from the concavity of $\tilde{V}^g_{t+1}$. \Halmos
\endproof

\newpage
 \setcounter{page}{1}
\renewcommand{\theequation}{\arabic{equation}}
\setcounter{equation}{0}
\setcounter{section}{0}
\setcounter{lemma}{0}
\setcounter{example}{0}
\setcounter{theorem}{0}
\setcounter{corollary}{0}
\setcounter{proposition}{0}

\end{document}